\documentclass[12pt]{amsart}

\usepackage{amssymb,amsfonts,amsmath,amsopn,amstext,amscd,latexsym,fullpage,verbatim}
\usepackage{amsthm}

\newtheorem{Proposition}{Proposition}[subsection]
\newtheorem{Theorem}[Proposition]{Theorem}
\newtheorem{Corollary}[Proposition]{Corollary}
\newtheorem{Lemma}[Proposition]{Lemma}
\newtheorem{Example}[Proposition]{Example}
\newtheorem{Remark}[Proposition]{Remark}

\numberwithin{equation}{section}

\newcommand{\bbG}{{\mathbb G}}

\newcommand{\bbN}{{\mathbb N}}

\newcommand{\bbQ}{{\mathbb Q}}

\newcommand{\bbZ}{{\mathbb Z}}

\newcommand{\bA}{{\bf A}}

\newcommand{\fra}{{\mathfrak a}}

\newcommand{\frd}{{\mathfrak d}}

\newcommand{\frg}{{\mathfrak g}}

\newcommand{\frl}{{\mathfrak l}}

\newcommand{\frp}{{\mathfrak p}}

\newcommand{\frt}{{\mathfrak t}}
\newcommand{\fru}{{\mathfrak u}}

\newcommand{\frx}{{\mathfrak x}}

\newcommand{\frz}{{\mathfrak z}}

\newcommand{\cC}{{\mathcal C}}

\newcommand{\cF}{{\mathcal F}}
\newcommand{\cG}{{\mathcal G}}
\newcommand{\cH}{{\mathcal H}}
\newcommand{\cI}{{\mathcal I}}

\newcommand{\cO}{{\mathcal O}}

\newcommand{\cX}{{\mathcal X}}
\newcommand{\cY}{{\mathcal Y}}

\newcommand{\Qp}{\mathbb Q_p} 
\newcommand{\Cp}{\mathbb C_p} 

\newcommand{\Lpn}{\Lambda/\pi^n\Lambda}
\newcommand{\Ln}{\Lambda^{(n)}}

\newcommand{\N}{\mathbb N}
\newcommand{\Z}{\mathbb Z}
\newcommand{\R}{\mathbb R}
\renewcommand{\P}{\mathbb P}

\newcommand {\Abb}[5]{\begin{eqnarray*} #1  #2  & \longrightarrow & #3 \\
              #4 & \longmapsto & #5 \end{eqnarray*} }
\newcommand {\aequi}{\mbox{$\Leftrightarrow$}}
\newcommand{\dej}{\mathfrak \Delta_{\overline{j}}}
\newcommand{\deju}{\mathfrak \Delta_{\underline{j}}}

\begin{document}

\title{Equivariant vector bundles on Drinfeld's upper half space}
\author{Sascha Orlik}
\address{Mathematisches Institut, Universit\"at Leipzig, Postfach 10 09 20,
D-04009 Leipzig,  Germany}
\email{orlik@math.uni-leipzig.de}
\date{}
\maketitle

\begin{abstract}
Let $\cX \subset \P_K^d$ be Drinfeld's upper half space over a finite extension $K$ of $\Qp.$ We construct for every ${\rm GL}_{d+1}$-equivariant vector bundle $\cF$ on $\P^d_K,$ a ${\rm GL}_{d+1}(K)$-equivariant filtration by closed subspaces on the $K$-Fr\'echet $H^0(\cX,\cF).$ This gives rise by duality to a filtration by locally analytic ${\rm GL}_{d+1}(K)$-representations
on the strong dual $H^0(\cX,\cF)'.$ The graded pieces of this filtration are locally analytic induced representations from locally algebraic ones with respect to maximal parabolic subgroups. This paper generalizes the cases of the canonical bundle due to Schneider and Teitelbaum \cite{ST1} and
that of the structure sheaf by Pohlkamp \cite{P}.
\end{abstract}

\normalsize

\tableofcontents

\section*{Introduction}
Let $K$ be a finite extension of $\Qp$.
We denote by $\cX$ Drinfeld's upper half space of dimension $d\geq1$ over $K$. Explicitly, $\cX$ is the complement of all  $K$-rational hyperplanes in projective space $\mathbb \P_K^d,$ i.e.,
$$ \cX=\mathbb \P_K^d\setminus \bigcup\nolimits_{H \varsubsetneq K^{d+1}}\mathbb \P(H).$$
The interest for studying the rigid analytic variety $\cX$ is manifold. The first is its connection to formal groups, since
it is the generic fibre of a formal scheme classifying certain $p$-divisible groups due  to Drinfeld \cite{D}. Another aspect is its role for the uniformization of certain Shimura varieties \cite{RZ}. On the other hand, as conjectured by Drinfeld, the \'etale coverings of $\cX$ realize the supercuspidal spectrum of the local Langlands correspondence by considering the $\ell$-adic cohomology of these spaces. Our interest is connected with the latter aspect. In \cite{S1} Schneider studies the cohomology of local systems on projective varieties which are uniformized by $\cX.$ For this purpose, he defines the notion of a $p$-adic holomorphic discrete series representation.
These representations can be realized by the space of rigid analytic holomorphic sections $\cF(\cX)$
of  ${\rm GL}_{d+1}$-equivariant vector bundles $\cF$ on $\P^d_K.$  They appear naturally in the theory of $p$-adic modular forms \cite{T1}, \cite{T2}. The dual of a holomorphic discrete series representation is a locally analytic ${\rm GL}_{d+1}(K)$-representation in the sense of  Schneider and Teitelbaum  \cite{ST3}. Those representations come up in the $p$-adic Langlands theory of Breuil and Schneider \cite{BS} as the locally analytic part of certain Banach space representations. Our aim in this paper is a description of $p$-adic holomorphic discrete series representations by considering a ${\rm GL}_{d+1}(K)$-equivariant filtration on them. We determine the graded pieces in terms of locally analytic ${\rm GL}_{d+1}(K)$-representations.

In \cite{SS} Schneider and Stuhler computed the \'etale and the de Rham cohomology of Drinfeld's upper half space.
The cohomology is equipped with an action of the general linear group $G={\rm GL}_{d+1}(K)$, which is induced by the natural action on $\cX.$ It turns out that
the  cohomology groups are duals
of certain smooth elliptic representations having Iwahori fix vectors. In particular, the top cohomology group is the dual of the Steinberg representation.
In contrast, the  space of holomorphic sections $\Omega^d(\cX)=H^0(\cX, \Omega^d)$ of the canonical bundle $\Omega^d$ on $\P^d_K$  is a much bigger object, it is a reflexive $K$-Fr\'echet space with a continuous $G$-action. Its strong dual $\Omega^d(\cX)'$ is a locally analytic $G$-representation.
In order to describe this latter space,
 Schneider and Teitelbaum construct in \cite{ST1} a $G$-equivariant decreasing filtration by closed $K$-Fr\'echet spaces
 $$\Omega^d(\cX)^0 \supset \Omega^d(\cX)^1 \supset \cdots \supset \Omega^d(\cX)^{d-1} \supset \Omega^d(\cX)^d \supset \Omega^d(\cX)^{d+1}=\{0\}$$ on $\Omega^d(\cX)^0=\Omega^d(\cX).$
 The definition of the filtration involves the geometry of $\cX$ being the complement of an hyperplane arrangement.
Further they construct  isomorphisms
$$I^{[j]} : (\Omega^d(\cX)^j/\Omega^d(\cX)^{j+1})' \stackrel{\thicksim}{\longrightarrow} C^{an}(G,P_{\underline{j}}; V_j')^{\deju =0}$$
of locally analytic $G$-representations, which they call (partial) boundary value maps.
Here, $P_{\underline{j}}$=$P_{(j,d+1-j)}\subset G$ is the (lower) standard-parabolic  subgroup attached to the decomposition $(j,d+1-j)$
of $d+1.$ The right hand side is a locally analytic induced representation. The  $P_{\underline{j}}$-representation
$V_j'$ is a locally algebraic representation. It is isomorphic to the tensor product ${\rm Sym}^j(K^{d+1-j}) \otimes {\rm St}_{d+1-j}$ of the irreducible algebraic ${\rm GL}_{d+1-j}$-representation ${\rm Sym}^j(K^{d+1-j})$ and the Steinberg representation ${\rm St}_{d+1-j}$ of ${\rm GL}_{d+1-j}(K).$
Here the factor ${\rm GL}_{d+1-j}(K)$ of the  Levi  subgroup $L_{(j,d+1-j)}= {\rm GL}_j(K) \times {\rm GL}_{d+1-j}(K)$ acts through the way just described. The action of ${\rm GL}_j(K)$ is given by the inverse of the determinant character. The operation of the unipotent radical of $P_{\underline{j}}$ on $V_j'$ is trivial.
Finally, $\deju$ denotes a system of differential equations which is here a submodule of a generalized Verma module.
In particular, the case $j=0$, i.e., the first subquotient of the above filtration  is isomorphic to $H^d_{dR}(\cX)$ and  yields the Steinberg representation of $G.$
Their paper presents consequently in a sense a generalization of the computation \cite{SS} since
it computes not only the top cohomology of $\cX.$
Further, it generalizes  pioneering work by Morita (e.g. \cite{Mo2}) who considered such representations in the ${\rm SL}_2$-case. We refer to the introduction of \cite{ST1} for a more comprehensive background on this topic.

Pohlkamp \cite{P} considers the other extreme, that of the structure sheaf $\Omega^0= \cO$ on $\P^d_K$. Again, by using a similar construction, Pohlkamp defines a
$G$-equivariant increasing filtration by closed $K$-Fr\'echet spaces
 $$K=\cO(\cX)_0 \subset \cO(\cX)_1 \subset \cdots \subset\cO(\cX)_{d-1} \subset \cO(\cX)_d$$  on $\cO(\cX)_d=H^0(\cX,\cO)$ together with isomorphisms
 $$(\cO(\cX)_j)/\cO(\cX)_{j-1})' \stackrel{\thicksim}{\longrightarrow} C^{an}(G,P_{\underline{d+1-j}}; W_j')^{\dej=0}$$
 of locally analytic $G$-representations.
Analogously to the  above case,  the $P_{\underline{d+1-j}}$-representation $W_j'$ is a tensor product of a Steinberg representation and an irreducible algebraic representation.

Our goal in this paper is to construct a decreasing $G$-equivariant filtration on $\cF(\cX)$ for all $G$-equivariant vector bundles  $\cF$  on $\cX,$ which are induced  by restriction of a homogeneous vector bundle on $\P^d_K\cong {G}/ {P_{(1,d)}}.$ The latter objects are defined by finite-dimensional algebraic representations of the parabolic subgroup ${P_{(1,d)}}.$
Our approach is different from \cite{ST1}, \cite{P}. We use local cohomology of coherent sheaves on rigid analytic varieties
as a technical ingredient.
In fact, $\cF(\cX)=H^0(\cX,\cF)$ appears in an exact sequence
$$0\rightarrow H^0(\P_K^d,\cF) \rightarrow H^0(\cX,\cF)\rightarrow H^1_{\cY}(\P_K^d,\cF) \rightarrow H^1(\P_K^d,\cF)\rightarrow 0.$$
We consider the $K$-Fr\'echet space  $H^1_{\cY}(\P_K^d,\cF),$ where ${\cY}\subset \P^d_K$ is the "closed" complement of $\cX$ in $\P_K^d.$
By a technique used in \cite{O3}, we are able to compute this latter module as a $G$-representation. Here, we use an acyclic resolution
of the constant sheaf $\Z$ on ${\cY}^{ad},$ where  $\cY^{ad} \stackrel{i}{\hookrightarrow} (\P^d_K)^{ad}$  is the closed complement of the adic space $\cX^{ad}$ in $(\P^d_K)^{ad}.$  By applying the functor ${\rm Hom}(i_\ast(\;-\;),\cF)$ to this complex, we get a spectral sequence converging to $H_{{\cY}^{ad}}^1((\P_K^d)^{ad},\cF^{ad})=H_{{\cY}}^1(\P_K^d,\cF).$
The canonical filtration on $H^1_{\cY}(\P_K^d,\cF)$ coming from this spectral sequence gives rise to  a decreasing filtration by closed $K$-Fr\'echet spaces
$$\cF(\cX)^0 \supset \cF(\cX)^{1} \supset \cdots \supset \cF(\cX)^{d-1} \supset \cF(\cX)^d =  H^0(\P^d,\cF)$$
on $\cF(\cX)^0=H^0(\cX,\cF).$ Our first  main theorem  is:

\medskip
{\it \noindent {\bf Theorem 1:} Let $\cF$ be a homogeneous vector bundle on $\P^d_K.$ For $j=0,\ldots,d-1,$ there are
extensions of locally analytic $G$-representa\-tions
$$0 \rightarrow v^{G}_{{P_{(j+1,1,\ldots,1)}}}(H^{d-j}(\P^d_K,\cF)') \rightarrow (\cF(\cX)^{j}/\cF(X)^{j+1})' \rightarrow  C^{an}(G,P_{\underline{j+1}};U_j')^{{\frd_j}=0} \rightarrow 0. $$ }

\noindent Here the module $v^{G}_{{P_{(j+1,1,\ldots,1)}}}(H^{d-j}(\P^d_K,\cF)')$ is a generalized Steinberg representation with coefficients in the finite-dimensional  algebraic $G$-module $H^{d-j}(\P^d_K,\cF)'.$
The $P_{\underline{j+1}}$-representation $U_j'$ is a tensor product $N_{j}' \otimes {\rm St}_{d-j}$ of an algebraic $P_{\underline{j+1}}$-representa\-tion $N_{j}'$ and the Steinberg representa\-tion ${\rm St}_{d-j}.$
The symbol $\frd_j$  indicates again  a system of differential equations depending on $N_j.$
Indeed, the representation $N_j$ is not uniquely determined. It is characterized  by the property that it generates the kernel of the natural  homomorphism $H^{d-j}_{\P^j_K} (\P^d_K, \cF) \rightarrow  H^{d-j}(\P^d_K,\cF)$ 
as a module with respect to the universal enveloping algebra $U(\frg)$ of the Lie algebra of $G.$
By enlarging the module $N_j$ we have to enlarge $\frd_j$, as well. In either case, the locally analytic $G$-representations  $C^{an}(G,P_{\underline{j+1}};U_j')^{{\frd_j}}$ remains the same.

In the case where $\cF$  arises from an irreducible representation of the Levi subgroup $L_{(1,d)},$  we can make our result more precise.
Let  $\lambda '=(\lambda_1,\ldots,\lambda_d) \in \Z^d$ be a dominant integral weight  of ${\rm GL}_d$  and a let $\lambda_0\in \Z.$ Set $\lambda:=(\lambda_0,\lambda_1,\ldots,\lambda_d)\in \Z^{d+1}.$ Denote by $\cF_\lambda$ the homogeneous vector bundle on $\P^d_K$ such that its fibre in the base point is the irreducible algebraic $L_{(1,d)}$-representation corresponding to  $\lambda.$ Put $w_j:=s_j\cdots s_1,$ where $s_i\in W$ is the (standard) simple reflection in the Weyl group $W\cong S_{d+1}$ of $G$. Then the above representation $N_j$  can be characterized as follows. By Bott \cite{Bo}  we know that there is at most one integer $i\geq 0$ with  $H^{i}(\P^d_K,\cF)\neq 0.$  Denote this integer by $i_0$ if it exists. Otherwise, there is an $i_0\leq d-1$ with
 $w_{i_0}\ast \lambda = w_{i_0+1}\ast \lambda,$  where $\ast$ is the dot operator of $W$ on the set of weights.
For $j=1,\ldots,d$, we set
$$\mu_{j,\lambda}:=\left\{ \begin{array}{cc} w_{j-1}\ast \lambda  &: j \leq i_0 \\  w_{j} \ast \lambda & : j > i_0 \end{array} \right..$$ 
Write
$\mu_{j,\lambda}=(\mu',\mu'')$  with $\mu'\in \bbZ^{j}$ and $\mu''\in \bbZ^{d-j+1}.$
For $j=1,\ldots,d$, let
\begin{eqnarray*}
\Psi_{j,\lambda} =  \bigcup\limits_{k=0}^{|\mu''|} &\Big\{&   \big(\mu''+(c_1,\ldots, c_{d-j+1}),\mu'-(d_{j},\ldots,d_1)\big) \mid {\textstyle \sum}_l c_l={\textstyle \sum}_l d_{l}=k, c_1=0 \\ & & \mbox{ or } d_1=0,\;  c_{l+1}\leq \mu_l'' -\mu_{l+1}'',  
\;l=1,\ldots,d-j,\; d_{l+1} \leq \mu'_{j-l} - \mu'_{j-l+1}, \,\\ & &  \; l=1,\ldots,j-1 \, \Big\} .
\end{eqnarray*}
Here $|\mu''|=\mu''_1 -\mu''_{d-j+1}.$ The elements in the finite set $\Psi_{j,\lambda}$ are dominant with respect to the Levi subgroup $L_{(d-j+1,j)}$ and $(\mu'',\mu')$ is its highest weight. Hence, for $\mu\in \Psi_{j,\lambda}$, we may consider the irreducible algebraic  $L_{(d-j+1,j)}$-representation $V_\mu$ attached to it.   

{\it \noindent {\bf Theorem 2:} Let $\cF=\cF_\lambda$ be the homogeneous
vector bundle on $\P^d_K$ with respect to the  dominant integral
weight $\lambda \in \Z^{d+1}$ of $L_{(1,d)}.$ Then we can choose
$N_j$ to be a quotient of $\bigoplus_{\mu\in \Psi_{d-j,\lambda} }V_{\mu}.$  }

\noindent We want to point out that for some weights $\lambda$, it happens that all irreducible constituents in the direct sum  apart from the module $V_\mu$ with $\mu=(\mu'',\mu')$ vanish under  the corresponding quotient map. But also the other extreme is possible, i.e., the quotient map can be an isomorphism, as well. 

Our filtration coincides  with the filtrations in \cite{ST1},
\cite{P}. More precisely,  in the case of the structure sheaf
$\cF=\cO$ the increasing filtration of Pohlkamp is related to our
decreasing filtration by $\cF(\cX)^j = \cO(\cX)_{d-j},
\;j=0,\ldots,d.$ In the case $\cF=\Omega^d$ the filtration of
Schneider and Teitelbaum is related to ours by a shift, i.e., we
have $\cF(\cX)^{i} = \Omega^d(\cX)^{i+1}$ for $i\geq 1.$ For $i=0,$ we
get an extension
$$0\rightarrow  \Omega^{d}(\cX)^1/\Omega^{d}(\cX)^2 \rightarrow \cF(\cX)^{0}/\cF(\cX)^{1} \rightarrow \Omega^{d}(\cX)^0/\Omega^{d}(\cX)^1 \rightarrow 0.$$
The dual sequence coincides with the corresponding  one  of Theorem 1.

The content of this paper is organized as follows. The first part  deals with  algebraic and analytic local cohomology of equivariant vector bundles $\cF$
on $\P^d_K.$ In the first section we
recall some facts on the restriction of these  bundles to $\cX$. Amongst other things, we explain the structure of
the strong dual $H^0(\cX,\cF)'$ as a locally analytic
$G$-representation.  In the following section we treat the algebraic local cohomology of $G$-equivariant vector bundles on $\P^d_K$.
We study the cohomology groups $H^{d-j}_{\P^j_K}(\P^d_K,\cF)$
as representations of $U(\frg)$ and $P_{\underline{j+1}}$. In Section 1.3 we turn to the analytic local cohomology groups  $H^{d-j}_{\P^j_K(\epsilon)}(\P^d_K,\cF)$ with support in the rigid analytic tube $\P^j_K(\epsilon).$ These groups
are naturally  equipped with a topology of a locally convex $K$-vector space. One of our focal points is to see that they are Hausdorff. Further, we prove a local duality theorem which is similar to
a result obtained by Morita \cite{Mo3}. It describes the dual of {\rm ker}$(H^{d-j}_{\P^j_K}(\P^d_K,\cF) \rightarrow H^{d-j}(\P^d_K,\cF) )$ for $\epsilon \to 0$, by means of analytic functions on certain polydiscs. In the final section of the first part we compute this kernel as representation of $U(\frg)$ and $P_{\underline{j+1}}$ 
when $\cF=\cF_\lambda$ is defined by a dominant integral weight $\lambda \in \Z^{d+1}$ of $L_{(1,d)}.$ Here we make use of the Grothendieck-Cousin complex with respect to the covering by Schubert cells. The second part of this paper deals with the computation of $H^0(\cX,\cF)$ as $G$-representation. 
First we repeat the construction of an acyclic resolution of the constant
sheaf $\Z$ on the closed complement $\cY^{ad}$ \cite{O3}.  In Section 2.2 we
evaluate the spectral sequence obtained by applying the functor
${\rm Hom}(i_\ast(\;-\;),\cF)$ to this acyclic complex. In Part 3
we compare our result to that of \cite{ST1} and \cite{P}.
Further, we provide with the cotangent
bundle $\Omega^1$ another example  for our computation. Finally, in the  Appendix we
present an alternative way for the computation avoiding adic spaces.
It is based purely on rigid analytic  varieties.

\vspace{0.7cm}
{\it Notation:} We denote  by $p$ a prime,  by $K\supset \bbQ_p$ a finite extension of the field of $p$-adic integers $\Qp$, by  $O_K$ its ring of integers and by  $\pi$ a uniformizer of $K$. Let $|\; \;|: K \rightarrow \R$ be the normalized norm, i.e., $|\pi|= \#(O_K/(\pi))^{-1}.$
We denote by $\Cp$ the completion of an algebraic closure $\overline{K}$ of $K$. Let $S:=K[X_0,\ldots, X_d]$ be the polynomial ring in $d+1$ indeterminates and denote by $\P^d_K:={\rm Proj}(S)$ the projective space over $K.$ If $Y \subset \P^d_K$ is a closed algebraic $K$-subvariety and $\cF$ is a sheaf on $\P^d_K$ we write $H_Y^\ast(\P^d_K,\cF)$ for the corresponding local cohomology. If $Y$ is a rigid analytic subvariety (resp. pseudo-adic subspace) of $(\P^d_K)^{rig}$ (resp.  $(\P^d_K)^{ad}$)  we also write  $H_Y^\ast(\P^d_K,\cF)$ instead of $H_Y^\ast((\P^d_K)^{rig},\cF^{rig})$ (resp. $H_Y^\ast((\P^d_K)^{ad},\cF^{ad})$)
to simplify matters.  For a locally convex $K$-vector space $V$, we denote by $V'$ its strong dual, i.e., the $K$-vector space of continuous linear forms
equipped with the strong topology of bounded convergence.

We use bold letters $\bf G, \bf P,\ldots $ to denote algebraic group schemes over $K$, whereas we use normal letters
$G,P, \ldots $ for their $K$-valued points of $p$-adic groups.
We use Gothic letters $\frg,\frp,\ldots $ for their Lie algebras. The corresponding enveloping algebras are denoted as usual by $U(\frg), U(\frp), \ldots
.$ Finally, we set ${\bf G}:= {\rm {\bf GL_{d+1}}}.$ If ${\bf H\subset G}$ is any closed linear algebraic subgroup and $R$ is a $O_K$-algebra, then we denote for simplicity by ${\bf H}(R)$ the
set of $R$-valued points of the schematic closure of ${\bf H}$ in ${\rm {\bf GL}_{{\bf d+1},O_K}}.$

\vspace{0.7cm}
{\it Acknowledgments:}
I am much obliged to M. Strauch for pointing out to me the topic treated in this paper and for all the
discussions we had during our stay at the IH\'ES in  2004.
I wish to thank M. Rapoport and P. Schneider for their numerous and interesting remarks on this paper.
Further, I would like to thank  Benjamin Schraen for indicating a mistake in a previous version.
Finally, thanks goes to Istv\'an Heckenberger for helpful discussions.

\newpage
\section{Local cohomology of equivariant vector bundles on $\P^d_K$}
\subsection{The rigid analytic variety $\cX$}
In this section we recall some geometric properties of Drinfeld's
upper half space $\cX$. We explain its rigid analytic structure
making it into a Stein space.  Furthermore, we treat briefly
$G$-equivariant vector bundles on $\cX$ which are induced by
homogeneous vector bundles on $\P^d_K.$ We discuss how we can
associate to to such a sheaf a locally analytic $G$-representation
in the sense of \cite{ST3}. In what follows, we denote for a variety
$X$ over $K$ by $X^{rig}$ the  rigid analytic variety attached to $X$ \cite{BGR}.

Let  $\epsilon \in \bigcup_{n\in \N}\sqrt[n]{|K^\times|} =
|\overline{K^\times}|$ be a $n$-th square root of
some absolute value in $|K^\times|$. Recall  the definition of an
open respectively closed $\epsilon$-neighborhood of a closed
$K$-subvariety $Y \subset \mathbb P_K^d$. Let
$f_1,\ldots,f_r \in S=K[X_0,\ldots, X_d]$ be finitely many
homogeneous polynomials with integral coefficients generating the
vanishing ideal of $Y$. We suppose that each polynomial has at least
one coefficient in $O_K^\times.$ Let $|\;\; |$ be the unique extension
of our fixed norm $|\;\;|$ on $K$  to $\mathbb C_p.$ A tuple $(z_0,\ldots,z_d)\in \mathbb A_K^{d+1}(\Cp)$ is called
unimodular if $|z_i|\leq 1 \; $ for $i=0,\ldots,d$, and $|z_i|=1$
for at least one $i$ with $0 \leq i \leq d.$ The open
$\epsilon$-neighborhood of $Y$ is defined by $$Y(\epsilon)=\Big\{z
\in \mathbb (\P_K^d)^{rig}\mid \mbox{ for any unimodular
representative }  \tilde{z} \mbox{ of $z,$  we have } $$
$$|f_j(\tilde{z})|\leq \epsilon
 \mbox{ for all } 1\leq j\leq r \Big\}.$$ This definition is
independent of the chosen unimodular representatives, so it is
well-defined. By using the standard covering $(D_+(X_i))_{i=0\ldots,
d}$ of $\P_K^d,$ one verifies that $Y(\epsilon)$ is a finite union
of $K$-affinoid spaces, cf. \cite{BGR} 7.2. In particular, it is a
quasi-compact open rigid analytic subspace of $(\P^d_K)^{rig}.$ On
the other hand, the set
$$Y^-(\epsilon)=\Big\{z \in \mathbb (\P_K^d)^{rig}\mid \mbox{  for any unimodular representative } \tilde{z} \mbox{ of $z,$  we have } $$
$$| f_j(\tilde{z})|< \epsilon
 \mbox{ for all } 1\leq j\leq r \Big\}$$
 is called the closed $\epsilon$-neighborhood of $Y.$ Again, it is an admissible open subset of $(\P^d_K)^{rig},$ but
 which is in general not quasi-compact.
 \begin{Remark}
{\rm We use the terminology open respectively closed, since the
corresponding neighborhoods for adic spaces \cite{H} are open
respectively closed in the adic space $(\P^d_K)^{ad}$. } \qed
 \end{Remark}
For a non-trivial linear $K$-subspace $U\subsetneq K^{d+1}=V,$ let $Y_U$
be the closed linear $K$-subvariety $\P(U)$ of $\P^d_K.$  Set $\epsilon_n:=|\pi^n|, n\in \N.$

\begin{Proposition}\label{Drinfeld} For every $n\in \N,$ both
$$(\P_K^d)^{rig}\setminus Y_U(\epsilon_n), \mbox{ for } U \subsetneq V,$$ and
$$\cY_n  :=  \bigcup_{U \varsubsetneq V} Y_U(\epsilon_n) \mbox{ resp. }   \cX_n^{-} :=   (\P^d_K)^{rig} \setminus \cY_n $$
are admissible open subsets of $(\P_K^d)^{rig},$ where $\cY_n$ is quasi-compact. The covering
$$\cX=\bigcup_{n \in \N}\cX^{-}_n$$ is admissible open and
$\cX$ is defined over $K.$
\end{Proposition}

\proof See \cite{SS} Proposition 1. \qed

\medskip
As in \cite{ST1} we also work with the closed
$\epsilon$-neighborhoods $Y^{-}_U(\epsilon).$ Similar to the above
proposition, we have for these spaces, the following statements.

\begin{Proposition}\label{Steinspace} For every $n\in \N,$ both
$$(\P_K^d)^{rig}\setminus Y_U^{-}(\epsilon_n), \mbox{ for } U \subsetneq V,$$ and
$$\cY^{-}_n  :=  \bigcup_{U \varsubsetneq V} Y^{-}_U(\epsilon_n)
\mbox{ resp. }
\cX_n  :=  (\P^d_K)^{rig} \setminus \cY^{-}_n $$
are admissible open subsets of $(\P_K^d)^{rig},$ where $\cX_n$ is quasi-compact. The covering
$$\cX=\bigcup_{n \in \N}\cX_n$$ is admissible open.
Furthermore, this covering induces on $\cX$ the structure of a Stein space. The $K$-algebra of global sections $\cO(\cX)$
is a $K$-Fr\'echet space. More precisely, we have $\cO(\cX) =\varprojlim_{n \in \N} \cO(\cX_n),$ where the $K$-algebras $\cO(\cX_n)$
are $K$-Banach spaces.
\end{Proposition}

\proof See \cite{SS} Proposition 4 resp. \cite{ST1} chapter 1. \qed

\bigskip
We follow the convention in \cite{ST1} and  consider the algebraic
action $m: {\bf G} \times \P^d_K \rightarrow \P^d_K$ of ${\bf G}$ on
$\P_K^d$ given by
$$g\cdot [q_0:\cdots :q_d]:=m(g,[q_0:\cdots :q_d]):= [q_0:\cdots :q_d]g^{-1}.$$
Let $\cF$ be a ${\bf G}$-equivariant  vector bundle  on $\P_K^d$. This is a vector bundle $\cF$ on $\P^d_K$ together with a  ${\bf G}$-linearization, i.e.,  an isomorphism of sheaves
\begin{eqnarray}\label{linearisation}
m^\ast(\cF) \stackrel{\sim}{\rightarrow} pr^\ast(\cF)
\end{eqnarray}
on ${\bf G} \times \P^d_K$ satisfying a certain cocycle condition,
cf.  \cite{MFK} Definition 1.6. Here $pr:{\bf G} \times \P^d_K
\rightarrow \P^d_K$ is the projection map onto the second factor.
We get by functoriality an induced ${\bf G}^{rig}$-equivariant vector
bundle on   $(\P_K^d)^{rig}$, which we denote for simplicity by
$\cF,$ as well.

Alternatively,  there is the following description of  ${\bf
G}$-equivariant vector bundles on $\P^d_K,$ cf. \cite{Bo}, \cite{Ja}
where they are called homogeneous vector bundles. Denote by ${\bf
P_{(1,d)}}$ the stabilizer of the base point $[1:0:\cdots:0] \in
\P^d_K(K),$ which is a parabolic subgroup of ${\bf G}.$ Let
$$\pi:{\bf G} \rightarrow {\bf G/P_{(1,d)}}$$ be the projection map
and identify ${\bf G /P_{(1,d)}}$ with $\P^d_K.$ Let $V$ be a
finite-dimensional algebraic representation of ${\bf P_{(1,d)}}.$
For a Zariski open subset $U\subset \P^d_K$, put
\begin{eqnarray*}\label{homvectorbundle}
\cF_V(U):= &\Big \{ & \mbox{algebraic morphisms } f:\; \pi^{-1}(U)
\rightarrow V \;\mid \; f(gp)\,=\,p^{-1}f(g)\; \mbox{ for all } \\ &
& g \in {\bf G}(\overline{K}), p\in {\bf
P_{(1,d)}}(\overline{K}) \Big\}.
\end{eqnarray*}
Then $\cF_V$ defines a homogeneous vector bundle on $\P^d_K$
with fibre $V.$ If $\cF$ is a ${\bf G}$-equivariant vector bundle with fibre $V$ in the base point, one has a natural identification $\cF\cong\cF_V.$

Our $p$-adic group $G$ stabilizes $\cX$. Therefore, we obtain an
induced action of $G$ on the $K$-vector space of rigid analytic holomorphic sections $\cF(\cX).$ Let
$\cO$ be the structure sheaf on $\P^d_K.$ Since $\cX$ is contained
in the rigid analytic variety attached to affine scheme $D_+(X_0) \cong \bA^d_K$, we may choose a
$K$-linear isomorphism
\begin{equation}\label{split}
 \cO(\cX)^{n} \stackrel{\sim}{\longrightarrow} \cF(\cX).
\end{equation}
Here the integer $n={\rm rk}(\cF) \in \N$ is the rank of $\cF.$ We transfer the
natural topology of the former one onto $\cF(\cX).$ The topology on
$\cF(\cX)$ is independent  of the chosen isomorphism. Thus
$\cF(\cX)$ inherits the structure of a $K$-Fr\'echet space.
Similarly, the sets $\cF(\cX_n)$ are $K$-Banach spaces and we get\footnote{Hence the $\cO(\cX)$-module $\cF(\cX)$ is  coadmissible in the sense of Schneider and Teitelbaum \cite{ST2}, cf. p. 152.}
$$\cF(\cX) = \varprojlim_n \cF(\cX_n).$$

Applying the same arguments as in \cite{ST1} Lemma 1.3, Proposition 1.4 and Proposition 2.1,
we conclude that $\cF(\cX)$ is a reflexive $K$-Fr\'echet space and its strong dual
$$\cF(\cX)' = \varinjlim_{n \in \N} \cF(\cX_n)'$$
is a locally convex inductive limit of duals of $K$-Banach spaces. Furthermore,
the action
$$G\times \cF(\cX) \rightarrow \cF(\cX) $$
is continuous and the orbit maps
\begin{eqnarray*}
G & \rightarrow & \cF(\cX)'  \\
g & \mapsto &  g\cdot  f  
\end{eqnarray*} 
are locally analytic for $f \in \cF(\cX)'.$ 
Thus, the strong dual $\cF(\cX)'$ is a locally analytic $G$-representation
in the sense of  Schneider and Teitelbaum  \cite{ST3}. By definition it is a barrelled locally
convex Hausdorff $K$-vector space together with a continuous action
of $G$ such that the orbit maps are locally analytic functions on
$G.$

\subsection{Algebraic local cohomology  I}
This section deals with the algebraic local cohomology  $H^\ast_{\P^j_K}(\P^d_K,\cF)$ of ${\bf G}$-equivariant vector bundles $\cF$ on 
$\P^d_K.$ We will study these $K$-vector spaces as representations of $\frg$ and of the parabolic subgroup fixing $\P^j_K.$

Denote by ${\bf B\subset G}$ the Borel subgroup of lower triangular
matrices and let  ${\bf U}$ be its unipotent radical. Let ${\bf
T\subset G}$ be the diagonal torus and denote by ${\bf
\overline{T}}$ its image in ${\rm \bf PGL}_{d+1}.$ For $0 \leq i \leq
d,$ let $\epsilon_i:{\bf T}\rightarrow {\bf \bbG_m}$ be the
character defined by $\epsilon_i({\rm diag}(t_1,\ldots,t_d))= t_i.$
Put $\alpha_{i,j}:=\epsilon_i - \epsilon_j$ for $i\neq j$, and
$\alpha_i:=\alpha_{i+1,i}$ for $0\leq i\leq d-1.$ Then
$$\Delta:=\{\alpha_i \,\mid\, 0\leq i \leq d-1\}$$ are the simple roots
and $$\Phi:=\{\alpha_{i,j}\,\mid\, 0\leq i\neq j \leq d-1\}$$ are
the roots of ${\bf G}$ with respect to ${\bf T\subset B}$. For a
decomposition $(i_1,\ldots,i_r)$ of $d+1,$ let ${\bf
P_{(i_1,\ldots,i_r)}}$ be the corresponding standard-parabolic
subgroup of ${\bf G}$, ${\bf U_{(i_1,\ldots,i_r)}}$ its unipotent radical, ${\bf
U^+_{(i_1,\ldots,i_r)}}$ its opposite unipotent radical and ${\bf
L_{(i_1,\ldots,i_r)}}$ its Levi component.

Fix  a  ${\bf G}$-equivariant vector bundle $\cF$ on
$\P^d_K.$ Let $F$ be a graded ${\bf G}$-module  which is projective
and of finite type over $S=K[X_0,\ldots,X_d],$ such that its
associated sheaf on $\P^d_K$ is just $\cF,$ cf. \cite{H} ch. 2, \S
5.  Then $\cF$ is  naturally a $\frg$-module, i.e., there is a
homomorphism of Lie algebras
\begin{eqnarray}\label{AktionLie}
\frg \rightarrow {\rm End}(\cF)
\end{eqnarray}
defined in  the following way. Restrict the linearization (\ref{linearisation}) to
$G^{(1)}\times \P^d_K,$ where $G^{(1)}$ is the first infinitesimal
neighborhood of the identity. Let $\frx \in \frg$ and let $f\in
\cF(U)$ be a section for a Zariski open subset $U\subset \P^d_K.$ 
Then
$$\frx \cdot f:= \frac{d}{dT}((1+ T\frx) \cdot f)_{T=0}.  $$
Further, there is the following Leibniz rule concerning the multiplication with functions $\Xi \in \cO_{\P^d_K}(U),$
\begin{eqnarray}\label{multLie}
\frx\cdot(\Xi \cdot f)= \Xi \cdot (\frx\cdot f) + (\frx\cdot \Xi) \cdot f .
\end{eqnarray}
Here we consider the structure sheaf $\cO=\cO_{\P^d_K}$ with its natural ${\bf G}$-linearization.
In this case  we can specify the action of $\frg$ on $\cO.$
Indeed, for a root $\alpha=\alpha_{i,j}\in \Phi,$ let
$$L_\alpha :=L_{(i,j)} \in \frg_\alpha$$ be the standard generator of the
weight space $\frg_\alpha$  in $\frg.$ Let $\mu \in
X^\ast(\overline{\bf T})$ be a character of the torus $\overline{\bf
T}.$ Write $\mu$ in the shape $\mu=\sum_{i=0}^d m_i \epsilon_i$ with
$\sum_{i=0}^d m_i=0.$ Define $\Xi_{\mu} \in \cO(\cX)$  by
$$\Xi_\mu(q_0,\ldots,q_d)=q_0^{m_0}\cdots q_d^{m_d}.$$
For these functions, the action  of $\frg$  is given by
\begin{equation}\label{Strukturgarbe}
L_{(i,j)}\cdot \Xi_\mu = m_j\cdot \Xi_{\mu + \alpha_{i,j}} 
\end{equation}
and
\begin{equation*}
t\cdot \Xi_\mu = (\sum\nolimits_ i m_i t_i)\cdot\Xi_\mu, \; t\in \frt.
\end{equation*}

Fix an integer $0 \leq j \leq d-1.$
Let $$\P^j_K=V(X_{j+1},\ldots,X_{d})\subset \P^d_K$$ be the closed $K$-subvariety
defined by the vanishing of the coordinates $X_{j+1},\ldots,X_{d}.$ The algebraic local cohomology modules $H^i_{\P^j_K}(\P^d_K,\cF), \; i\in \N$,
sit in a long exact sequence
$$\cdots \rightarrow H^{i-1}(\P^d_K\setminus \P^j_K,\cF) \rightarrow H^i_{\P^j_K}(\P^d_K,\cF) \rightarrow H^i(\P^d_K,\cF) \rightarrow H^i(\P^d_K\setminus \P^j_K,\cF) \rightarrow \cdots$$
Let $\cH^{\ast}_{\P^j_K}(\P^d_K,\cF)$ be the local cohomology sheaf with support in the closed subvariety $\P^j_K.$
It is related to the local cohomology groups by a spectral sequence (cf. \cite{SGA2}, Theorem 2.6.)
$$E_2^{p,q}=H^p(\P^d_K, \cH^{q}_{\P^j_K}(\P^d_K,\cF)) \Longrightarrow H^{p+q}_{\P^j_K}(\P^d_K,\cF).$$
Since $\P^j_K$ and $\P^d_K$ are both smooth, the local cohomology
groups $\cH^i_{\P^j_K}(\P^d_K,\cF)$ vanish for $i\neq d-j$ \cite{SGA2}.
It follows that
$$H^p(\P^d_K, \cH^{d-j}_{\P^j_K}(\P^d_K,\cF)) \cong H^{p+d-j}_{\P^j_K}(\P^d_K,\cF) $$
for all $p\in \N.$
In particular,
$$ H^{i}_{\P^j_K}(\P^d_K,\cF)=0 \, \mbox{ for } i< d-j.$$
On the other hand, the cohomology groups $H^{\ast}(\P^d_K \setminus
\P^j_K,\cF)$ can be computed by  the $\check{{\rm C}}$ech complex
\begin{align}
\nonumber
\bigoplus\limits_{j+1\leq k \leq d} \cF(D_+(X_k)) & \rightarrow  \bigoplus\limits_{j+1 \leq k_1<k_2 \leq d} \cF(D_+(X_{k_1} \cdot X_{k_2})) & \rightarrow  \cdots   \rightarrow & \cF(D_+(X_{j+1}\cdots X_d)) \\ \label{Cechalg} \\ =
\nonumber\bigoplus\limits_{j+1\leq k \leq d} (F_{X_k})^0 & \rightarrow   \bigoplus\limits_{j+1 \leq k_1 < k_2 \leq d} (F_{X_{k_1}\cdot X_{k_2}})^0 &\rightarrow  \cdots   \rightarrow &   (F_{X_{j+1}\cdots X_d})^0.
\end{align}
Here for a homogeneous polynomial $f\in S $, the set $D_+(f) \subset \P^d_K$  denotes as usual the Zariski open subset of $\P^d_K,$ where $f$ does not vanish.
The symbol $^0$ indicates the degree zero contribution of a graded module.

Alternatively, we may compute  $H^{d-j-1}(\P^d_K\setminus \P^j_K,\cF)$ for $j \leq d-2,$ by the inductive limit
\begin{eqnarray}\label{IndLim}
\varinjlim_{n\in \N} \, (F/ (X_{j+1}^n,\ldots, X_d^n) F)^0,
\end{eqnarray}
cf. \cite{EGAIII} Prop. 2.1.5.
In the case $j=d-1$, we have merely an exact sequence
$$0 \rightarrow H^{0}(\P^d_K,\cF) \rightarrow H^{0}(\P^d_K\setminus \P^{d-1}_K,\cF) \rightarrow \varinjlim_{n\in \N} (F/X_d^n F)^0 \rightarrow 0.  $$
Here, the (twisted) degree of a coset $[f] \in F/ (X_{j+1}^n,\ldots, X_d^n)$ is  by definition $$\deg([f])-n\cdot(d-j),$$
where $\deg([f])$ is the ordinary degree of $[f],$ cf. loc.cit. 2.1.
Sometimes, when we write down elements of $\varinjlim_{n\in \N} F/ (X_{j+1}^n,\ldots, X_d^n) F,$
we use generalized fractions. This terminology arises from the natural embedding of $K$-vector spaces 
$$\varinjlim_{n\in \N} \, (F/ (X_{j+1}^n,\ldots, X_d^n) F)^0 \hookrightarrow (F_{X_{j+1}\cdot \ldots \cdot X_d})^0 .$$ 
For $[f] \in (F/ (X_{j+1}^n,\ldots, X_d^n)F)^0,$ we use the symbol
$$ \left [ [f] \atop {X_{j+1}^n\cdot \ldots \cdot X_d^n} \right]$$
for its image in $(F_{X_{j+1}\cdot \ldots \cdot X_d})^0 $. The transition maps are then simply given by
$$\left [ [f] \atop {X_{j+1}^n \cdot \ldots \cdot X_d^n} \right] \mapsto \left [ [X_{j+1}\cdot \ldots \cdot  X_d \cdot f] \atop {X_{j+1}^{n+1}\cdot \ldots \cdot X_d^{n+1}} \right]. $$
We deduce from (\ref{Cechalg}) that
$H^{\ast}(\P^d_K \setminus \P^j_K,\cF)= 0 \mbox{ for } i\geq d-j$ and consequently
$$H^{i}_{\P^j_K}(\P^d_K,\cF)= H^{i}(\P^d_K,\cF) $$
for $i>d-j.$
In the case $i=d-j,$ we get an exact sequence
\begin{eqnarray}\label{exaktesequenzalg}
\nonumber 0 \rightarrow &  H^{d-j-1}(\P^d_K,\cF) & \rightarrow H^{d-j-1}(\P^d_K\setminus \P^j_K,\cF) \rightarrow  H^{d-j}_{\P^j_K}(\P^d_K,\cF) \\  \rightarrow   & H^{d-j}(\P^d_K,\cF)&  \rightarrow 0.
\end{eqnarray}
There is a natural algebraic action of ${\bf P_{(j+1,d-j)}}$ on each of the entries in this sequence, since the parabolic subgroup stabilizes $\P^j_K.$ In particular, this sequence is equivariant with respect to this action.
Furthermore, by (\ref{AktionLie}) the Lie algebra $\frg$ acts by functoriality  on all the cohomology groups, so that the sequence is equivariant for $\frg,$ too.  We set
$$\tilde{H}^{d-j}_{\P^j_K} (\P^d_K ,\cF):={\rm ker}\,\Big(H^{d-j}_{\P^j_K}(\P^d_K,\cF) \rightarrow H^{d-j}(\P^d_K,\cF)\Big)$$
which is consequently a ${\bf P_{(j+1,d-j)}} \ltimes U(\frg)$-module. Here the semi-direct product is defined via the adjoint action of ${\bf P_{(j+1,d-j)}}$ on $\frg.$ Indeed, for a section $f$ of $\cF$ and $\frz\in \frg$ resp. $p\in P_{(j+1,d-j)},$ we compute 
$$ \frac{d}{dT}((1+\frz T)(p\cdot f)) = \frac{d}{dT}(p\cdot (p^{-1}\cdot (1+\frz T) \cdot p)\cdot f))= p\cdot\frac{d}{dT}((1+(p^{-1} \cdot \frz \cdot p) T)\cdot  f)).$$
This compatibility transfers by functoriality onto the cohomology groups.

\begin{Lemma}\label{finite-dimensional} There exists a finite-dimensional  ${\bf P_{(j+1,d-j)}}$-invariant $K$-subspace
$$N_j\subset \tilde{H}^{d-j}_{\P^j_K} (\P^d_K ,\cF)$$
which generates $\tilde{H}^{d-j}_{\P^j_K} (\P^d_K ,\cF)$ as  $U(\frg)$-module.
\end{Lemma}

\proof Consider the formula (\ref{IndLim}). The parabolic subgroup ${\bf P_{(j+1,d-j)}}$ acts on each entry appearing in the inductive limit separately. Each entry is a finite-dimensional $K$-vector space.
Let 
$S/(X_{j+1}\ldots X_d) \rightarrow S$ be the $K$-linear section of the projection, given by $X_i\mapsto X_i$ for $i\leq j,$ and $X_i \mapsto 0$ for $i\geq j+1.$ This map induces a $K$-linear  section
$F/(X_{j+1},\ldots, X_d) F \rightarrow F.$ Denote by $F'$ the image of this section which forms consequently a system of representatives of $F/(X_{j+1},\ldots, X_d)F.$  Then $(F/(X_{j+1},\ldots, X_d)F)^0$ may be identified with the homogeneous elements $f \in F_1:=F'$ of degree $d-j.$ Similarly, $(F/(X_{j+1}^2,\ldots, X_d^2)F)^0$ may
be identified with the  homogeneous  elements $$f\in F_2:=F' \oplus \bigoplus_{k \geq j+1}^d X_k \cdot F'\oplus \bigoplus_{k,l \geq j+1 \atop k\neq l}^d X_k\cdot X_l \cdot F'\oplus \cdots \oplus X_{j+1}\cdots X_d\cdot F'$$ of degree
$2(d-j)$, etc. Under this identification the outer term $X_{j+1}\cdots X_d\cdot  F'$ coincides with the image of the first transition map in (\ref{IndLim}). Since $F$ is a finitely generated graded $S$-module, there is an integer $n \in \N$, such that any homogeneous representative $f \in F'$ of degree $n(d-j)$ is divisible by some monomial $X_0^{k_0}\cdot X_1^{k_1}\cdot \cdots X_j^{k_j}$ of degree $d-j.$
Set 
$$N_j={\rm im} \Big((F/(X_{j+1}^{n-1},\ldots, X_d^{n-1}))^0 \longrightarrow (F_{X_{j+1}\cdots X_d})^0\Big).$$ 
Thus
$$N_j=\Big\{  \left [ [f] \atop {X_{j+1}^{n-1}\cdot \ldots \cdot X_d^{n-1}} \right] \mid  f \in F_{n-1} \mbox{ homogeneous of degree } (n-1)(d-j)\Big\}. $$ 
We claim that this finite-dimensional  ${\bf P_{(j+1,d-j)}}$-invariant $K$-subspace satisfies the condition of our lemma.
In fact, let $f\in F_n$ be a homogeneous element of degree $n(d-j).$
By assumption, we may assume that
$$f=X_0^{k_0}\cdot X_1^{k_1}\cdot \cdots X_j^{k_j}\cdot X_{j+1}^{k_{j+1}}\cdots X_d^{k_d} \cdot g$$
with  $g\in F_{n-1}$ and $\sum_l k_l=d-j.$ Further, we may assume that $k_i\geq 1$ for at least one $i\leq j.$
Consider the identity
$$L_{(i,j+1)}\cdot \left[ [g] \atop X_{j+1}^{n-1}\cdot  \ldots \cdot X_d^{n-1} \right] =  \left[ {[L_{(i,j+1)}\cdot  g]} \atop X_{j+1}^{n-1}\cdot  \ldots \cdot X_d^{n-1} \right]- (n-1) \left[ [X_i\cdot X_{j+2}\cdots \cdot X_d \cdot g] \atop X_{j+1}^n\cdot \ldots \cdot X_d^n \right]  .$$
The left hand side and the first summand are contained in $U(\frg)\cdot N_j.$  It follows that $\left[ [X_i\cdot X_{j+2}\cdots \cdot X_d \cdot g] \atop X_{j+1}^n\cdot \ldots \cdot X_d^n \right]$ is contained in $U(\frg)\cdot N_j.$
By induction on the indices $l\leq j$ with $k_l \geq 1,$ we see that $\left[ [f] \atop X_{j+1}^{n}\cdot X_{j+2}^{n}\cdot \ldots \cdot X_d^{n} \right]\in  U(\frg)\cdot N_j.$ The case where $F\in F_{n'}$ with $n'> n$ follows inductively, as well. 
\qed

\begin{Remark}
{\rm Alternatively we can prove Lemma \ref{finite-dimensional}  by using Corollary \ref{irreduciblelagbraicrepr} in section 1.4.
This section is independent of the results in 1.2 and 1.3. In the case where ${\bf U =  U_{(1,d)}}$ acts trivially on the fibre $V$ of $\cF$, it produces  an  explicit candidate. In the general case, we know that the fix point set $V^{\bf U} \neq 0$ is non-trivial since ${\bf U}$ is unipotent, cf. \cite{Ja} ch.I, 2.14 (8). Consider the exact sequence
$$0\rightarrow V^{\bf U} \rightarrow V \rightarrow V/V^{\bf U} \rightarrow 0$$
of algebraic ${\bf P_{(1,d)}}$-modules, which induces an exact sequence of homogeneous vector bundles $$0\rightarrow \cF_{V^{\bf U}} \rightarrow \cF_V \rightarrow \cF_{V/V^{\bf U}} \rightarrow 0$$ on $\P^d_K$.  We get an  equivariant long exact sequence
\begin{eqnarray*}
0\rightarrow & H^{d-j}_{\P^j_K}(\P^d_K,\cF_{V^{\bf U}}) &  \rightarrow H^{d-j}_{\P^j_K}(\P^d_K,\cF_V) \rightarrow H^{d-j}_{\P^j_K}(\P^d_K,\cF_{V/V^{\bf U}}) \\
 \rightarrow  & H^{d-j+1}_{\P^j_K}(\P^d_K,\cF_{V^{\bf U}}) & \rightarrow \cdots
\end{eqnarray*}
of  ${\bf P_{(j+1,d-j)}} \ltimes U(\frg)$-modules. The groups $ H^{d-j}_{\P^j_K}(\P^d_K,\cF_W)$ and  $\tilde{H}^{d-j}_{\P^j_K}(\P^d_K,\cF_W), \, W\in\{V,V^{\bf U},V/V^{\bf U}\},$ differ only by the finite-dimensional $K$-vector space  $H^{d-j}(\P^d_K,\cF_W).$  By Corollary \ref{irreduciblelagbraicrepr} and by induction on the dimension of $V$,  there
are  finite-dimensional ${\bf P_{(i+1,d-i)}}$-submodules  of the outer terms generating them as $U(\frg)$-modules. But then the statement is true for the middle term $H^{d-j}_{\P^j_K}(\P^d_K,\cF_V)$ of the exact sequence and thus for $\tilde{H}^{d-j}_{\P^j_K}(\P^d_K,\cF_V)$ }. \qed
\end{Remark}

\subsection{Analytic local cohomology}
In the following we study the analytic local cohomology groups  $H^{d-i}_{\P^i_K(\epsilon)}(\P^d_K,\cF)$ as topological $K$-vector spaces. We shall prove  a local duality theorem which  describes the topological dual of  ${\rm ker}\big(H^{d-i}_{\P^i_K(\epsilon)}(\P^d_K,\cF) \rightarrow H^{d-i}(\P^d_K,\cF)\big)$ for $\epsilon \to 0$, by means of analytic functions on certain polydiscs.

Let $X$ be a rigid analytic variety over $K$ and consider a coherent sheaf $\cG$  on $X$. Let $U\subset X$ be an admissible open subset and let $Y:=X\setminus U$ be its set theoretical complement. Then the local algebraic cohomology groups $H^\ast_Y(X,\cG)$  are defined by the right derived functors
of $${\rm ker}\big(\Gamma(X,\cG) \rightarrow  \Gamma(U,\cG)\big).$$  If $X$ is a separated rigid analytic variety of countable type one can equip these cohomology groups with a locally convex topology as follows, cf. \cite{vP}.
For a rigid analytic variety $X$  of countable type, the space of global sections $\cG(X)=\Gamma(X,\cG)$ has a natural structure of a $K$-Fr\'echet space. 
If $X$ is an arbitrary separated rigid analytic variety of countable type with an admissible covering $X=\bigcup_i X_i$ by affinoids resp. by (quasi-) Stein
spaces \cite{K2}, one considers the corresponding $\check{{\rm C}}$ech complex
\begin{eqnarray*}
\prod_{i} \cG(X_i) \rightarrow  \prod_{i < j} \cG(X_i \cap X_j) \rightarrow \cdots  
\end{eqnarray*}
computing $H^{\ast}(X,\cG)$.
All contributions are $K$-Fr\'echet spaces, in particular, they are locally convex Hausdorff $K$-vector spaces.
Hence they induce on the cohomology groups  $H^{\ast}(X,\cG)$ in a natural way a locally convex topology. This topology does not depend on the covering $X=\bigcup X_i$, cf. \cite{Ba}
Lemma 1.32. We point out that the topology on the cohomology is in general not Hausdorff.   
Finally, we consider the long exact cohomology sequence
\begin{eqnarray*}
\cdots & \longrightarrow & H^{i}_{Y}(X,\cG)  \longrightarrow  H^{i}(X,\cG) \longrightarrow H^{i}(U,\cG) \\
 & \stackrel{\delta^i}{\longrightarrow} & H^{i+1}_{Y}(X,\cG) \longrightarrow \cdots
\end{eqnarray*} 
The cohomology groups $H^{i}_{Y}(X,\cG)$ are equipped with the finest locally convex topology such that the boundary maps $\delta^i$ become continuous. 
It turns out that the long exact cohomology sequence is then  even topological\footnote{ A topological exact sequence (or topological complex) of topological vector spaces is an  algebraic exact sequence (complex) $\cdots \rightarrow E^{i-1} \rightarrow E^i \rightarrow E^{i+1} \rightarrow \cdots$, such that all homomorphisms are continuous.  }  exact, cf. Lemma 5.1 in \cite{S2}.

We are interested in the analytic cohomology groups $H^\ast_{\P^j_K(\epsilon_n)}(\P^d_K,\cF)$ where $\cF$ is our fixed homogenous vector bundle.
Recall that $\epsilon_n = |\pi^n|,\; n\in \N.$ By GAGA (cf. \cite{K2} \S 4),  we know that
$$H^i(\P^d_K,\cF)=H^i((\P^d_K)^{rig},\cF)$$ for $i \geq 0.$
The cohomology group $H^i_{\P^j_K(\epsilon_n)}(\P^d_K,\cF)$  sits in the long exact cohomology sequence
\begin{eqnarray*}
\cdots & \rightarrow & H^{i-1}((\P^d_K)^{rig}\setminus \P^j_K(\epsilon_n),\cF) \rightarrow H^i_{\P^j_K(\epsilon_n)}(\P^d_K,\cF) \rightarrow H^i(\P^d_K,\cF) \\
& \rightarrow & H^i((\P^d_K)^{rig} \setminus \P^j_K(\epsilon_n),\cF) \rightarrow \cdots .
\end{eqnarray*}
As for the computation of $H^i((\P^d_K)^{rig}\setminus \P^j_K(\epsilon_n), \cF)$, we consider the $\Check{\rm C}$ech complex \bigskip
\begin{multline*}
\bigoplus_{j+1\leq k \leq d} \cF(D_+(X_k)_{\epsilon_n}^-) \rightarrow \!\!\!\!\!\bigoplus_{j+1 \leq k_1<k_2 \leq d} \cF(D_+(X_{k_1})_{\epsilon_n}^-\cap D_+(X_{k_2})_{\epsilon_n}^-)  \rightarrow \cdots  \\ \cdots \rightarrow  \cF(D_+(X_{j+1})_{\epsilon_n}^-\cap \cdots \cap D_+(X_d)_{\epsilon_n}^-)
\end{multline*}
with respect to the covering of Stein spaces 
$$(\P^d_K)^{rig}\setminus \P^j_K(\epsilon_n)=\bigcup_{k=j+1}^d D_+(X_k)_{\epsilon_n}^-,$$ where
$$D_+(X_k)_{\epsilon_n}^-:= \Big\{[x_0:\ldots:x_d] \in (\P^d_K)^{rig}\mid\, |x_k| > |x_l|\cdot \epsilon_n \;\;  \forall \, l \Big\}.$$ 
For $1\geq \epsilon >0$, we set 
$$D_+(X_k)_{\epsilon}:= \Big\{[x_0:\ldots:x_d] \in (\P^d_K)^{rig}\mid\, |x_k| \geq |x_l|\cdot \epsilon \;\;  \forall \, l \Big\}.$$
These are affinoid rigid analytic varieties and we can write
$$(\P^d_K)^{rig}\setminus (\P^j_K)(\epsilon)^-=\bigcup_{k=j+1}^d D_+(X_k)_{\epsilon}.$$
Thus, we get an admissible  covering
$$(\P^d_K)^{rig}\setminus \P^j_K(\epsilon_n)= \bigcup_{\epsilon \to \epsilon_n \atop \epsilon_n < \epsilon \in |\overline{K^\ast}|} (\P^d_K)^{rig}\setminus \P^j_K(\epsilon)^-$$
by quasi-compact admissible open subsets. 
Consider the $\Check{\rm C}$ech complex computing $H^\ast((\P^d_K)^{rig}\setminus \P^j_K(\epsilon)^-,\cF):$\medskip
\begin{multline}\label{Cechanal}
\bigoplus_{j+1\leq k \leq d} \cF(D_+(X_k)_{\epsilon}) \rightarrow \!\!\!\!\!\bigoplus_{j+1 \leq k_1<k_2 \leq d} \cF(D_+(X_{k_1})_{\epsilon}\cap D_+(X_{k_2})_{\epsilon}) \rightarrow \cdots  \\ \cdots \rightarrow \cF(D_+(X_{j+1})_{\epsilon}\cap \cdots \cap D_+(X_d)_{\epsilon})
\end{multline}

\medskip
\begin{Lemma}
The cohomology groups   $H^i((\P^d_K)^{rig}\setminus \P^j_K(\epsilon)^-, \cF)$ (resp. $H^i_{\P^j_K(\epsilon)^-}(\P^d_K,\cF)),$ $i,j=0,\ldots, d,$ are  $K$-Banach spaces in which the algebraic cohomology  $H^i(\P^d_K\setminus \P^j_K, \cF)$ (resp. $H^i_{\P^j_K}(\P^d_K,\cF)$)
is a dense subspace. 
\end{Lemma}

\proof First we treat the case $\cF=\cO.$ Set $\epsilon'= \frac{1}{\epsilon}.$
Consider the Gauss-Norm $|\;\, |_{\epsilon'}$ on the homogenous  localization $(K[X_0,\ldots,X_d]_{X_0\cdots X_d})^0=\cO(D_+(X_0\cdot\ldots \cdot X_d))$  given as follows. Let $f\in K[X_0,\ldots,X_d]$ be a homogeneous polynomial of degree $n(d+1)\,, n\in \N.$
Write $$f= \sum_{i_0 + \cdots + i_d=n(d+1)} a_{i_0\cdots i_d}X_0^{i_0}\cdot \ldots \cdot X_d^{i_d}.$$
Then $$|\frac{f}{(X_0 \cdots X_d)^n}|_{\epsilon'}= \max_{i_0,\ldots ,i_d} {|a_{i_0 \cdots i_d}|(\epsilon')^{r(i_0,\cdots, i_d)}}$$ 
Here $r(i_0,\cdots, i_d)=\sum_{i_j\geq n} (i_j-n).$
Then for every subset $\{i_1,\ldots,i_{r+1}\} \subset \{0,\ldots ,d\},$ the restriction  maps $$\cO(D_+(X_{i_1})\cap \cdots \cap \widehat{D_+(X_{i_j})}\cap \cdots \cap D_+(X_{i_{r+1}})) \rightarrow \cO(D_+(X_{i_1})\cap \cdots \cap D_+(X_{i_{r+1}}))$$ are isometries. Here the symbol $\widehat{D_+(X_{i_j})}$ indicates that we omit the open subset $D_+(X_{i_j})$ from the intersection. The images of the differentials in (\ref{Cechalg}) are closed, in particular the differentials are strict\footnote{Recall that a homomorphism $f: V \rightarrow W$ of topological vector spaces is strict if the induced homomorphism $V/{\rm ker } f \rightarrow {\rm im} f$ of topological vector spaces with the inherited topologies is a homeomorphism, cf. \cite{BGR}.}  . In fact, it suffices to show that the image of the maps 
$$ \bigoplus_{j=1}^{r+1}\cO(D_+(X_{i_1})\cap \cdots \cap \widehat{D_+(X_{i_j})}\cap \cdots \cap D_+(X_{i_{r+1}})) \stackrel{\delta}{\rightarrow} \cO(D_+(X_{i_1})\cap \cdots \cap D_+(X_{i_{r+1}}))$$
are closed. But
$${\rm  im}(\delta) = \Big\{ \sum_{i_0+\cdots + i_d=0} a_{i_0 \cdots i_d} X_0^{i_0}\cdots X_d^{i_d} \mid a_{i_0 \cdots i_d}  = 0\;  \mbox{ if } i_{k} < 0 \mbox{ for some } k \not\in \{i_1, \ldots,i_{r+1}\}$$ 
$$ \mbox{ resp. if $i_j < 0$ for all } j=1,\ldots, r+1  \Big\}.$$
Further, the completion of (\ref{Cechalg}) with respect to the norm $|\;\,|_{\epsilon'}$ is exactly (\ref{Cechanal}). By \cite{BGR} section 1.2, Cor. 6, the completion functor is exact for strict homomorphism. Hence the statement of our lemma follows in the case $\cF=\cO.$

For an arbitrary vector bundle $\cF$, we use the fact that it splits on the affine sets $D_+(X_i)$ as a direct sum of ${\rm rk }\, \cF$ copies of $\cO$.
Again the complex (\ref{Cechanal}) is the completion of (\ref{Cechalg}) and the images of the differentials are closed.
\qed

Now we treat the situation  of the open tubes $\P^d_K(\epsilon_n)\subset (\P^d_K)^{rig}.$ We shall see that $H^i((\P^d_K)^{rig}\setminus \P^j_K(\epsilon_n), \cF)$ respectively $H^i_{\P^j_K(\epsilon_n)}(\P^d_K,\cF)$ are naturally  $K$-Fr\'echet spaces in which the algebraic cohomology $H^i(\P^d_K\setminus \P^j_K, \cF)$ respectively $H^i_{\P^j_K}(\P^d_K,\cF)$
is a dense subset.  More precisely, we can write these cohomology groups as projective limits of $K$-Banach spaces:
\begin{Lemma}\label{Frechetspaces} 
We have
$$H^i((\P^d_K)^{rig}\setminus \P^j_K(\epsilon_n), \cF)= \varprojlim_{\epsilon \to \epsilon_n \atop \epsilon_n < \epsilon \in |\overline{K^\ast}|} H^i((\P^d_K)^{rig}\setminus \P^j_K(\epsilon)^-, \cF)$$
respectively
$$H^i_{\P^j_K(\epsilon_n)}(\P^d_K,\cF)= \varprojlim_{\epsilon \to \epsilon_n \atop \epsilon_n < \epsilon \in |\overline{K^\ast}|} H^i_{\P^j_K(\epsilon)^-}(\P^d_K,\cF).$$
\end{Lemma}

\proof In fact, the compatibility with the projective limit follows from
the following propositions. Here, the density condition follows from the previous lemma. \qed

\begin{Proposition}\label{projlimI}
Let $\cG$ be a coherent sheaf on a rigid analytic variety  $X.$ Consider a  decreasing family of  subsets  $Y_1 \supset Y_2 \supset \cdots \supset Y_k \supset Y_{k+1} \supset \cdots $ in $X$,
such that every subset $X\setminus Y_i$  is admissible open in $X.$
Set $Y:=\bigcap_{k \in \N} Y_k$ and assume that  $X\setminus Y$ is admissible open in $X$, as well.
Suppose that all cohomology groups $H^{i-1}_{Y_k}(X,\cG)$ are $K$-Fr\'echet spaces, such that the images of the transition maps $H^{i-1}_{Y_{k+1}}(X,\cG) \rightarrow H^{i-1}_{Y_k}(X,\cG)$ are dense for $k\in \N$.
Then there is a topological isomorphism $$\varprojlim\nolimits_{k \in \bbN} H^i_{Y_k}(X,\cG) \cong H^i_{Y}(X,\cG)$$
of $K$-Fr\'echet spaces.
\end{Proposition}

\proof By the same reasoning as in  Proposition 4 on \S 2 of \cite{SS} (cf. also Proposition \ref{Ext=projlim}), we have a short exact sequence
$$0\rightarrow \varprojlim_k\nolimits^{(1)} H^{i-1}_{Y_k}(X,\cG) \rightarrow H^i_{Y}(X,\cG) \rightarrow \varprojlim_k H^i_{Y_k}(X,\cG)\rightarrow 0.$$
But the projective system $(H^{i-1}_{Y_k}(X,\cG))_{k\in \N}$ of $K$-Fr\'echet spaces  has the topological Mittag-Leffler property by our condition on the density, cf. \cite{EGAIII} 13.2.4.
Thus we get by loc.cit. 13.2.3 an algebraic isomorphism $p:H^i_{Y}(X,\cG) \stackrel{\sim}{\rightarrow}  \varprojlim\nolimits_k H^i_{Y_k}(X,\cG).$ But $p$ is continuous and $\varprojlim\nolimits_k H^i_{Y_k}(X,\cG)$ is a $K$-Fr\'echet space. It follows from the bijectivity that $H^i_{Y}(X,\cG)$ has to be Hausdorff. Since it is a quotient of a $K$-Fr\'echet space it has to be a $K$-Fr\'echet space, as well. Now the claim follows from the open mapping theorem \cite{S2} 8.6.
\qed

Analogously one proves the "dual" version of this Proposition.

\begin{Proposition}\label{projlimII}
Let $\cG$ be a coherent sheaf on a rigid analytic variety  $X.$ Let $U\subset X$ be an admissible open subset and
consider an  increasing family of open admissible subsets  $U_1 \subset U_2 \subset \cdots \subset U_k \subset U_{k+1} \subset \cdots \subset U$ of $U$  with $\bigcup_{k \in \N} U_k = U.$
Suppose that all cohomology groups $H^{i-1}(U_{k},\cG)$  are $K$-Fr\'echet spaces and that the images of the transition maps $H^{i-1}(U_{k+1},\cG) \rightarrow H^{i-1}(U_k,\cG) $ are dense for  $k\in \N$.
Then there is a topological isomorphism  $$\varprojlim\nolimits_{k \in \bbN} H^i(U_k,\cG) = H^i(U,\cG).$$
\end{Proposition}

\begin{Remark} {\rm
In \cite{SS} Corollary 5 the authors consider a similar question concerning the compatibility with projective limits. They
deal with constant coefficients in which all the cohomology groups are finitely generated modules over an artinian ring, so that
the usual Mittag-Leffler property holds. \qed }

\end{Remark}

\begin{Remark} {\rm  An alternative way for proving Lemma \ref{Frechetspaces} is to apply the following Lemma to the $\Check{{\rm C}}$ech complex (\ref{Cechanal}).

\begin{Lemma}\label{TopML}   Let $0\rightarrow V^1_n \rightarrow V^2_n \rightarrow V^3_n \rightarrow 0,\; n \in \N,$  be a projective system of topological  exact sequences of $K$-Banach spaces (or more generally of $K$-Fr\'echet spaces). Suppose that the  transition maps $V^1_{n+1} \rightarrow V^1_n$, $n\in\N$,  have dense image. Then the sequence $$0\rightarrow \varprojlim_n V^1_n \rightarrow \varprojlim_n V^2_n \rightarrow \varprojlim_n V^3_n \rightarrow 0$$ is topological exact, too.
\end{Lemma}

{\noindent Proof:}  The exactness follows from the topological Mittag-Leffler property, cf. \cite{EGAIII}, 13.2.4. \qed }

\qed
\end{Remark}

It follows from Lemma \ref{Frechetspaces} that
\begin{equation}\label{lokKoho}
 H^{i}_{\P^j_K(\epsilon_n)}(\P^d_K,\cF)=0\; \mbox{ for } i<d-j
 \end{equation}
\mbox{and}
\begin{equation*}
H^{i}_{\P^j_K(\epsilon_n)}(\P^d_K,\cF)= H^{i}(\P^d_K,\cF)\; \mbox{ for } i>d-j.
\end{equation*}

\noindent Put $G_0={\bf G}(O_K).$ For any positive integer $n \in \N$, we consider the reduction map
\begin{equation}\label{reduction_map}
p_n:G_0\rightarrow {\bf G}(O_K/{(\pi^n)}).
\end{equation} 
Put
$$P_{(j+1,d-j)}^n:= p_n^{-1}\big({\bf P_{(j+1,d-j)}}(O_K/{(\pi^n)})\big).$$
This is a compact open subgroup of $G_0$ which stabilizes $\P^j_K(\epsilon_n).$  Again, as in the algebraic setting, we have
an exact $P_{(j+1,d-j)}^n \ltimes U(\frg)$-equivariant topological complex
\begin{eqnarray*}
0 & \rightarrow &  H^{d-j-1}(\P^d_K,\cF) \rightarrow H^{d-j-1}((\P^d_K)^{rig}\setminus \P^j_K(\epsilon_n),\cF) \rightarrow H^{d-j}_{\P^j_K(\epsilon_n)}(\P^d_K,\cF) \\
& \rightarrow &   H^{d-j}(\P^d_K,\cF) \rightarrow 0.
\end{eqnarray*} 

\begin{Proposition}
The action $P_{(j+1,d-j)}^n \times H^{d-j}_{\P^j_K(\epsilon_n)}(\P^d_K,\cF) \rightarrow  H^{d-j}_{\P^j_K(\epsilon_n)}(\P^d_K,\cF)$ is continuous.
\end{Proposition}
\begin{proof}
We follow the proof of Lemma 1.3 in \cite{ST1}. Since $H^{d-j}_{\P^j_K(\epsilon_n)}(\P^d_K,\cF)$ is a $K$-Fr\'echet space it is by the same reasoning as there enough to show that for $m>n,$ the orbit maps (into $K$-Banach spaces)  $P_{(j+1,d-j)}^m \rightarrow H^{d-j}_{\P^j_K(\epsilon_n)^-}(\P^d_K,\cF)$ are locally analytic. We may assume that $\cF=\cO,$  cf. Prop. 2.1' in loc.cit.  The cohomology group $H^{d-j}_{\P^j_K(\epsilon)^-}(\P^d_K,\cO)$ is a quotient of $H^0(g\cdot(D_+(X_{j+1})_{\epsilon}\cap \cdots \cap D_+(X_d)_{\epsilon}),\cO)$ for all $g\in P_{(j+1,d-j)}^m$. Let $g\in P_{(j+1,d-j)}^m$ and $F \in H^0(D_+(X_{j+1})_{\epsilon}\cap \cdots \cap D_+(X_d)_{\epsilon},\cO).$  We choose an open neighborhood $Q$ of $g$ in $P_{(j+1,d-j)}^m$ such that $h\cdot (D_+(X_{j+1})_{\epsilon}\cap \cdots \cap D_+(X_d)_{\epsilon}) =  g\cdot (D_+(X_{j+1})_{\epsilon}\cap \cdots \cap D_+(X_d)_{\epsilon})\, \forall h\in Q$. Then it suffices to see that the induced map
$Q \rightarrow H^0(p\cdot (D_+(X_{j+1})_{\epsilon}\cap \cdots \cap D_+(X_d)_{\epsilon}),\cO), \; h\mapsto h F,$ is locally analytic. We may assume that $g=1.$ Then  we apply the same argument as  in Prop. 2.1' loc.cit. 
\end{proof}

\begin{Corollary}\label{locallyanalytic}
The dual space $H^{d-j}_{\P^j_K(\epsilon_n)}(\P^d_K,\cF)'$  is a locally analytic $P_{(j+1,d-j)}^n$-representation.
\end{Corollary}
\begin{proof}
The dual space $H^{d-j}_{\P^j_K(\epsilon_n)}(\P^d_K,\cF)'$ is by Lemma \ref{Frechetspaces} the locally convex inductive limit 
$$H^{d-j}_{\P^j_K(\epsilon_n)}(\P^d_K,\cF)' =  \varinjlim_{\epsilon \to \epsilon_n \atop \epsilon_n < \epsilon \in |\overline{K^\ast}|} H^i_{\P^j_K(\epsilon)^-}(\P^d_K,\cF)'$$
of duals of $K$-Banach spaces.
In the proof of the previous proposition we have seen that the orbit maps $P_{(j+1,d-j)}^m \rightarrow H^i_{\P^j_K(\epsilon)^-}(\P^d_K,\cF)$ are
locally analytic. Thus the orbit maps on the dual space are locally analytic. The claim follows.
\end{proof}

\noindent Set
$$\tilde{H}^{d-j}_{\P^j_K(\epsilon_n)}(\P^d_K,\cF):= {\rm ker}\,\big(H^{d-j}_{\P^j_K(\epsilon_n)}(\P^d_K,\cF) \rightarrow H^{d-j}(\P^d_K,\cF)\big).$$ This $K$-Fr\'echet space has the structure of a  $P_{(j+1,d-j)}^n \ltimes U(\frg)$-module in which the algebraic cohomology $\tilde{H}^{d-j}_{\P^j_K}(\P^d_K,\cF)$ is a dense subspace.

We apply Lemma \ref{finite-dimensional} to obtain a $P_{(j+1,d-j)}$-invariant finite-dimensional $K$-subspace
$$N_j\subset \tilde{H}^{d-j}_{\P^j_K} (\P^d_K ,\cF)$$
which generates
$\tilde{H}^{d-j}_{\P^j_K} (\P^d_K ,\cF)$ as $U(\frg)$-module. Thus $\tilde{H}^{d-j}_{\P^j_K} (\P^d_K ,\cF)$
is a quotient of a generalized Verma module \cite{Le}. More precisely, there is an epimorphism
$$\varphi_j: U(\frg) \otimes_{U(\frp_{(j+1,d-j)})} N_j \rightarrow \tilde{H}^{d-j}_{\P^j_K} (\P^d_K ,\cF)$$
of $U(\frg)$-modules. 
Since the universal enveloping algebra splits into a tensor product $U(\frg)=U(\fru^+_{(j+1,d-j)})\otimes_K U(\frp_{(j+1,d-j)}),$ we may regard $\varphi_j$ as an epimorphism
\begin{eqnarray}\label{phij}
\varphi_j: U(\fru^+_{(j+1,d-j)})\otimes_K N_j \rightarrow \tilde{H}^{d-j}_{\P^j_K} (\P^d_K ,\cF).
\end{eqnarray}
Denote by $\frd_j={\rm ker}(\varphi_j)$ the kernel of this map.

Consider the affine algebraic group ${\bf U^+_{(j+1,d-j)}}.$ The Levi subgroup ${\bf L_{(j+1,d-j)}}$ stabilizes ${\bf U^+_{(j+1,d-j)}}$
with respect to the action of conjugation. Let
$$\Phi_j=\{\beta_1,\ldots,\beta_r\}$$ be the set of roots of  $\fru^+_{(j+1,d-j)}.$
The $K$-algebra  $\cO({\bf U^+_{(j+1,d-j)}})$ of algebraic functions on ${\bf U^+_{(j+1,d-j)}}$ may be viewed as  the polynomial $K$-algebra in the indeterminates $X_{\beta_1},\ldots,X_{\beta_r}.$
Consider the $L_{(j+1,d-j)}\cdot U^+_{(j+1,d-j)}$-equivariant pairing
\begin{eqnarray}\label{pairing1}  \cO({\bf U^+_{(j+1,d-j)}}) \times  U(\fru^+_{(j+1,d-j)}) &\rightarrow& K \\
\nonumber (f,\frz) &\mapsto & \frz \cdot f(1).
\end{eqnarray}
This is a non-degenerate pairing and induces therefore a $K$-linear $L_{(j+1,d-j)}\cdot U^+_{(j+1,d-j)}$-equivariant injection
$$\cO({\bf U^+_{(j+1,d-j)}}) \hookrightarrow {\rm Hom}_K(U(\fru^+_{(j+1,d-j)}),K) . $$

\noindent More concretely, this map  is given by
$$X_{\beta_1}^{i_1}\cdots X_{\beta_r}^{i_r} \mapsto (i_1)! \cdots (i_r)! \cdot (L^{i_1}_{\beta_1}\cdots L^{i_r}_{\beta_r})^\ast$$
where $$\big\{ (L^{i_1}_{\beta_1}\cdots L^{i_r}_{\beta_r})^\ast \mid
(i_1,\ldots,i_r)\in \N_0^r \big\}$$ is the dual basis of $\{
L^{i_1}_{\beta_1}\cdots L^{i_r}_{\beta_r}\mid(i_1,\ldots,i_r)\in
\N_0^r \}.$

\noindent Put
$$U^{+,n}_{(j+1,d-j)} = {\rm ker}\,\big({\bf U^{+}_{(j+1,d-j)} }(O_K)\rightarrow {\bf U^{+}_{(j+1,d-j)}}(O_K/{(\pi^n)})\big).  $$
Thus we have the identity
$$P_{(j+1,d-j)}^n= {\bf P}_{(j+1,d-j)}(O_K)\cdot U^{+,n}_{(j+1,d-j)}.$$
Further, we may interpret $U^{+,n}_{(j+1,d-j)} \subset U^+_{(j+1,d-j)}$ as an open $K$-affinoid polydisc, since  all entries $x$ in $U^{+,n}_{(j+1,d-j)}$ apart from the diagonal have norm $|x| \leq |\pi^n|.$
Hence the ring of $K$-analytic functions $\cO(U^{+,n}_{(j+1,d-j)})$ is a $K$-Banach algebra.
The pairing $(\ref{pairing1})$ extends by continuity to a non-degenerate ${\bf L}_{(j+1,d-j)}(O_K)\cdot U^{+,n}_{(j+1,d-j)}$-equivariant pairing
\begin{eqnarray}
\cO(U^{+,n}_{(j+1,d-j)}) \times  U(\fru^+_{(j+1,d-j)}) \rightarrow K
\end{eqnarray}

\noindent which in turn extends to a $P_{(j+1,d-j)}^n$-equivariant pairing
\begin{equation}\label{pairing2}
(\;,\;):  (\cO(U^{+,n}_{(j+1,d-j)}) \otimes N_j')  \times  (U(\fru^+_{(j+1,d-j)}) \otimes N_j) \rightarrow K.
\end{equation}
$$ (f\otimes\phi,\frz\otimes n) \mapsto \phi(n)\cdot\frz \cdot f(1) $$

\noindent Here, the subgroup ${\bf U}_{(j+1,d-j)}(O_K) \subset P_{(j+1,d-j)}^n$ acts by definition  trivially  on the $K$-Banach space 
$\cO(U^+_{(j+1,d-j)})$ respectively on $U(\fru^+_{(j+1,d-j)}).$
We put
$$\cO(U^{+,n}_{(j+1,d-j)} ,N_j')^{\frd_j}:=\Big\{ f \in \cO(U^{+,n}_{(j+1,d-j)})\otimes N_j' \,\mid (f,\frd_j)=0 \Big\}. $$

\noindent We obtain an equivariant injection
$$\cO(U^{+,n}_{(j+1,d-j)},N_j')^{\frd_j} \hookrightarrow {\rm Hom}_K(U(\fru^+_{(j+1,d-j)})\otimes N_j/ \frd_j,K)\cong{\rm Hom}_K(\tilde{H}^{d-j}_{\P^j_K} (\P^d_K ,\cF),K).$$
On the other hand, we have an injection of the duals
$$\tilde{H}^{d-j}_{\P^j_K(\epsilon_n)}(\P^d_K,\cF)' \hookrightarrow {\rm Hom}_K(\tilde{H}^{d-j}_{\P^j_K} (\P^d_K ,\cF),K),$$
since $\tilde{H}^{d-j}_{\P^j_K} (\P^d_K ,\cF)$ is dense in $\tilde{H}^{d-j}_{\P^j_K(\epsilon_n)}(\P^d_K,\cF).$ The following proposition says that for $n\rightarrow \infty,$ these two topological $K$-vector spaces coincide in  ${\rm Hom}_K(\tilde{H}^{d-j}_{\P^j_K} (\P^d_K ,\cF),K).$ It is based on the same principle as the duality theorem of  Morita, cf. \cite{Mo3} Theorem 2.
\begin{Proposition}\label{Duality}  For $n\in \N$ tending to infinity, we get an isomorphism of (Hausdorff) locally convex $K$-vector spaces
$$\varinjlim_{n\in \N} \cO(U^{+,n}_{(j+1,d-j)},N_j')^{\frd_j} \stackrel{\sim}{\longrightarrow}  \varinjlim_{n \in \N} \tilde{H}^{d-j}_{\P^j_K(\epsilon_n)}(\P^d_K,\cF)' $$
 compatible with the action of $\varprojlim_n P_{(j+1,d-j)}^n =  {\bf P_{(j+1,d-j)}}(O_K)$.
\end{Proposition}

\proof
Recall that we can express the $K$-Fr\'echet space $H^{d-j}_{\P^j_K(\epsilon_n)}(\P^d_K, \cF)$  by Lemma \ref{Frechetspaces}  as the projective limit of
the  $K$-Banach spaces $H^{d-j}_{\P^j_K(\epsilon)^-}(\P^d_K, \cF),$ where $\epsilon \to \epsilon_n$ and $\epsilon_n < \epsilon.$
Since $H^{d-j}(\P^d_K,\cF)$ is finite-dimensional, we see that $\tilde{H}^{d-j}_{\P^j_K(\epsilon_n)^-}(\P^d_K, \cF)$ is closed in $H^{d-j}_{\P^j_K(\epsilon_n)^-}(\P^d_K, \cF).$ We deduce the same compatibility  for the $K$-Fr\'echet space
$\tilde{H}^{d-j}_{\P^j_K(\epsilon_n)}(\P^d_K,\cF),$ i.e.,
$$\tilde{H}^{d-j}_{\P^j_K(\epsilon_n)}(\P^d_K, \cF)= \varprojlim_{\epsilon \to \epsilon_n \atop \epsilon_n < \epsilon \in |\overline{K^\ast}|} \tilde{H}^{d-j}_{\P^j_K(\epsilon)^-}(\P^d_K, \cF).$$
Therefore, we can replace the $K$-Fr\'echet spaces in the statement by the $K$-Banach spaces
$\tilde{H}^{d-j}_{\P^j_K(\epsilon_n)^-}(\P^d_K,\cF).$ We set for $\epsilon \in |\overline{K^\ast}|,$
\begin{eqnarray*}
U(\fru^+_{(j+1,d-j)})_{\epsilon}  := & \Big\{& \sum_{(i_1,\ldots,i_r) \in \N^r_0} a_{i_1,\ldots,i_r} L^{i_1}_{\beta_1}\cdots L^{i_r}_{\beta_r} \mid a_{i_1,\ldots,i_r} \in K,\; \\
& & |(i_1)! \cdots (i_r)! \cdot  a_{i_1,\ldots,i_r}|\epsilon^{i_1+\cdots + i_r} \rightarrow 0, i_1+\cdots + i_r \rightarrow \infty \Big\}.
\end{eqnarray*}
This is a $K$-Banach algebra in which the universal enveloping algebra $U(\fru^+_{(j+1,d-j)})$ is a dense subset. We get an epimorphism
of $K$-Banach spaces (use  \cite{BGR} Cor. 6 in 1.2), 
\begin{equation}\label{epi}
U(\fru^+_{(j+1,d-j)})_{\frac{1}{\epsilon_n}} \otimes N_j \longrightarrow   \tilde{H}^{d-j}_{\P^j_K(\epsilon_n)^-}(\P^d_K,\cF).
\end{equation}
On the other hand, let
\begin{eqnarray*}\cO_b(U^{+,n}_{(j+1,d-j)}):=\Big\{\sum_{(i_1,\ldots,i_r) \in \N^r_0}\!\!\!\!\!\!\!\! a_{i_1,\ldots,i_r} X^{i_1}_{\beta_1}\cdots X^{i_r}_{\beta_r} \mid a_{i_1,\ldots,i_r} \in K,\; \!\!\!\!\!\sup_{(i_1,\ldots,i_r)}\!\! |a_{i_1,\ldots,i_r}|\epsilon_n^{i_1+\cdots + i_r} < \infty \Big\}
\end{eqnarray*}
resp.
$$\cO_b(U^{+,n}_{(j+1,d-j)} ,N_j')^{\frd_j}:=\Big\{ f \in \cO_b(U^{+,n}_{(j+1,d-j)})\otimes N_j'\,\mid (f,\frd_j)=0 \Big\} $$
be the $K$-Banach spaces of bounded functions on $U^{+,n}_{(j+1,d-j)}$. Then
$$\varinjlim\nolimits_n \cO_b(U^{+,n}_{(j+1,d-j)})  = \varinjlim\nolimits_n \cO(U^{+,n}_{(j+1,d-j)})$$
resp.
$$\varinjlim\nolimits_n \cO_b(U^{+,n}_{(j+1,d-j)},N_j')^{\frd_j}  = \varinjlim\nolimits_n \cO(U^{+,n}_{(j+1,d-j)},N_j')^{\frd_j}$$
are identities of locally convex $K$-vector spaces.  But
$\cO_b(U^{+,n}_{(j+1,d-j)})$ is the topological dual of
$U(\fru^+_{(j+1,d-j)})_{\frac{1}{\epsilon_n}}$ (cf. Example in
\cite{S2} ch. I, \S 3). We deduce from (\ref{pairing2}) together
with (\ref{epi}) that  $\cO_b(U^{+,n}_{(j+1,d-j)} ,
N_j')^{\frd_j}$ is the topological dual of
$\tilde{H}^{d-j}_{\P^j_K(\epsilon_n)^-}(\P^d_K, \cF).$  The claim
follows now from \cite{Mo1}  Theorem 3.4  respectively  \cite{S2}
Prop. 16.10 on the duality of projective limits of $K$-Fr\'echet
spaces and injective limits of $K$-Banach spaces with compact
transition maps. \qed

\begin{Remark}\label{Bemerkung} 
{\rm The inductive limit $\varinjlim_{n \in \N} \tilde{H}^{d-j}_{\P^j_K(\epsilon_n)}(\P^d_K,\cF)' $ identifies by Proposition \ref{projlimI} with
the strong dual  of the analytic local cohomology group $\tilde{H}^{d-j}_{(\P^j_K)^{rig}}((\P^d_K)^{rig},\cF).$ In particular, the action of ${\bf P_{(j+1,d-j)}}(O_K)$ on  $\varinjlim_{n \in \N} \tilde{H}^{d-j}_{\P^j_K(\epsilon_n)}(\P^d_K,\cF)'$ extends to one of
$P_{(j+1,d-j)}.$ On the other hand, the expression  $\varinjlim_{n\in \N} \cO(U^{+,n}_{(j+1,d-j)},N_j')^{\frd_j}$  can been thought as the stalk in the point $\P^j_K$ of the Grassmannian ${\rm Gr}_{j+1}(K^{d+1})$ of a certain "sheaf". Here the action extends to one of $P_{(j+1,d-j)},$ as well.  It is easily seen that the isomorphism of Proposition \ref{Duality} is even $P_{(j+1,d-j)}$-equivariant, where the unipotent radical $U_{(j+1,d-j)}$ acts on the pairings (\ref{pairing2}) via $N_j'$. Furthermore, it follows from \ref{locallyanalytic} that the map is even an isomorphism of locally analytic $P_{(j+1,d-j)}$-representations.  
 \qed}
\end{Remark}

\subsection{Algebraic local cohomology II}

Let $\cF$ be a homogeneous vector bundle on $\P^d_K$ which arises by a representation of the Levi subgroup ${\bf L_{(1,d)}}$
of  ${\bf P_{(1,d)}},$ i.e., such that the unipotent radical ${\bf U_{(1,d)}}$ acts trivially on the fibre.
In this section we  compute explicit formulas for the $K$-vector spaces.
$$\tilde{H}^i_{\P^{d-i}_K}(\P^d_K,\cF)= {\rm ker}\,\big(H^{i}_{\P^{d-i}_K}(\P^d_K,\cF) \rightarrow
 H^{i}(\P^d_K,\cF)\big)$$ as representations of ${\bf P_{(d-i+1,i)}} \ltimes U(\frg).$
First we consider the local cohomology groups with respect to the closed subschemes $V(X_0,\ldots,X_{i-1}) \subset \P^d_K$ defined by the vanishing  of the first $i$ coordinate functions.
Note that the stabilizer of this subvariety is the upper triangular parabolic subgroup ${\bf P^+_{(i,d+1-i)}} \subset {\bf G}.$ Afterwards, we use the conjugacy
of   $V(X_0,\ldots,X_{i-1})$ and  $V(X_{d-i+1},\ldots,X_d)$ within  $\P^d_K$ via the action of ${\bf G}$ on $\P^d_K.$
The reason is that we follow the notation used by Kempf in \cite{Ke}.

Let $\pi:{\bf G}  \rightarrow {\bf G/P_{(1,d)}}$ be the projection map and identify ${\bf G /P_{(1,d)}} \cong \P^d_K$ as described in  section one.
Let $$\lambda'=(\lambda_1 \geq \ldots \geq \lambda_d)\in \Z^d$$ be a dominant integral weight of ${\rm {\bf GL_d}}.$
This gives rise to a finite-dimensional irreducible algebraic  representation $V_{\lambda'}$ of ${\rm {\bf GL_d}}.$
Let $\lambda_0 \in \Z$ be arbitrary and put
$$\lambda= (\lambda_0,\lambda_1,\ldots,\lambda_d)\in \Z^{d+1}.$$
 Extend the action of ${\bf L_{(1,d)}}$ on $V_{\lambda'}$ to one of ${\bf P_{(1,d)}}$ on the same space $V_\lambda:= V_{\lambda'}$, such that ${\bf U_{(1,d)}}$ acts trivially on it and such that ${\bf \bbG_m}$ acts via multiplication $$(x,v)\mapsto x^{\lambda_0}\cdot v,$$ 
$v\in V_{\lambda}\otimes_K \overline{K}, \; x\in {\bf \bbG_m}(\overline{K}).$ Denote by $\cF_\lambda= \cF_{V_\lambda}$ the corresponding homogeneous vector bundle on $\P^d_K$, cf. (\ref{homvectorbundle}). 
Furthermore, if we add to $\lambda$ the tuple $r\cdot(1,\ldots,1),\,r\in \N,$ then the ${\bf G}$-linearization on $\cF_\lambda$ is twisted by  $\det^{\otimes r}.$  Finally, we point out that
$H^0(D_+(X_0),\cF_\lambda)$ is isomorphic to $K[\frac{X_1}{X_0},\ldots,\frac{X_d}{X_0}] \otimes V_\lambda$ as 
${\bf P^+_{(1,d)}} \ltimes U(\frg)$-module.

\begin{Example}\label{Exampleirred}{\rm  The following identifications can be seen by the procedure used in \cite{Ja}, part II, 2.16. Note that we work with the contragredient identification
of ${\bf G/P_{(1,d)}}$ with $\P^d_K.$
\begin{enumerate}
\item Let $\lambda=(0,\ldots,0).$ Then $\cF_\lambda =\cO$ is the structure sheaf on $\P^d_K$.

\medskip
\item Let $\lambda=(r,0,\ldots,0),\; r \in \Z.$ Then $\cF_\lambda =\cO(r)$ is a twisted sheaf.

\medskip
\item Let $\lambda=(-1,1,0,\ldots,0).$ Then $\cF_\lambda =\Omega^1$ is the cotangent sheaf  on $\P^d_K$.
\medskip
\item Let $\lambda=(-d,1,\ldots,1).$ Then $\cF_\lambda =\Omega^d$ is the canonical bundle  on $\P^d_K$. \qed
\end{enumerate} }
\end{Example}
Let $W$ be the Weyl group of ${\bf G}$.
Set
$$w_i:=s_i\cdot s_{i-1}\cdot \cdots \cdot s_1 \in W,$$ where  $s_i\in W$ is the simple reflection with respect to the simple root $\alpha_i\in \Delta.$ We put $w_0=1.$
Let $W_{{\bf L_{(1,d)}}}$ be the Weyl group of ${\bf L_{(1,d)}}.$
Then the reflections $w_i\, ,i=0,\ldots,d-1,$ are just the representatives of shortest length in $W$ with respect to  $W/W_{\bf L_{(1,d)}}.$
Let ${\bf B^+ \subset G}$ be the Borel subgroup of upper triangular matrices.  The Schubert cells in $\P^d_K$ are given
by
\begin{eqnarray*}
X_{w_0}:={\bf B^+}w_0 {\bf P_{(1,d)}/ P_{(1,d)}} & = & D_+(X_0)  \\ \\
X_{w_1}:={\bf B^+}w_1 {\bf P_{(1,d)}/ P_{(1,d)}} & = & V(X_0) \cap D_+(X_1)  \\ \\
X_{w_2}:={\bf B^+}w_2 {\bf P_{(1,d)}/ P_{(1,d)}} & = & V(X_0,X_1) \cap D_+(X_2)  \\
& \vdots & \\
X_{w_d}:={\bf B^+}w_d {\bf P_{(1,d)}/ P_{(1,d)}} & = & V(X_0,X_1,\ldots, X_{d-1}) \cap D_+(X_d) = \{[0:0: \cdots : 1] \}
\end{eqnarray*}
The corresponding Schubert varieties $\overline{X}_{w_i},$ i.e., the Zariski closures of the cells $X_{w_i}$ are just the closed subschemes $V(X_0,\ldots,X_{i-1}) \subset \P^d_K ,\; 0\leq i\leq d.$

Denote by $\ast$ the dot action of $W$ on $X^\ast({\bf T})_{\bbQ}$  given by
$$w\ast \chi= w(\chi+\rho)-\rho,$$
where $\rho=\frac{1}{2}\sum_{\alpha \in \Phi^-} \alpha.$ Note that the set of negative roots $\Phi^-$ correspond to the set of positive roots with respect to the Borel subgroup ${\bf B^+}.$
We get
\begin{eqnarray*}
w_0 \ast \lambda & = &  \lambda \\ \\
w_1 \ast \lambda & = &  (\lambda_1-1,\lambda_0 +1 , \lambda_2,\ldots,\lambda_d) \\ \\
w_2 \ast \lambda & = &  (\lambda_1-1,\lambda_2-1,\lambda_0 +2 , \lambda_3,\ldots,\lambda_d) \\
&\vdots & \\
w_i \ast \lambda & = &  (\lambda_1-1,\lambda_2-1,\ldots, \lambda_i-1,\lambda_0 +i , \lambda_{i+1},\ldots,\lambda_d) \\
& \vdots & \\
w_d \ast \lambda & = &  (\lambda_1-1 ,\lambda_2-1,\ldots,\lambda_d-1,\lambda_0+d)\,.
\end{eqnarray*}
In particular, there is at most one integer $0\leq i \leq d$, such that $w_{i}\ast \lambda$ is dominant with respect to ${\bf B^+}.$
This integer is characterized by the non-vanishing of $H^{i}(\P^d_K,\cF_\lambda)$, cf. \cite{Bo} Theorem IV'.
We denote this integer by $i_0$  if it exists. Otherwise, there is a unique integer $i_0 < d$ with $w_{i_0}\ast \lambda = w_{i_0+1}\ast \lambda.$
We get
\begin{equation}\label{domordkl}
w_{i}\ast \lambda \succ  w_{i+1}\ast \lambda
\end{equation}
for all $i\geq i_0$ (resp. $i > i_0$ if $w_{i_0}\ast \lambda = w_{i_0+1}\ast \lambda$), and
\begin{equation}\label{domordgr}
w_{i}\ast \lambda \prec w_{i+1}\ast \lambda
\end{equation}
for all $i<i_0$, with respect to the dominance order $\succ$ on $X^\ast({\bf T})_{\bbQ}.$
We put
$$\mu_{i,\lambda}:=\left\{ \begin{array}{cc} w_{i-1}\ast \lambda  &: i \leq i_0 \\  w_{i} \ast \lambda & : i > i_0. \end{array} \right.$$
This is a ${\bf L_{(i,d-i+1)} }$-dominant  weight. Let $V_{i,\lambda}$ be the finite dimensional ${\bf L_{(i,d-i+1)} }$-module
with highest weight $\mu_{i,\lambda}.$ By considering the trivial action of ${\bf U^+_{(i,d-i+1)}}$ on it, we may view it as a ${\bf P^+_{(i,d-i+1)}}$-module.

\begin{Proposition}\label{Propositionweight} 
For $1\leq i \leq d,$ the ${\bf P^+_{(i,d-i+1)}} \ltimes U(\frg)$-module
$\tilde{H}^{i}_{\overline{X}_{w_i}}(\P^d_K,\cF_\lambda)$ is 
a quotient of the   ${\bf P^+_{(i,d-i+1)}}$-module\footnote{As for the ${\bf P^+_{(i,d-i+1)}}$-module  structure on $\bigoplus K\cdot X_0^{k_0}X_1^{k_1}\cdots X_d^{k_d}$  we refer to 3.1. In the proof we will see that we can realize $V_{i,\lambda}$ as a submodule of $\tilde{H}^{i}_{\overline{X}_{w_i}}(\P^d_K,\cF_\lambda)$.}  $\bigoplus_{k_0,\ldots, k_{i-1} \leq 0  \atop {k_{i},\ldots, k_d \geq 0  \atop {k_0+\cdots + k_d=0}}} K\cdot X_0^{k_0}X_1^{k_1}\cdots X_d^{k_d} \otimes  V_{i,\lambda} .$ 
\end{Proposition}

\proof Set $\cF=\cF_\lambda.$ We consider the Grothendieck-Cousin complex of $\cF$ with respect to the covering $(X_{w_i})_{i=0,\ldots,d}$ of $\P^d_K,$ i.e., the complex
$$0 \rightarrow H^{0}_{X_{w_0}}(\P^d_K,\cF) \stackrel{\delta_0}{\rightarrow} H^{1}_{X_{w_1}}(\P^d_K,\cF) \stackrel{\delta_1}{\rightarrow} \cdots  \stackrel{\delta_{d-1}}{\rightarrow}  H^{d}_{X_{w_d}}(\P^d_K,\cF) \rightarrow 0.$$
The $i$-th cohomology of this complex yields exactly $H^i(\P^d_K,\cF),$ cf. \cite{Ke} Theorem 8.7.  Furthermore, it is compatible with the action of ${\bf B^+}$ and
$\frg.$
We have for each $0\leq i \leq d-1,$ an exact sequence 
\begin{eqnarray*}
0 & \rightarrow & H^{i}_{\overline{X}_{w_i}}(\P^d_K,\cF) \rightarrow H^{i}_{X_{w_i}}(\P^d_K,\cF) \rightarrow H^{i+1}_{\overline{X}_{w_{i+1}}}(\P^d_K,\cF) \\  & \rightarrow & H^{i+1}_{\overline{X}_{w_i}}(\P^d_K,\cF) \rightarrow 0.
\end{eqnarray*}
Since $H^{i+1}_{\overline{X}_{w_{i+1}}}(\P^d_K,\cF)\rightarrow H^{i+1}_{X_{w_{i+1}}}(\P^d_K,\cF)$ is injective, we see that $H^{i}_{\overline{X}_{w_i}}(\P^d_K,\cF)$
is the kernel of $H^i_{X_{w_i}}(\P^d_K,\cF) \stackrel{\delta_i}{\rightarrow} H^{i+1}_{X_{w_{i+1}}}(\P^d_K,\cF).$
It follows that
$$\tilde{H}^i_{\overline{X}_{w_i}}(\P^d_K,\cF) = {\rm ker}\,( H^i_{\overline{X}_{w_i}}(\P^d_K,\cF) \rightarrow H^i(\P^d_K,\cF))$$
is isomorphic to the image of the boundary homomorphism
$$H^{i-1}_{X_{w_{i-1}}}(\P^d_K,\cF) \stackrel{\delta_{i-1}}{\rightarrow} H^{i}_{X_{w_i}}(\P^d_K,\cF) $$
resp.  to the cokernel of
$$H_{X_{w_{i-2}}}^{i-2}(\P^d_K,\cF) \stackrel{\delta_{i-2}}{\rightarrow} H^{i-1}_{X_{w_{i-1}}}(\P^d_K,\cF)$$
if $H^{i-1}(\P^d_K,\cF) \neq 0.$
Here we put $H^{-1}_{X_{w_{-1}}}(\P^d_K,\cF)= H^{0}(\P^d_K,\cF).$
By excision we get the following chain of ${\bf T} \ltimes \frg$-isomorphisms
\begin{eqnarray*}
H^{i}_{X_{w_i}}(\P^d_K,\cF) & \cong & H^{i}_{V(X_0,\ldots, X_{i-1})\cap D_+(X_i)}(D_+(X_i),\cF) \\
& \cong & H^{i}_{V(X_0,\ldots, X_{i-1})\cap D_+(X_i)}(D_+(X_i),\cO) \otimes \cF(w_i),
\end{eqnarray*} 
where $\cF(w_i)$ denotes the fibre of $\cF$ in $w_i\cdot {\bf P_{(1,d)}}.$ A simple computation gives
$$H^{i}_{V(X_0,\ldots, X_{i-1})\cap D_+(X_i)}(D_+(X_i),\cO) =\bigoplus_{k_0,\ldots, k_{i-1} < 0  \atop {k_{i+1},\ldots, k_d \geq 0  \atop {k_i \in \Z \atop k_0+\cdots + k_d=0}}} K\cdot X_0^{k_0}X_1^{k_1}\cdots X_d^{k_d}.$$
We can rewrite this expression as follows. Recall  that for a root $\alpha=\alpha_{k,l} \in \Phi$ the symbol $L_{(k,l)}$ denotes the standard generator
of the weight space $\frg_\alpha\subset \frg$. Analogously, we put 
$$X_{(k,l)}:=X_{\alpha_{k,l}}:=\frac{X_k}{X_l}\in K[X_0,\ldots,X_d]_{X_0\cdots X_d}.$$ 
Then we get (compare \cite{Ke} Corollary 11.10)
$$ H^{i}_{X_{w_i}}(\P^d_K,\cO)= K[X_{(i+1,i)},\ldots, X_{(d,i)}] \otimes \!\!\!\!\! \sum_{(n_0,\ldots,n_{i-1})\in \N^i}\!\!\!\!\! L_{(i,0)}^{n_0}\cdot L_{(i,1)}^{n_1}\cdots L_{(i,i-1)}^{n_{i-1}} \cdot \frac{X_i^i}{X_0\cdots X_{i-1}}.$$
Thus we obtain
\begin{eqnarray*}
H^{i}_{X_{w_i}}(\P^d_K,\cF)  & =  &  H^{i}_{X_{w_i}}(\P^d_K,\cO)\otimes \cF(w_i) \\ &=& K[X_{(i+1,i)},\ldots, X_{(d,i)}] \otimes \!\!\!\!\! \sum_{(n_0,\ldots,n_{i-1})} \!\!\!\!\!\!\!\! L_{(i,0)}^{n_0}\cdots L_{(i,i-1)}^{n_{i-1}} \cdot \frac{X_i^i}{X_0\cdots X_{i-1}} \otimes \cF(w_i).
\end{eqnarray*} 
The weights of  $H^{i}_{X_{w_i}}(\P^d_K,\cO)$ are given by
$$\Big\{w_i\ast (0,\ldots,0) - \,n_0 \cdot \alpha_{0,i} - \cdots - \,n_{i-1}\cdot  \alpha_{i-1,i} -\, n_i\cdot  \alpha_{i,i+1} - \cdots - \,n_{d-1} \cdot \alpha_{i,d}\,\mid\, n_0,\ldots,n_{d-1}\in \N \Big\}.$$ Here the highest weight is $w_i\ast (0,\ldots,0).$
The highest weight of the fibre $\cF(w_i)$ is given by  $w_i\cdot \lambda.$
We conclude that the highest weight of $H^{i}_{X_{w_i}}(\P^d_K,\cF)$ is given by $w_i\ast \lambda.$
A highest weight vector is given by $$v_{i,\lambda}=\frac{X_i^i}{X_0\cdots X_{i-1}}\otimes w_i(v_\lambda),$$ where
$v_\lambda=v_{0,\lambda}$ is a highest weight vector of $V_\lambda.$

Reconsider  the homomorphism 
$$H^{i-1}_{X_{w_{i-1}}}(\P^d_K,\cF) \stackrel{\delta_{i-1}}{\rightarrow} H^{i}_{X_{w_i}}(\P^d_K,\cF).$$
If $i<i_0$ then by (\ref{domordgr}) the highest weight of the image is $w_{i-1}\ast \lambda.$ If  $i\geq i_0$ then by (\ref{domordkl}) the highest weight of the image is $w_{i}\ast \lambda.$ Thus the highest weight of the image is $\mu_{i,\lambda}\in X^\ast({\bf T})$.  
Furthermore, the weight vectors $\delta_{i-1}(v_{i-1,\lambda})$ and $v_{i,\lambda}$ differ by the factor $(\frac{X_i}{X_{i-1}})^{\lambda_0-\lambda_i+i}$. More precisely, we conclude by weight reasons
\begin{eqnarray}
\nonumber\delta_{i-1}(v_{i-1,\lambda})=(\frac{X_i}{X_{i-1}})^{\lambda_0-\lambda_i+i} v_{i,\lambda} & \mbox{ if } & w_{i-1}\ast\lambda \preceq w_i\ast\lambda ,\\ 
\label{weights}\delta_{i-1}((\frac{X_i}{X_{i-1}})^{\lambda_0-\lambda_i+i} v_{i-1,\lambda}) = v_{i,\lambda} & \mbox{ if } & w_{i-1}\ast\lambda \succeq w_i\ast\lambda .
\end{eqnarray}
Since  $\tilde{H}^{i}_{\overline{X}_{w_i}}(\P^d_K,\cF)$ is also a
${\bf P^+_{(i,d-i+1)}}$-module,  it follows that it contains
the irreducible algebraic ${\bf L_{(i,d-i+1)}}= {\rm {\bf GL_i}} \times {\rm {\bf GL_{d-i+1}} }$-representation $V_{i,\lambda}$ corresponding to the highest weight $\mu_{i,\lambda}.$  One checks  that ${\rm im}(\delta_{i-1}) = \tilde{H}^{i}_{\overline{X}_{w_i}}(\P^d_K,\cF)$ is equal to\footnote{Note that $U(L_{(i-1,0)},\ldots,L_{(i-1,i-2)})$ leaves $V_{i,\lambda}$ invariant,
since $V_{i,\lambda}$ is a ${\bf P^+_{(i,d-i+1)}}$-module.}
\begin{equation}\label{expression}
U(L_{(i-1,0)},\ldots,L_{(i-1,i-2)})(K[{\textstyle \frac{X_{i}}{X_{i-1}},\ldots, \frac{X_d}{X_{i-1}}] }\cdot  V_{i,\lambda})=\bigoplus_{k_0,\ldots, k_{i-1} \leq 0  \atop {k_{i},\ldots, k_d \geq 0  \atop {k_0+\cdots + k_d=0}}}\!\!\!\!\! K\cdot X_0^{k_0}\cdots X_d^{k_d} \cdot  V_{i,\lambda} .
\end{equation}
Here $U(L_{(i-1,0)},\ldots,L_{(i-1,i-2)})$ denotes the subalgebra of $U(\frg)$ generated by $L_{(i-1,0)},\ldots,$ $L_{(i-1,i-2)}.$
Indeed, the above expression is contained in the image. As for the other inclusion, we
note that $$V_\lambda=U({\rm Lie}(R_u(L_{(1,d)} \cap B)))v_{\lambda},$$ where $R_u(L_{(1,d)} \cap B)$ denotes the unipotent radical of $L_{(1,d)} \cap B $, cf. \cite{Ja} p. 204.
Thus, if $L_{(k,l)} $ is a root of $\frg$ contained in ${\rm Lie}(R_u(L_{(1,d)} \cap B))$ then $k> l,k\geq 2$ and  $l\geq 1.$
Let $w_{i-1}\ast \lambda \prec w_i\ast \lambda.$ Then 
\begin{eqnarray*}
& & \frac{X_{i-1}^{i-1}}{X_0\cdots X_{i-2}}\cdot w_{i-1}(L_{(k,l)} v_\lambda)=
\frac{X_{i-1}^{i-1}}{X_0\cdots X_{i-2}}\cdot L_{(w_{i-1}(k),w_{i-1}(l))} w_{i-1}(v_\lambda)\\ & = &  L_{(w_{i-1}(k),w_{i-1}(l))}(v_{i-1,\lambda}) - \big(L_{(w_{i-1}(k),w_{i-1}(l))} \frac{X_{i-1}^{i-1}}{X_0\cdots X_{i-2}}\big)\cdot w_{i-1}(v_\lambda).
\end{eqnarray*}
Since $k>l\geq 1,$ we conclude that $w_{i-1}(k)>w_{i-1}(l).$ If $w_{i-1}(l) \not\in\{0,\ldots,i-2\}$ or $ w_{i-1}(k)\in\{0,\ldots,i-2\}$ we deduce that  $(L_{(w_{i-1}(k),w_{i-1}(l))} \frac{X_{i-1}^{i-1}}{X_0\cdots X_{i-2}})\cdot w_{i-1}(v_\lambda)=0.$ Otherwise, we get  $$\big(L_{(w_{i-1}(k),w_{i-1}(l))} \frac{X_{i-1}^{i-1}}{X_0\cdots X_{i-2}}\big)\cdot w_{i-1}(v_\lambda) = -\frac{X_{w_{i-1}(k)}}{X_{w_{i-1}(l)}}\cdot v_{i-1,\lambda}.$$ In any case, since $L_{(w_{i-1}(k),w_{i-1}(l))}(v_{i-1,\lambda})$ is contained in (\ref{expression}), we see that
$\frac{X_{i-1}^{i-1}}{X_0\cdots X_{i-2}}\cdot w_{i-1}(L_{(k,l)} v_\lambda)$ is contained in (\ref{expression}) as well.
The case $w_{i-1}\ast \lambda \succ w_i\ast \lambda$ is treated similarly by using identity (\ref{weights}). 
\qed

\begin{Example}
Let $\lambda=(0,\ldots,0)$. Then $\mu_{1,\lambda}=w_1\ast \lambda =(-1,1,0\ldots,0)$ and $V_{1,\lambda}=K\frac{X_1}{X_0}\oplus \cdots \oplus K\frac{X_d}{X_0}$. So $\tilde{H}_{\overline{X}_{w_1}}^1(\P^d_K,\cO)$ is a quotient of $K[\frac{X_1}{X_0},\ldots, \frac{X_d}{X_0}] \otimes V_{1,\lambda}$ with non-trivial kernel.
On the other hand, if $\lambda=(-1,1,0,\ldots,0)$ then  $\mu_{1,\lambda}=w_0\ast \lambda =\lambda$ and  $V_{1,\lambda}=V_\lambda = K d(\frac{X_1}{X_0})\oplus \cdots \oplus K d(\frac{X_d}{X_0}).$ In this situation the map $K[\frac{X_1}{X_0},\ldots, \frac{X_d}{X_0}] \otimes V_{1,\lambda}\rightarrow  \tilde{H}_{\overline{X}_{w_1}}^1(\P^d_K,\Omega^1)$ is an isomorphism. 
\end{Example}

We shall determine a ${\bf P^+_{(i,d-i+1)}}$-submodule of  (\ref{expression})  generating it as  $U(\frg)$-module.  We make use of the following statement which can be found
in various descriptions in \cite{FH}.

\begin{Lemma}\label{tensorproduct}
Let $\nu \in (\bbZ^n)_+=\{(\nu_1,\ldots,\nu_n) \mid \nu_1\geq \nu_2 \geq \cdots \geq \nu_n \}$ be a dominant weight, $n\in\N$. Let $V_\nu$ be the irreducible algebraic representation of ${\rm GL}_n$ 
of highest weight $\nu$. For $k\in \N,$ we consider the irreducible algebraic representation $V_k$ of highest weight $(k,0,\ldots,0)\in (\bbZ^n)_+,$ i.e., $V_k \cong {\rm Sym}^k(K^n).$
Then\footnote{Note that the condition $c_{i+1} \leq \nu_i- \nu_{i+1},\, i=1,\ldots,n-1$, implies that $\nu+(c_1,\ldots,c_n)\in (\bbZ^n)_+$} 
$$\label{} V_k \otimes V_\nu = \bigoplus_{(c_{1},\ldots,c_n)\in \N_0^n \atop c_{i+1} \leq \nu_i- \nu_{i+1}, \, i=1,\ldots,n-1}  V_{\nu+(c_1,\ldots,c_n)} .$$
\end{Lemma}

\proof By tensoring $V_\nu$ with $\det^{-\nu_n}$, we may 
suppose that $\nu_n=0$. Set $a_i:= \nu_i - \nu_{i+1},\; i=1,\ldots,n-1.$ Then we have
$V_\nu= \Gamma_{a_1,\ldots,a_{n-1}}$ and $V_k=\Gamma_{k,0\ldots,0}$ with the terminology in \cite{FH} \S 15.
By loc.cit. Prop. 15.25 a) we deduce the decomposition  (as ${\rm SL}_n$-modules)
$$  V_k \otimes V_\nu  = \bigoplus_{b_1,\ldots,b_{n-1}}  \Gamma_{b_1,\ldots,b_{n-1}},$$
where the sum is over all  non-negative integers $b_1,\ldots,b_{n-1}$ for which there are non-negative integers $c_1,\ldots,c_n$ with
$\sum_i c_i =k$, with $c_{i+1} \leq a_i$  and with $b_i=a_i+c_i-c_{i+1}$ for $1\leq i \leq n-1$.
The highest weight of $\Gamma_{b_1,\ldots,b_{n-1}}$ is $\nu+(c_1-c_n,c_2-c_n,\ldots,c_{n-1}-c_n,0).$
Since we deal with ${\rm GL}_n$-modules of total weight $\nu_1+\cdots + \nu_n + k$, we have to replace
$V_{\nu+(c_1-c_n,c_2-c_n,\ldots,c_{n-1}-c_n,0)}$ by $V_{\nu +(c_1,\ldots,c_n)}$. 
\qed

We start to investigate the extreme cases $i=1$ and $i=d.$
For $k\geq 0$, let $K[\frac{X_1}{X_0},\cdots,\frac{X_d}{X_0}]_k$ be the set of polynomials of degree $k$ in the indetermines $\frac{X_1}{X_0}, \ldots , \frac{X_d}{X_0}$. Then we get identifications
$${\textstyle K[\frac{X_1}{X_0},\cdots,\frac{X_d}{X_0}]_k\cong{\rm Sym}^k(K\frac{X_1}{X_0}\oplus \cdots \oplus K \frac{X_d}{X_0}) \cong V_{(-k,k,0\ldots,0)} }.$$
For an integral vector $\nu=(\nu_1,\ldots, \nu_n)\in (\bbZ^n)_+$, we set $$|\nu|=\nu_1-\nu_n.$$

\begin{Lemma}
(i) Let $\mu_{1,\lambda}=(\mu_0,\mu')$ with $\mu_0\in \bbZ$ and $\mu'\in (\bbZ^d)_+.$ Then  $\tilde{H}^{1}_{\overline{X}_{w_1}}(\P^d_K,\cF_\lambda)$ $= K[\frac{X_1}{X_0},\cdots,\frac{X_d}{X_0}] \cdot V_{1,\lambda}$ is generated as $U(\frg)$-module by 
$\oplus_{k \leq |\mu'|}  K[\frac{X_1}{X_0},\cdots,\frac{X_d}{X_0}]_k  \cdot V_{1,\lambda}.$

(ii)  The module   $\tilde{H}^{d}_{\overline{X}_{w_d}}(\P^d_K,\cF_\lambda)= K[\frac{X_d}{X_0},\cdots,\frac{X_d}{X_{d-1}}] \cdot V_{d,\lambda}$ is generated as $U(\frg)$-module by  $V_{d,\lambda}.$ 
\end{Lemma}

\proof  (i) We identify
$K[\frac{X_1}{X_0},\ldots,\frac{X_d}{X_0}] \otimes V_{1,\lambda}$ with $\bigoplus_{k\geq 0} V_{(-k,k,0\ldots,0)} \otimes V_{1,\lambda}  .$
We may apply the previous lemma  with $n=d$ and $\nu = \mu '$. We deduce that the number of irreducible summands in $V_{(-k,k,0\ldots,0)} \otimes V_{1,\lambda}$  is the same for $k \geq a_1+a_2+\cdots + a_{n-1}$ with $a_i=\mu'_i-\mu'_{i+1}.$ 
But the latter sum is exactly $\mu'_1-\mu'_d=|\mu'|.$ 
The Lie algebra $\frg$ maps $V_{(-k,k,0\ldots,0)}\cdot V_{1,\lambda}$ to $V_{(-k-1,k+1,0\ldots,0)}\cdot V_{1,\lambda}$
and its image is again a ${\bf P^+_{(1,d)}}$-module.  
The claim follows since   irreducible submodules are mapped to different irreducible submodules by weight reasons.

(ii) By (\ref{expression}) it is enough to show that $K[\frac{X_d}{X_{d-1}}]\cdot V_{d,\lambda} \subset U(\frg)\cdot V_{d,\lambda}.$ We consider the following two cases.

\noindent Case 1: $w_d\ast \lambda \preccurlyeq  w_{d-1}\ast \lambda$. Then $v_{d,\lambda}$ is a highest weight vector of $V_{d,\lambda}$. We compute 
\begin{eqnarray*}
\frac{X_d}{X_i}\cdot  v_{d,\lambda} & = & \frac{X_d}{X_i}\cdot \frac{X_d}{X_0\cdots X_{d-1}}\cdot  w_d(v_\lambda)
= - (L_{(d,i)}\frac{X_d}{X_0\cdots X_{d-1}}) \cdot  w_d(v_\lambda) \\ & = & - L_{(d,i)}(\frac{X_d}{X_0\cdots X_{d-1}} \cdot w_d(v_\lambda)) + \frac{X_d}{X_0\cdots X_{d-1}} \cdot L_{(d,i)}  w_d(v_\lambda).
\end{eqnarray*}
But $L_{(d,i)} w_d(v_\lambda)= w_d(L_{(0,i+1)}v_\lambda)=0$ since $L_{(0,i+1)}v_\lambda=0.$ 
Thus we obtain $\frac{X_d}{X_i}\cdot v_{d,\lambda} = -L_{(d,i)} v_{d,\lambda} .$
On the other hand, the module $V_{d,\lambda}$ is equal to $U({\rm Lie}(R_u(L_{(d,1)} \cap B)))v_{d,\lambda}$.
If $L_{(k,l)}$ is a root contained in ${\rm Lie}(R_u(L_{(d,1)} \cap B))$ we necessarily have $l<k<d.$ We deduce that
\begin{eqnarray*}
\frac{X_d}{X_{d-1}}\cdot L_{(k,l)} v_{d,\lambda} & = & L_{(k,l)}(\frac{X_d}{X_{d-1}}\cdot v_{d,\lambda})= - L_{(k,l)} (L_{(d,d-1)}v_{d,\lambda}) \\ 
& = &  - L_{(d,d-1)} (L_{(k,l)} v_{d,\lambda})
+[L_{(k,l)},L_{(d,i)}]\cdot v_{d,\lambda} \\
& = &  - L_{(d,d-1)} (L_{(k,l)} v_{d,\lambda}) +
\delta_{i,k} L_{(d,l)}\cdot v_{d,\lambda}
\end{eqnarray*}
is contained in  $U(\frg)\cdot V_{d,\lambda}.$ The case of polynomials of higher degree is treated in the same way.

\noindent Case 2: $w_d\ast \lambda \succ w_{d-1}\ast \lambda$. Then $\delta_{d-1}(v_{d-1,\lambda})$ is a highest weight vector of $V_{d,\lambda}$.
We get $\delta_{d-1}(v_{d-1,\lambda})=(\frac{X_d}{X_{d-1}})^n \cdot v_{d,\lambda} $ for some $n>0,$ cf. (\ref{weights}).
We compute 
\begin{eqnarray*}
& & \frac{X_d}{X_i}\cdot \delta_{d-1}(v_{d-1,\lambda}) = \frac{X_d}{X_i}\cdot (\frac{X_d}{X_{d-1}})^n \cdot v_{d, \lambda} \\
& = & (\frac{X_d}{X_{d-1}})^n\cdot \frac{X_d}{X_i} \cdot v_{d,\lambda} = - (\frac{X_d}{X_{d-1}})^n \cdot L_{(d,i)} v_{d,\lambda} 
 \mbox{ (by case 1) } \\
& = & -L_{(d,i)}((\frac{X_d}{X_{d-1}})^n \cdot v_{d,\lambda}) + (L_{(d,i)}(\frac{X_d}{X_{d-1}})^n)  \cdot v_{d,\lambda} \\
& = &
-L_{(d,i)}(\delta_{d-1}(v_{d-1,\lambda})) + (L_{(d,i)}(\frac{X_d}{X_{d-1}})^n)  \cdot v_{d,\lambda}.
\end{eqnarray*}
But $L_{(d,i)}(\frac{X_d}{X_{d-1}})^n \cdot v_{d,\lambda}= 0$ if $i\neq d-1$. If $i=d-1$ then 
\begin{eqnarray*}
L_{(d,i)}(\frac{X_d}{X_{d-1}})^n  \cdot v_{d,\lambda} & = & 
-n (\frac{X_d}{X_{d-1}})^{n+1} \cdot v_{d,\lambda}  \\ & = & -n \frac{X_d}{X_i}\cdot \delta_{d-1}(v_{d-1,\lambda}),
\end{eqnarray*} so $(n+1)\frac{X_d}{X_i} \cdot \delta_{d-1}(v_{d-1,\lambda}) =  - L_{(d,i)} \delta_{d-1}(v_{d-1,\lambda}).$
If $w\in V_{d,\lambda}$ is an arbitrary vector we argue as in case 1.
\qed

Now we treat the general case which is a mixture  between the above extreme cases. Note that 
$K[X_{(m,n)}\mid m\geq i,\; n\leq i-1]=\bigoplus_{k_0,\ldots, k_{i-1} \leq 0  \atop {k_{i},\ldots, k_d \geq 0  \atop {k_0+\cdots + k_d=0}}} K\cdot X_0^{k_0}X_1^{k_1}\cdots X_d^{k_d}.$ 

\begin{Lemma}\label{Lemma_boxtimes} Write $\mu_{i,\lambda}=(\mu',\mu'')$ with $\mu'\in (\bbZ^i)_+$ and $\mu''\in (\bbZ^{d+1-i})_+.$ Then 
the ${\bf P^+_{(i,d-i+1)}}\ltimes U(\frg)$-module $\tilde{H}^{i}_{\overline{X}_{w_i}}(\P^d_K,\cF_\lambda)= K[X_{(m,n)}\mid m\geq i,\; n\leq i-1] \cdot  V_{i,\lambda}$ is generated by 
the ${\bf P^+_{(i,d-i+1)}}$-submodule
$\oplus_{k\leq |\mu''|} K[X_{(m,n)}\mid m\geq i,\; n\leq i-1]_k \cdot  V_{i,\lambda}.$  
\end{Lemma}

\proof 
Write $V_{i,\lambda}=V_{\mu'} \boxtimes V_{\mu''}.$ We may identify $K[X_{(m,n)}\mid m\geq i,\; n\leq i-1]_k$ with
the outer tensor product representation  $V_{(0,\ldots,0,-k)} \boxtimes V_{(k,0\ldots,0)}$ of ${\bf L^+_{(i,d-i+1)}}.$  We get
$$K[X_{(m,n)}\mid m\geq i,\; n\leq i-1]_k \otimes  V_{i,\lambda} \cong (V_{(0,\ldots,0,-k)} \otimes V_{\mu'}) \boxtimes (V_{(k,0\ldots,0)} \otimes V_{\mu''}). $$
By the proof of Proposition \ref{Propositionweight} we saw that $K[X_{(k,l)}\mid k\geq i,\; l\leq i-1] \cdot V_{i,\lambda}$  is already  generated as $U(\frg)$-module by $K[\frac{X_{i}}{X_{i-1}},\ldots, \frac{X_d}{X_{i-1}}] \cdot V_{i,\lambda}.$   Now we apply the proof of part (i) of the previous lemma  to deduce the claim. \qed

For 
$\mu_{i,\lambda}=(\mu',\mu'')$  with $\mu'\in (\bbZ^i)_+$ and $\mu''\in (\bbZ^{d-i+1})_+ ,$
we define
\begin{eqnarray*}
\Phi_{i,\lambda} =  \bigcup\limits_{k=0}^{|\mu''|} &\Big\{&   \big(\mu'-(d_i,\ldots,d_1),\mu''+(c_1,\ldots, c_{d-i+1})\big) \mid {\textstyle \sum}_j c_j= \textstyle{\sum}_j  d_{j}=k, c_1=0 \\ & & \mbox{ or } d_1=0,\;  c_{j+1}\leq \mu_j'' -\mu_{j+1}'', \;
j=1,\ldots,d-i,\; d_{j+1} \leq \mu'_{i-j} - \mu'_{i-j+1}, \,\\ & &  j=1,\ldots,i-1 \, \Big\}.
\end{eqnarray*}
We let
$$M_{i,\lambda}:= \bigoplus_{\mu\in \Phi_{i,\lambda}} V_{\mu} \;\;\;\;\; \subset \;\;\;\;\;  K[X_{(m,n)}\mid m\geq i,\; n\leq i-1] \otimes  V_{i,\lambda} $$
be the sum of the irreducible  ${\bf P^+_{(i,d-i+1)}}$-modules with respect to the weights appearing in $\Phi_{i,\lambda}.$
Let 
$$p_i: K[X_{(m,n)}\mid m\geq i,\; n\leq i-1] \otimes  V_{i,\lambda} \rightarrow \tilde{H}^{i}_{\overline{X}_{w_i}}(\P^d_K,\cF_\lambda)$$ be the quotient map.

\begin{Corollary}
Then for $1\leq i \leq d,$ the ${\bf P^+_{(i,d-i+1)}} \ltimes U(\frg)$-module
$\tilde{H}^{i}_{\overline{X}_{w_i}}(\P^d_K,\cF_\lambda)$ is generated as $U(\frg)$-module by $p_i(M_{i,\lambda}).$
\end{Corollary}

\proof We write $V_{(0,\ldots,0,-k)} \otimes V_{\mu'}= (V_{(0,\ldots,0,-k)} \otimes V_{\mu'})^{\ast \ast} = (V_{(k,\ldots,0,0)} \otimes V_{(\mu')^\ast})^\ast$ where $(\mu')^\ast=(-\mu_i,\ldots,-\mu_1)$ and $^\ast$ indicates the dual representation, cf. \cite{FH} Ex. 15.50.
By Lemma \ref{tensorproduct} we have a decomposition 
$$V_{(k,\ldots,0,0)} \otimes V_{(\mu')^\ast}= \bigoplus_{(d_{1},\ldots,d_i)\in \N_0^i \atop d_{j+1} \leq (\mu')^\ast_j- (\mu')^\ast_{j+1}, \, j=2,\ldots,i } \!\!\!\!V_{(\mu')^\ast+(d_1,\ldots,d_i)}.$$ 
Thus we get 
$$V_{(0,\ldots,0,-k)} \otimes V_{\mu'}= \big(\!\!\!\!\!\!\!\!\bigoplus_{(d_{1},\ldots,d_i)\in \N_0^i \atop d_{j+1} \leq (\mu')^\ast_j- (\mu')^\ast_{j+1}, \, j=2,\ldots,i} \!\!\!\!\!\!\!\! \!\!\!\!V_{(\mu')^\ast+(d_1,\ldots,d_i)}\big)^\ast=\bigoplus_{(d_{1},\ldots,d_i)\in \N_0^i \atop d_{j+1} \leq (\mu')^\ast_j- (\mu')^\ast_{j+1}, \, j=2,\ldots,i}  \!\!\!\!\!\!\!\!\!\!\!\!V_{\mu'-(d_i,\ldots,d_1)}$$ But $(\mu')^\ast_j - (\mu')^\ast_{j+1}=
 \mu'_{i-j} - \mu'_{i-j+1}$. 
Finally, if $c_1 >0$ and $d_1>0$ then 
$$p_i(V_{\mu'-(d_i,\ldots,d_1)}\boxtimes V_{\mu''+(c_1,\ldots,c_{d-i+1})})= \frg \cdot p_i( V_{\mu'-(d_i,\ldots,d_1-1)}\boxtimes V_{\mu''+(c_i-1,\ldots,c_{d-1+1})}).$$
The claim follows.
\qed

\begin{Remark}
We point out that for some weights $\lambda$, some of the irreducible submodules $V_\mu \subset M_{i,\lambda}, \mu \in \Phi_{i,\lambda},$ apart from $\mu_{i,\lambda}$   are mapped to zero under the quotient map $p_i$. We refer to section 3 for examples.
\end{Remark}

Now we translate the above result for the computation of $\tilde{H}^i_{\P^{d-i}_K}(\P^d_K,\cF).$
Consider the block matrix
$$ z_i:=\left( \begin{array}{cc} 0  & I_i \\  I_{d+1-i} & 0 \end{array} \right) \in G,$$
where $I_j\in {\rm GL}_{j}$ denotes the $j\times j$-identity matrix. Then $V(X_0,\ldots X_{i-1})=\overline{X}_{w_i}$ is transformed into $V(X_{d-i+1},\ldots ,X_d)=\P^{d-i}_K$
under the action of $z_i.$ We have $$z_i \cdot {\bf P_{(d-i+1,i)} } \cdot z_i^{-1} = {\bf P^+_{(i,d+1-i)}}$$ and on the Levi subgroups the conjugacy map
is given by
$${\bf L_{(d-i+1,i)}} \ni\left( \begin{array}{cc} A  & 0 \\  0 & B \end{array} \right) \mapsto \left( \begin{array}{cc} B & 0 \\  0 & A \end{array} \right) \in {\bf L_{(i,d-i+1)}}.$$  Thus we get an isomorphism
$$\tilde{H}^{i}_{\P^{d-i}_K}(\P^d_K,\cF) \stackrel{\sim}{\rightarrow} \tilde{H}^{i}_{\overline{X}_{w_i}}(\P^d_K,\cF)$$
compatible with the action of the parabolic subgroups.
Hence in order to determine the  ${\bf P_{(d-i+1,i)}}\ltimes U(\frg)$-representation of  $\tilde{H}^{i}_{\P^{d-i}_K}(\P^d_K,\cF)$
in terms of highest weight vectors, we have to apply $z_i$ - regarded as an element in the Weyl group W - to them.
Clearly the dominant weights  are respected by this transformation. We set
\begin{eqnarray*}
\Psi_{i,\lambda}  & = & z_i^{-1}\cdot \Phi_{i,\lambda} \\ & = & \bigcup\limits_{k=0}^{|\mu''|} \Big\{   (\mu''+(c_1,\ldots, c_{d-i+1}),\mu'-(d_i,\ldots,d_1)) \mid \textstyle{\sum}_j c_j=\textstyle{\sum}_j d_{j}=k, c_1=0 \\ & & \mbox{ or } d_1=0,\;  c_{j+1}\leq \mu_i'' -\mu_{j+1}'', \; 
j=1,\ldots,d-i,\; d_{j+1} \leq \mu'_{i-j} - \mu'_{i-j+1}, \,\\ & &  j=1,\ldots,i-1 \, \Big\}
\end{eqnarray*} 
and 
$$N_{i,\lambda}=  \bigoplus_{\mu\in \Psi_{i,\lambda}} V_\mu \;\;\;\;\; \subset \;\;\;\;\;   K[X_{(m,n)}\mid m\leq d-i,\; n\geq d-i+1] \otimes  V_{z_i^{-1}\mu_{i,\lambda}}.$$
Let 
$$q_i: K[X_{(m,n)}\mid m\leq d-i,\; n\geq d-i+1] \otimes  V_{z_i^{-1}\mu_{i,\lambda}}  \rightarrow \tilde{H}^{i}_{\P^{d-i}_K}(\P^d_K,\cF_\lambda)$$ be the quotient map. We obtain:

\begin{Corollary}\label{irreduciblelagbraicrepr} 
For $1\leq i \leq d,$ the ${\bf P_{(d-i+1,i)}} \ltimes U(\frg)$-module
$\tilde{H}^{i}_{\P^{d-i}_K}(\P^d_K,\cF_\lambda)$ is generated by $q_i(N_{i,\lambda}).$
\end{Corollary}

\section{The $G$-representation $H^0(\cX,\cF)$}
\subsection{The fundamental complex }
In this section we recall the construction \cite{O3} of an acyclic resolution of the constant sheaf $\Z$ on the boundary of $\cX$ considered as an object in the category of pseudo-adic spaces \cite{H}.

In order to determine the structure of $\cF(\cX)'=H^0(\cX,\cF)'$ as a locally analytic $G$-representation, we
proceed as follows. Let $${\cY}=(\P^d_K)^{rig} \setminus \cX  $$ be the set-theoretical complement of $\cX.$ Consider the  topological exact sequence
of locally convex $K$-vector spaces with continuous $G$-action
$$0\rightarrow H^0(\P^d_K,\cF) \rightarrow H^0(\cX,\cF)\rightarrow H^1_{\cY}(\P^d_K,\cF) \rightarrow H^1(\P^d_K,\cF)\rightarrow 0.$$
Note that the higher cohomology groups $H^i(\cX,\cF), \; i>0,$ vanish since $\cX$ is a Stein space \cite{K2}.
The $G$-representations $H^0(\P^d_K,\cF),\;H^1(\P^d_K,\cF)$ are finite-dimensional algebraic. A more delicate problem
is to understand the structure of the $G$-representation $H^1_{\cY}(\P^d_K,\cF)$ which is a $K$-Fr\'echet space. More precisely, it is by Proposition \ref{projlimI} and Proposition \ref{Steinspace} a projective limit of $K$-Banach spaces 
$$H^1_{\cY}(\P^d_K,\cF)= \varprojlim\nolimits_n H^1_{{\cY}_n^-}(\P^d_K,\cF). $$
In \cite{O3} we constructed
 acyclic resolutions of overconvergent \'etale sheaves on the boundary of period domains. We want to apply this construction to our situation. The construction makes  use of Huber's adic spaces \cite{H}. In the appendix we give an alternative approach
avoiding these spaces.
In the following, the symbol $X^{ad}$ indicates the adic space attached to a scheme $X$ or to a rigid analytic variety $X$ defined over $K.$

We take the complement of $\cX$ in the category  of adic spaces, i.e., we  set
$${\cY}^{ad}:=(\P^d_K)^{ad}\setminus \cX^{ad}.$$
This is a closed pseudo-adic subspace of $\P^d_K.$
Let $\{e_0,\ldots,e_d\}$ be the standard basis of $V=K^{d+1}.$
For any $\alpha_i\in \Delta,$ put
$$V_i=\bigoplus^i_{j=0} K\cdot e_j \; \mbox{ and } \; Y_i=\P(V_i)$$
For any subset $I\subset \Delta$ with  $\Delta\setminus I=\{\alpha_{i_1} < \ldots < \alpha_{i_r}\},$ let $Y_I$ be the closed $K$-subvariety of $\P^d_K$ defined by
$$Y_I=\P(V_{i_1}).$$
Furthermore, let $P_I$ be the lower parabolic subgroup of $G,$ such that  $I$
coincides with the simple roots appearing in the Levi factor of $P_I$. Hence the group $P_I$ stabilizes $Y_I.$
We obtain
\begin{eqnarray}\label{Yad=}
{\cY}^{ad}=\bigcup_{I \subset \Delta} \bigcup_{g\in G/P_I}g\cdot Y_I^{ad} = \bigcup_{g\in G}g\cdot Y_{ \Delta \setminus \{\alpha_{d-1}\}}^{ad}.
\end{eqnarray}

\noindent For any compact open subset $ W\subset G/P_I,$ put
$$Z_I^W:= \bigcup_{g\in W} gY_I^{ad}.$$
We  proved in \cite{O3}, Lemma 3.2,
that $Z_I^W$ is a closed pseudo-adic subspace of $(\P^d_K)^{ad}.$
By (\ref{Yad=})  it follows that
$${\cY}^{ad}=\bigcup_{I\subset \Delta \atop |\Delta\setminus I|=1}Z_I^{G/P_I}=Z_{\Delta\setminus \{\alpha_{d-1}\}}^{G/P_{\Delta\setminus{\{\alpha_{d-1} \}}}}.$$

Starting from the constant \'etale sheaf $\Z$ on ${\cY}^{ad}$ we constructed a sheaf of locally constant sections on the same space. We recall  the definition.
Consider the natural closed embeddings of pseudo-adic spaces
$$\Phi_{g,I} :gY_I^{ad} \longrightarrow {\cY}^{ad}$$
resp.
$$\tilde{\Phi}_{g,I,W} :gY_I^{ad} \longrightarrow Z_I^W$$
resp.
$$\Psi_{I,W} :Z^W_I \longrightarrow {\cY}^{ad}.$$
Put
$$ \Z_{g,I}:=(\Phi_{g,I})_\ast(\Phi_{g,I}^\ast\, \Z)$$
resp.
$$ \Z_{Z_I^W}:=(\Psi_{I,W})_\ast(\Psi_{I,W}^\ast\, \Z)$$
and let
$$\tilde{\Phi}_{g,I,W}^{\#} :\Z_{Z_I^W} \longrightarrow \Z_{g,I}$$
be the natural homomorphism given by restriction.
Let
${\mathcal C}_I$ be the category of compact open disjoint coverings of $G/P_I$
where the morphisms are given by the refinement-order.
For a covering $c =(W_j)_{j} \in {\mathcal C}_I,$ we denote by $\Z_c$ the sheaf  on $\cY^{ad}$ defined by

\noindent $\begin{array}{lcl}\Z_c(U) & := & \Big\{ (s_g)_g \in
\displaystyle\prod_{g\in G/P_I} \Z_{g,I}(U) \mid \;\mbox{ there are
sections }\; s_j\; \in \Z_{Z_I^{W_j}}(U),  \mbox{  such }
\\ & & \\ & &
 \mbox{ that }\tilde{\Phi}^{\#}_{g,I,W_j}(s_j) = s_g \; \mbox{for all } g \in W_j \Big\}.
\end{array} $

\noindent Note that
$\Z_c$ is just the image of the natural morphism of sheaves
$$\bigoplus_{j \in A}\; \Z_{Z_I^{W_j}} \hookrightarrow \prod_{g\in G/P_I} \Z_{g,I}.$$
We put
\begin{eqnarray}\label{prod'=limInd}
\sideset{}{'}\prod_{g\in G/P_I} \Z_{g,I} = \varinjlim_{c\in {\mathcal C}_I} \Z_c,
\end{eqnarray}
We obtain the following  complex of sheaves on ${\cY}^{ad},$

\begin{eqnarray}
\nonumber 0 \rightarrow \Z \rightarrow\!\!\! \bigoplus_{I \subset
\Delta \atop |\Delta\setminus I|=1} \sideset{}{'}\prod_{g\in G/P_I} \Z_{g,I}
\rightarrow \!\!\!\bigoplus_{I \subset \Delta \atop |\Delta\setminus I|=2}
\sideset{}{'}\prod_{g\in G/P_I} \Z_{g,I} \rightarrow \cdots \rightarrow \!\!\!\bigoplus_{I \subset \Delta \atop |\Delta\setminus I|=i}
\sideset{}{'}\prod_{g\in G/P_I} \Z_{g,I} \rightarrow \cdots \\ \label{complex} \\ \nonumber \cdots\rightarrow
\!\!\!\bigoplus_{I \subset \Delta \atop |\Delta\setminus I|=d-1} \sideset{}{'}\prod_{g\in G/P_I} \Z_{g,I}
\rightarrow \sideset{}{'}\prod_{g\in G/P_\emptyset} \Z_{g,\emptyset} \rightarrow 0.
\end{eqnarray}

\begin{Theorem}\label{Theorem1}
The complex (\ref{complex})  is acyclic.
\end{Theorem}

\proof This is Theorem 3.3 in \cite{O3}. Strictly speaking, we treated in loc.cit. the case of the constant \'etale sheaf $\Z/n\Z.$
But the proof  is the same.    \qed

\subsection{Evaluation of the spectral sequence}\label{Berechnung}

In this section we evaluate the spectral sequence which is induced
by the complex (\ref{complex}) applied to ${\rm
Ext}^\ast(i_\ast(-),\; \cF).$ Here $i: \cY^{ad} \hookrightarrow
(\P^d_K)^{ad}$ denotes the closed embedding.

By \cite{SGA2} Proposition 2.3
bis., we conclude that
$${\rm Ext}^\ast(i_\ast(\Z_{{\cY}^{ad}}),\cF) =
H^\ast_{{\cY}^{ad}}(\P^d_K,\cF),$$  Further, we have $H^\ast_{{\cY}^{ad}}(\P^d_K,\cF) =
H^\ast_{{\cY}}(\P^d_K,\cF)$ since the topoi of $\cX$ and $\cX^{ad}$
are equivalent, cf. \cite{H}, Prop. 2.1.4.
Recall that $G_0={\bf G}(O_K)$ denotes the compact $p$-adic group of $O_K$-valued points of ${\bf G}$.

\begin{Proposition}\label{Ext=projlim} For all $I\subset \Delta,$ there is an isomorphism
$${\rm Ext}^\ast(i_\ast(\sideset{}{'}\prod\limits_{g\in G/P_I} \Z_{g,I}), \cF)
= \varprojlim_{n\in \N}\bigoplus_{g\in G_0/P_I^n} H_{gY_I(\epsilon_n)}^{\ast}(\P^d_K, \cF).$$
\end{Proposition}

\proof Consider the  family
$$\Big\{gP_I^n\mid g\in G_0, n\in \N\Big\}$$
of compact open subsets in $G/P_I$ which yields cofinal coverings in $\cC_I.$
We obtain by (\ref{prod'=limInd}) the identity
$$\sideset{}{'}\prod\limits_{g\in G/P_I} \Z_{g,I}
=\varinjlim_{c\in {\mathcal C}_I} \Z_c = \varinjlim_{n\in \N} \bigoplus_{g\in G_0/P^n_I} \Z_{Z^{gP^n_I}_I}.$$
Choose an injective resolution $\cI^\bullet$ of $\cF$. 
We get
\begin{eqnarray*}
{\rm Ext}^i(i_\ast(\sideset{}{'}\prod\limits_{g\in G/P_I} \Z_{g,I}),\cF) &=& H^i({\rm Hom}(i_\ast(\sideset{}{'}\prod\limits_{g\in G/P_I} \Z_{g,I}), \cI^\bullet)) \\
= H^i({\rm Hom}(\varinjlim_{n\in \N} \bigoplus_{g\in G_0/P^n_I} i_\ast(\Z_{Z^{gP^n_I}_I}), \cI^\bullet)) &  = & H^i(\varprojlim_{n\in \N} \bigoplus_{g\in G_0/P^n_I} {\rm Hom}(i_\ast(\Z_{Z^{gP^n_I}_I}), \cI^\bullet)) \\
= H^i(\varprojlim_{n\in \N} \bigoplus_{g\in G_0/P^n_I} H^0_{Z^{gP^n_I}_I}(\P^d_K,\cI^\bullet)) &  &  .
\end{eqnarray*}
We make use of the following lemma. Here $\varprojlim\nolimits_{n\in \N}^{(r)}$ is the r-th right derived functor of $\varprojlim\nolimits_{n\in \N}$.
\begin{Lemma}
Let  $\cI$  be an injective sheaf on $(\P^d_K)^{ad}.$ Then 
$$\varprojlim\nolimits_{n\in \N}^{(r)} \bigoplus_{g\in G_0/P^n_I} H^0_{Z^{gP^n_I}_I}(\P^d_K,\cI)=0 \mbox{ for } r\geq 1. $$
\end{Lemma}
\begin{proof}
It suffices to show that the projective  systems
$$\Big(\bigoplus_{g\in G_0/P^n_I} H^0(\P^d_K,\cI)\Big)_{n\in \N} $$
and
$$\Big(\bigoplus_{g\in G_0/P^n_I} H^0((\P^d_K)^{ad}\setminus {Z^{gP^n_I}_I},\cI)\Big)_{n\in \N} $$
are $\varprojlim\nolimits_{n\in \N}$-acyclic. Clearly  the maps $G_0/P^n_I \rightarrow G_0/P^m_I$ are surjective for $n\geq m$. It follows that the first projective system is $\varprojlim\nolimits_{n\in \N}$-acyclic. Since $\cI$ is injective, we have surjections
$$H^0((\P^d_K)^{ad}\setminus {Z^{hP^n_I}_I},\cI) \rightarrow H^0((\P^d_K)^{ad}\setminus {Z^{gP^m_I}_I},\cI) $$
for $n\geq m$ and $hP^n_I \subset gP^m_I.$
Thus we see that the transition maps $$\bigoplus_{h\in G_0/P^n_I} H^0((\P^d_K)^{ad}\setminus {Z^{hP^n_I}_I},\cI) \rightarrow \bigoplus_{g\in G_0/P^m_I} H^0((\P^d_K)^{ad}\setminus {Z^{gP^m_I}_I},\cI)$$ are surjective, as well. The claim follows.
\end{proof}
\noindent Thus we get by applying  a spectral sequence argument (note that $\varprojlim^{(r)}=0 $ for $r\geq 2$ \cite{Je}) short exact sequences, $i\in \N,$
$$0\rightarrow {\varprojlim_{n}}^{(1)}\!\!\!\!\! \bigoplus_{g\in G_0/P^n_I} \!\!\!\!\!  H^{i-1}_{Z^{gP^n_I}_I}(\P^d_K,\cF) \rightarrow  {\rm Ext}^i(i_\ast(\sideset{}{'}\prod\limits_{g\in G/P_I} \Z_{g,I}),\cF) \!\rightarrow \varprojlim_{n} \!\!\bigoplus_{g\in G_0/P^n_I} \!\!\!\!\!  H^i_{Z^{gP^n_I}_I}(\P^d_K,\cF) \rightarrow 0 .$$

\begin{Lemma} The projective system
$\Big(\bigoplus_{g\in G_0/P^n_I} H^{i-1}_{Z^{gP^n_I}_I}(\P^d_K,\cF)\Big)_{n\in \N}$ consists of $K$-Fr\'echet spaces and satisfies the
(topological) Mittag-Leffler property for all $i\geq 1$ (cf. \cite{EGAIII} 13.2.4).
\end{Lemma}
\begin{proof}
By the same methods as those used in Propositions \ref{Drinfeld} and \ref{Steinspace}, we can choose a decreasing sequence of admissible open subsets $$\cdots \supset D_{m-1} \supset D_m \supset D_{m+1}\supset \cdots$$
in $(\P^d_K)^{rig}$ with  $$\bigcap \nolimits_{m} (D_m)^{ad} = Z_I^{P_I^n},$$
and such that the complements $(\P^d_K)^{rig}\setminus D_m$ are admissible open.
In fact, we can choose these subsets to be coverings of the shape
$$D_m= \bigcup_{h\in R_m} h\cdot Y_I^-(\epsilon_m), $$ where
$R_m \subset P_I^n$ are  finite subsets. Here we may assume that $R_m\subset R_{m+1}$  and $1\in R_m$ for all $m.$ 
By translation with $g\in  G_0$ we obtain admissible open subsets $g\cdot D_m$ of $(\P^d_K)^{rig}$  with
\begin{equation}\label{Translation}
\bigcap\nolimits_{m} g\cdot (D_m)^{ad} =  Z_I^{gP_I^n}.
\end{equation}
We shall see that the cohomology groups  $H^\ast((\P^d_K)^{rig}\setminus  D_m,\cF)$ are $K$-Fr\'echet spaces and that the  transition maps
$$H^\ast((\P^d_K)^{rig}\setminus  D_{m+1},\cF) \rightarrow  H^\ast((\P^d_K)^{rig}\setminus  D_m,\cF)$$ have dense image.
Let $\Delta\setminus I=\{\alpha_{i_1}< \cdots < \alpha_{i_r} \}.$ If $\Delta\setminus I=\{\alpha_{{d-1}}\}$, i.e., $Y_I\subset \P^d_K$ is a hyperplane,  then the covering $((\P^d_K)^{rig}\setminus D_m)_m$ is of the type  considered in Proposition \ref{Steinspace} and the statement is a priori clear since it holds for Stein spaces.
In general, we may write $$Y_I^-(\epsilon_m)= \bigcap_{j>i_1} H_j^-(\epsilon_m)$$ where $H_j$ is the hyperplane $V(X_j)\subset \P^d_K.$
Then 
\begin{equation}\label{set}
(\P^d_K)^{rig} \setminus D_m = \bigcap_{h\in R_m} h\cdot \Big(\bigcup_{j>i_1} (\P^d_K)^{rig}\setminus H_j^-(\epsilon_m)\Big).
\end{equation}
For a hyperplane $H\subset \P^d_K$, let $\ell_{H}\in S$ be a unimodular linear  polynomial with $V(\ell_H)=H.$
Thus a point $z\in (\P^d_K)^{rig}$ is contained in (\ref{set}) if for all $h\in R_m$, there is an index $j>i_1$ with $|\ell_{h\cdot H_j}(z)|\geq \epsilon_m.$
For each $h\in R_m$, let $j_h>i_1$ be some integer. Set 
$$U_{(j_h)_{h}}=\{ z\in (\P^d_K)^{rig} \mid |\ell_{h\cdot H_{j_h}}(z)|\geq \epsilon_m \;\; \forall h\in R_m\}.$$
Then for varying $(j_h)_h$, the sets $U_{(j_h)_h}$ form  an open  covering of $(\P^d_K)^{rig}\setminus D_m$ consisting of $K$-affinoid subsets. In fact, let
$$\cH_{(j_h)_h}=\{h\cdot H_{j_h}\mid h\in R_m\}.$$ By the  
same reasoning as in \cite{SS} Prop. 4, we see that the $K$-algebra $\cO(U_{(j_h)_h})$ of analytic functions on $U_{(j_h)_h}$ is isomorphic to the $K$-affinoid algebra 
$$K\langle T_{H_i,H}, T_{H,H'} \mid 0\leq i \leq d,\, H,H' \in \cH_{(j_h)_h}\}  \rangle / I_m,$$ 
where $I_m$ is the closed ideal generated by the elements 
\begin{eqnarray*}
& & T_{H,H}-\pi^m,\;  H\in \cH_{(j_h)_h} \\
& &  T_{H,H'}\cdot  T_{H',H''} - \pi^m \cdot T_{H,H''},\; H',H''\in \cH_{(j_h)_h}, H\in \cH_{(j_h)_h} \cup \{H_i \mid 0\leq i\leq d \} \\
& & T_{H,H_{j_1}}-\sum_{i=0}^d \lambda_i\cdot T_{H_i,H_{j_1}} \mbox{ if } \ell_H(z)=\sum_{i=0}^d \lambda_i z_i,\;H\in \cH_{(j_h)_h}.
\end{eqnarray*}
The isomorphism is given by $T_{H,H'}\mapsto \pi^m \ell_H/\ell_{H'}\in \cO(U_{(j_h)_{h}}).$
 
From now on, it suffices to treat the case $\cF=\cO.$ We consider the $\check{{\rm C}}$ech  complex with respect to the covering $U_{(j_h)_h}$ for varying $(j_h)_h$. The analytic functions on the intersections of the various sets  $U_{(j_h)_h}$ are described in the same manner. It is checked that the boundary maps are closed.
In particular, the cohomology groups $H^\ast((\P^d_K)^{rig}\setminus D_m,\cO)$ are $K$-Fr\'echet spaces. Furthermore the transition maps
are dense. Thus by Proposition \ref{projlimII} we get\footnote{Note that $H^\ast((\P^d_K)^{ad}\setminus  D^{ad}_{m+1},\cF)= H^\ast((\P^d_K)^{rig}\setminus  D_{m+1},\cF)$ resp. $H^\ast((\P^d_K)^{ad}\setminus Z_I^{P^n_I}, \cF)= H^\ast((\P^d_K)^{rig}\setminus P^n_I\cdot Y_I^{rig}, \cF)$ since the corresponding topoi are equivalent cf. \cite{H} Prop. 2.1.4. }
\begin{eqnarray}\label{limproj}
H^\ast((\P^d_K)^{ad}\setminus Z_I^{P^n_I}, \cF)= \varprojlim\nolimits_m H^\ast((\P^d_K)^{ad}\setminus  D^{ad}_m,\cF).
\end{eqnarray}
Furthermore, the transition maps
$$H^\ast((\P^d_K)^{ad}\setminus  Z_I^{P_I^{n+1}} ,\cF) \rightarrow  H^\ast((\P^d_K)^{ad}\setminus   Z_I^{P_I^n},\cF)$$ are dense.
Thus, taking (\ref{Translation}) into account, our projective system  consists of $K$-Fr\'echet spaces and satisfies the topological Mittag-Leffler property.
\end{proof}

\noindent We deduce  from  \cite{EGAIII} 13.2.4 that $$\varprojlim_{n\in \N}\nolimits^{(1)}  \Big(\bigoplus_{g\in G_0/P^n_I} H^{i-1}_{Z^{gP^n_I}_I}(\P^d_K,\cF)\Big)_{n\in \N} = 0.$$
We obtain the identity
$${\rm Ext}^i(i_\ast(\sideset{}{'}\prod\limits_{g\in G/P_I} \Z_{g,I}),\cF)\, \cong \, \varprojlim_{n\in \N} \bigoplus_{g\in G_0/P^n_I} H^i_{Z^{gP^n_I}_I}(\P^d_K, \cF) . $$
On the other hand, we have  $\bigcap_{n\in \N} Z_I^{P_I^n} =  \bigcap_{n\in \N} Y_I(\epsilon_n)^{ad} = Y_I^{ad}.$ Again, by applying Proposition \ref{projlimI},  we deduce the identity $\varprojlim_{n}  H^\ast_{Z^{P^n_I}_I}(\P^d_K, \cF) = \varprojlim_{n}  H^\ast_{Y_I(\epsilon_n)^{ad}}(\P^d_K, \cF).$ We get
\begin{eqnarray*}
\varprojlim_{n} \bigoplus_{g\in G_0/P^n_I} H^\ast_{Z^{gP^n_I}_I}(\P^d_K, \cF) & = &  \varprojlim_{n}\bigoplus_{g\in G_0/P_I^n}  H^\ast_{gY_I(\epsilon_n)^{ad}}(\P^d_K,\cF) \\
& = & \varprojlim_{n}\bigoplus_{g\in G_0/P_I^n} H_{gY_I(\epsilon_n)}^\ast(\P^d_K, \cF).   
\end{eqnarray*}
\noindent Thus the statement of our proposition is proved.
\qed

Consider the spectral sequence
\begin{multline}\label{ss}
E_1^{-p,q} =  {\rm Ext}^q(\!\!\!\!\!\bigoplus_{I \subset \Delta
\atop |\Delta\setminus I|=p+1}\!\!\!\!\! i_\ast(\sideset{}{'}\prod_{g\in
G/P_I}\limits \Z_{g,I}),\cF)  \Rightarrow {\rm
Ext}^{-p+q}(i_\ast(\Z_{{\cY}^{ad}}), \cF)=
H^{-p+q}_{{\cY}^{ad}}(\P_K^d, \cF)
\end{multline}
induced by the acyclic complex (\ref{complex}), cf. Theorem \ref{Theorem1}.
By applying  the previous proposition to it we compute for the rows $E_1^{\bullet,q},\; q\in \N,$ the following complexes of $K$-Fr\'echet spaces.
Note that we have
$$H_{\P_K^{d-j}(\epsilon_n)}^{\ast}(\P^d_K,\cF)=H_{\P_K^{d-j}(\epsilon_n)}^{j}(\P^d_K,\cF) \oplus \bigoplus_{k=j+1}^d
H^{k}(\P^d_K,\cF)$$
by (\ref{lokKoho}).
\begin{eqnarray*}
E_1^{\bullet,d} &\!\!\!\!\!: & \!\!\!\!\!\varprojlim_{n}\bigoplus_{g\in G_0/P_\emptyset^n} H_{g\P_K^0(\epsilon_n)}^d(\P^d_K, \cF)\rightarrow \bigoplus_{I\subset \Delta \atop { \#I=1 }} \varprojlim_{n}\bigoplus_{g \in G_0/P_I^n} M_{g,I}^d 
\rightarrow \bigoplus_{I\subset \Delta \atop { \#I=2}} \varprojlim_{n}\bigoplus_{g \in G_0/P_I^n} M_{g,I}^d \\ & & \rightarrow  \ldots \rightarrow \varprojlim_{n}\bigoplus_{g \in G_0/P_{(1,d)}^n} M_{g,I}^d,\\
\mbox{where}& & \\ \\
& & M_{g,I}^d=\left\{\begin{array}{ccc}  H_{g\P_K^0(\epsilon_n)}^d(\P^d_K, \cF) & ; & \alpha_0 \notin I \\  \\
H^d(\P^d_K,\cF) & ; & \alpha_0 \in I \end{array} \right.,
\end{eqnarray*}
\begin{eqnarray*} \\
E_1^{\bullet,d-1} & \!\!\!\!\!: & \!\!\!\!\!\varprojlim_{n}\!\!\!\bigoplus_{g\in G_0/P_{(2,1,\ldots,1)}^n} \!\!\!\!\!\!\!\!\!\! H_{g\P_K^1(\epsilon_n)}^{d-1}(\P^d_K, \cF)\rightarrow \bigoplus_{I\subset \Delta \atop {\#I=2 \atop \alpha_0\in I }} \varprojlim_{n}\bigoplus_{g \in G_0/P_I^n} M^{d-1}_{g,I}
\rightarrow \bigoplus_{I\subset \Delta \atop { \#I=3 \atop \alpha_0\in I}} \varprojlim_{n}\bigoplus_{g \in G_0/P_I^n} M^{d-1}_{g,I} \\ & &  \rightarrow  \ldots  \rightarrow \varprojlim_{n}\bigoplus_{g \in G_0/P_{(2,d-1)}^n} M^{d-1}_{g,I},\\
\mbox{where} & & \\ \\
& & M_{g,I}^{d-1}=\left\{\begin{array}{ccc}  H_{g\P_K^1(\epsilon_n)}^{d-1}(\P^d_K, \cF) & ; & \alpha_1 \notin I \\  \\
H^{d-1}(\P^d_K,\cF) & ; & \alpha_1 \in I \end{array} \right.,
\end{eqnarray*}
\begin{eqnarray*}
& & \centerline{$\vdots$}\\ \\
E_1^{\bullet,j} & \!\!\!\!\!: &\!\!\!\!\!\varprojlim_{n}\!\!\!\bigoplus_{g\in G_0/P_{(d+1-j,1,\ldots,1)}^n} \!\!\!\!\!\!\!\!\!\!\!H_{g\P_K^{d-j}(\epsilon_n)}^{j}(\P^d_K, \cF)\rightarrow \!\!\!\!\!\!\!\!\!\!\! \bigoplus_{I\subset \Delta \atop { \#I=d-j+1 \atop \alpha_0,\ldots, \alpha_{d-j-1}\in I}} \!\!\!\!\!\!\!\!\varprojlim_{n}\bigoplus_{g \in G_0/P_I^n} M_{g,I}^{j}
\rightarrow \!\!\!\!\!\!\!\!\!\!\!\bigoplus_{I\subset \Delta \atop { \#I=d-j+2 \atop \alpha_0\ldots \alpha_{d-j-1} \in I}}\!\!\!\!\!\!\!\varprojlim_{n}\bigoplus_{g \in G_0/P_I^n} M_{g,I}^j \\ & & \rightarrow \ldots  \rightarrow \varprojlim_{n}\bigoplus_{g \in G_0/P_{(d+1-j,j)}^n}  M_{g,I}^j,\\
\mbox{where} & &\\ \\
& & M_{g,I}^{j}=\left\{\begin{array}{ccc}  H_{g\P_K^{d-j}(\epsilon_n)}^{j}(\P^d_K, \cF) & ; & \alpha_{d-j} \notin I \\  \\
H^{j}(\P^d_K,\cF) & ; & \alpha_{d-j} \in I \end{array} \right., 
\end{eqnarray*}
\begin{eqnarray*}\\ \\
& & \centerline{$\vdots$}  \\ \\
E_1^{0,1} &\!\!\!\!\!: &  \varprojlim_{n}\bigoplus_{g\in G_0/P_{(d,1)}^n} H_{g\P_K^{d-1}(\epsilon_n)}^1(\P^d_K, \cF).
\end{eqnarray*}

\noindent Here, the very left term in each row $E^{\bullet,j}_1$ sits in degree $-j+1.$
We can rewrite these complexes in terms of induced representations. Here we abbreviate 
$$(d+1-j,1^j):=(d+1-j,1,\ldots,1)$$
for any decomposition $(d+1-j,1,\ldots,1)$ of $d+1.$

\begin{eqnarray*}
E_1^{\bullet,d} & \!\!\!\!\!  : & \!\!\!\!\!\varprojlim_{n} {\rm Ind}^{G_0}_{P_\emptyset^n} H_{\P_K^0(\epsilon_n)}^d(\P^d_K, \cF)\rightarrow  \varprojlim_{n} \bigoplus_{I\subset \Delta \atop { \#I=1}} {\rm Ind}^{G_0}_{P_I^n} M_{g,I}^d
\rightarrow \varprojlim_{n}\bigoplus_{I\subset \Delta \atop { \#I=2 }} {\rm Ind}^{G_0}_{P_I^n} M_{g,I}^d \\ & &  \rightarrow \ldots \rightarrow \varprojlim_{n} {\rm Ind}^{G_0}_{P_{(1,d)}^n} M_{g,I}^d,  \\ \\
E_1^{\bullet,d-1} &\!\!\!\!\! : &\!\!\!\!\! \varprojlim_{n} {\rm Ind}^{G_0}_{P_{(2,1^{d-1})}^n} H_{\P_K^1(\epsilon_n)}^{d-1}(\P^d_K, \cF)\rightarrow  \varprojlim_{n} \bigoplus_{I\subset \Delta \atop { \#I=2 \atop \alpha_0\in I}} {\rm Ind}^{G_0}_{P_I^n} M_{g,I}^{d-1}
\rightarrow \varprojlim_{n} \bigoplus_{I\subset \Delta \atop { \#I=3 \atop \alpha_0\in I}} {\rm Ind}^{{\bf} G_0}_{P_I^n} M_{g,I}^{d-1} \\ & & \rightarrow \ldots  \rightarrow \varprojlim_{n} {\rm Ind}^{G_0}_{P_{(2,d-1)}^n} M_{g,I}^{d-1} \\
& & \centerline{$\vdots$}\\ \\ 
E_1^{\bullet,j} &\!\!\!\!\!: &\!\!\!\!\!\varprojlim_{n} {\rm Ind}^{G_0}_{P_{(d+1-j,1^j)}^n} \!\!\!\! H_{\P_K^{d-j}(\epsilon_n)}^{j}(\P^d_K, \cF)\rightarrow  \varprojlim_{n}\!\!\!\!\!\!\!\!\!\!\!\bigoplus_{I\subset \Delta \atop { \#I=d-j+1 \atop \alpha_0,\ldots, \alpha_{d-j-1}\in I }}\!\!\!\!\!\!\!\!\!\!\! {\rm Ind}^{G_0}_{P_I^n} M_{g,I}^j
\rightarrow \varprojlim_{n}\!\!\!\!\!\! \bigoplus_{I\subset \Delta \atop { \#I=d-j+2 \atop \alpha_0\ldots \alpha_{d-j-1} \in I}}   \!\!\!\!\!\!{\rm Ind}^{G_0}_{P_I^n} M_{g,I}^j \\ & &  \rightarrow \ldots  \rightarrow \varprojlim_{n}{\rm Ind}^{G_0}_{P_{(d+1-j,j)}^n}  M_{g,I}^j\\ \\
& & \centerline{$\vdots$} \\ \\ 
E_1^{0,1} & \!\!\!\!\!: & \!\!\!\!\!\varprojlim_{n} {\rm Ind}^{G_0}_{P_{(d,1)}^n} H_{\P_K^{d-1}(\epsilon_n)}^1(\P^d_K, \cF).
\end{eqnarray*}

\medskip
\begin{Proposition}\label{E_1 acyclic}
Each of the complexes $E^{\bullet,j}_1, \; j=1,\ldots,d$,  is acyclic apart from the very left and right position.
\end{Proposition}

\proof  We can write each of the complexes in the shape $E^{\bullet,j}_1=\varprojlim_n K_{j,n}^\bullet,$
where $K_{j,n}^\bullet$ is a complex of $K$-Fr\'echet spaces which appears in a short exact sequence of complexes
of $K$-Fr\'echet spaces
\begin{eqnarray}\label{kurzeexakteSequenz}
0\rightarrow K_{j,n}^{\bullet,'} \rightarrow K_{j,n}^\bullet \rightarrow K_{j,n}^{\bullet,''} \rightarrow 0 .
\end{eqnarray}
Here, $K_{j,n}^{\bullet,'}$ is the complex
$${\rm Ind}^{G_0}_{P_{(d+1-j,1^j)}^n} \tilde{H}^j_{\P^{d-j}_K(\epsilon_n)}(\P^d_K,\cF) \rightarrow \!\!\!\!\!\!\!\!\!\!\!\!\!\bigoplus_{I\subset \Delta \atop {\#I=d-j+1 \atop \alpha_0,\ldots, ,\alpha_{d-j-1}\in I, \alpha_{d-j} \notin I }}\!\!\!\!\!\!\!\!\!\!\!\!\! {\rm Ind}^{G_0}_{P_I^n} \tilde{H}^j_{\P^{d-j}_K(\epsilon_n)}(\P^d_K,\cF) $$
$$\rightarrow\!\!\!\!\!\!\!\!\!\!\!\!\!\bigoplus_{I\subset \Delta \atop { \#I=d-j+2 \atop \alpha_0,\ldots ,\alpha_{d-j-1} \in I, \alpha_{d-j} \notin I}} \!\!\!\!\!\!\!\!\!\!\!\!\! {\rm Ind}^{G_0}_{P_I^n} \tilde{H}^j_{\P^{d-j}_K(\epsilon_n)}(\P^d_K,\cF) \rightarrow \ldots  \rightarrow {\rm Ind}^{G_0}_{P_{(d+1-j,j)}^n}   \tilde{H}^j_{\P^{d-j}_K(\epsilon_n)}(\P^d_K,\cF).$$

\noindent Furthermore, the complex $K_{j,n}^{\bullet,''}$ is given by
\begin{eqnarray*}
&{\rm Ind}^{G_0}_{P_{(d+1-j,1^j)}^n} \!\!\!\! H^{j}(\P^d_K, \cF) & \rightarrow\!\!\!\! \bigoplus_{I\subset \Delta \atop { \#I=d-j+1 \atop \alpha_0,\ldots, \alpha_{d-j-1}\in I }} \!\!\!\!  {\rm Ind}^{G_0}_{P_I^n}  H^{j}(\P^d_K, \cF) 
\rightarrow \!\!\!\!  \bigoplus_{I\subset \Delta \atop { \#I=d-j+2 \atop \alpha_0,\ldots ,\alpha_{d-j-1} \in I}}\!\!\!\!  {\rm Ind}^{G_0}_{P_I^n}  H^{j}(\P^d_K, \cF) \rightarrow \\ & \ldots & \rightarrow  \bigoplus_{I\subset \Delta \atop {\#(\Delta\setminus I)=1 \atop \alpha_0,\ldots ,\alpha_{d-j-1} \in I}}  {\rm Ind}^{G_0}_{P_{I}^n}   H^{j}(\P^d_K, \cF) .
\end{eqnarray*}
Since $H^{j}(\P^d_K, \cF)$ is a $G$-module, this complex is isomorphic to
$$\Big({\rm Ind}^{G_0}_{P_{(d+1-j,1^j)}^n} \!\!\!\! K \rightarrow \bigoplus_{I\subset \Delta \atop { \#I=d-j+1 \atop \alpha_0,\ldots, \alpha_{d-j-1}\in I }} {\rm Ind}^{G_0}_{P_I^n} K
 \rightarrow \ldots  \rightarrow  \bigoplus_{I\subset \Delta \atop {\#(\Delta\setminus I)=1 \atop \alpha_0,\ldots ,\alpha_{d-j-1} \in I}} {\rm Ind}^{G_0}_{P_{I}^n} K\Big) \otimes  H^{j}(\P^d_K, \cF).$$
It suffices to prove that the complexes  $\varprojlim_n K_{j,n}^{\bullet,'}$ and  $\varprojlim_n K_{j,n}^{\bullet,''}$ are acyclic apart from the very left and right position. We deduce this property for $K_{j,n}^{\bullet,''}$ by applying the first part of the following lemma to its dual.

\begin{Lemma}\label{Steinberg} (i) For each integer $0 \leq j \leq d-1,$ the following complex is acyclic apart from the very left and right position:\medskip
$$\Big({\rm Ind}^{G_0}_{P_{(d+1-j,1^j)}^n} \!\!\!\! K \rightarrow \bigoplus_{I\subset \Delta \atop { \#I=d-j+1 \atop \alpha_0,\ldots, \alpha_{d-j-1}\in I }} {\rm Ind}^{G_0}_{P_I^n} K
 \rightarrow \ldots  \rightarrow  \bigoplus_{I\subset \Delta \atop {\#(\Delta\setminus I)=1 \atop \alpha_0,\ldots ,\alpha_{d-j-1} \in I}} {\rm Ind}^{G_0}_{P_{I}^n} K\Big) \otimes  H^{j}(\P^d_K, \cF).$$
(ii) Let $W$ be any ${\bf P_{(d+1-j,j)}}(O_K/(\pi^n))$-module. Consider $W$ as a $P_{(d+1-j,j)}^n$-module via the inflation map (\ref{reduction_map}). Then for each integer $0 \leq j \leq d-1,$ the following complex is acyclic apart from the very left and right position:
$${\rm Ind}^{G_0}_{P_{(d+1-j,j)}^n}  W \rightarrow  \ldots  \rightarrow \!\!\!\!\!\!\!\!\!\!\!\!\!\!\bigoplus_{I\subset \Delta \atop { \#I=d-j+2 \atop \alpha_0,\ldots,\alpha_{d-j-1} \in I,\alpha_{d-j}\not\in I}}\!\!\!\!\!\!\!\!\!\!\!\!\!\! {\rm Ind}^{G_0}_{P_I^n} W \rightarrow \!\!\!\!\!\!\!\!\!\!\!\!\!\! \bigoplus_{I\subset \Delta \atop { \#I=d-j+1 \atop \alpha_0,\ldots,\alpha_{d-j-1} \in I, \alpha_{d-j}\not\in I}}\!\!\!\!\!\!\!\!\!\!\!\!\!\! {\rm Ind}^{G_0}_{P_I^n} W  \rightarrow  {\rm Ind}^{G_0}_{P_{(d+1-j,1^j)}^n} W $$
\end{Lemma}

\proof For a subset $I\subset \Delta,$ there is a natural isomorphism  
$$G_0/P_I^n \stackrel{\sim}{\longrightarrow} {\bf G}(O_K/(\pi^n))/{\bf P_I}(O_K/(\pi^n)).$$ 
Via this identification the representation ${\rm Ind}^{G_0}_{P_I^n} W$ coincides with ${\rm Ind}^{{\bf G}(O_K/(\pi^n))}_{{\bf P_I}(O_K/(\pi^n))} W.$  Then statement (i) follows from Theorem 2.5, ch. III in \cite{OR}.
In fact, loc.cit. treats the dual complex in the case $n=1$, but the proof for $n>1$ is the same. 
The proof of part (ii) works by the same reasoning as in loc.cit. In particular, it does not depend on the coefficient system.  \qed

\noindent Since the complex $K_{j,n}^{\bullet,''}$ consists of finite-dimensional $K$-vector spaces, we conclude by the Mittag-Leffler condition that
$\varprojlim_n K_{j,n}^{\bullet,''}$ is acyclic apart from the very left and right position.
The dual complex $(\varprojlim_n K_{j,n}^{\bullet,''})'$ of  $\varprojlim_n K_{j,n}^{\bullet,''}$ is given by
\begin{eqnarray*}
 \varinjlim_n {\rm Ind}^{G_0}_{P_{(d+1-j,1^j)}^n}  H^{j}(\P^d_K, \cF)' & \leftarrow & \!\!\!\!\!\!\!\!\!\!\!\!\!  \bigoplus_{I\subset \Delta \atop { \#I=d-j+1 \atop \alpha_0,\ldots, \alpha_{d-j-1}\in I }}\!\!\!\!\!\!\!\!\!\! \varinjlim_n  {\rm Ind}^{G_0}_{P_I^n}  H^{j}(\P^d_K, \cF)' 
\leftarrow \!\!\!\!\!\!\!\!\!\!   \bigoplus\limits_{I\subset \Delta \atop { \#I=d-j+2 \atop \alpha_0,\ldots ,\alpha_{d-j-1} \in I}}\!\!\!\!\!\!\!\!\!\! \varinjlim_n  {\rm Ind}^{G_0}_{P_I^n}  H^{j}(\P^d_K, \cF)' \\ \!\!\!\!\!\!\!\!\!\!\!\!\!   \leftarrow  \ldots & \leftarrow & \!\!\!\!\!\!\!\!\!\!\!\! \bigoplus_{I\subset \Delta \atop {\#(\Delta\setminus I)=1 \atop \alpha_0,\ldots ,\alpha_{d-j-1} \in I}} \!\!\!\!\!\!\!\!\!\! \varinjlim_n  {\rm Ind}^{G_0}_{P_{I}^n}   H^{j}(\P^d_K, \cF)' \\ \\ 
= {\rm Ind}^{\infty,G}_{P_{(d+1-j,1^j)}}H^{j}(\P^d_K, \cF)' & \leftarrow & \!\!\!\!\!\!\!\!\!\!\!\! \bigoplus_{I\subset \Delta \atop { \#I=d-j+1 \atop \alpha_0,\ldots, \alpha_{d-j-1}\in I }} \!\!\!\!\!\!\!\!\!\! {\rm Ind}^{\infty, G}_{P_I}  H^{j}(\P^d_K, \cF)' 
\leftarrow  \!\!\!\!\!\!\!\!\!\!   \bigoplus\limits_{I\subset \Delta \atop { \#I=d-j+2
\atop \alpha_0,\ldots ,\alpha_{d-j-1} \in I}}  \!\!\!\!\!\!\!\!\!\!  {\rm Ind}^{\infty,
G}_{P_I}  H^{j}(\P^d_K, \cF)' \\   \leftarrow \ldots &  \leftarrow & \!\!\!\!\!\!\!\!\!\!\!\! \bigoplus_{I\subset \Delta \atop {\#(\Delta\setminus I)=1 \atop \alpha_0,\ldots ,\alpha_{d-j-1} \in I}}\!\!\!\!\!\!\!\!\!\!  {\rm
Ind}^{\infty,G}_{P_{I}}   H^{j}(\P^d_K, \cF)'.
\end{eqnarray*}
Here ${\rm Ind}^{\infty,G}_P$ denotes the (unnormalized) smooth
induction functor for a parabolic subgroup $P \subset G,$ cf. \cite{Ca}. Again, since
$H^{j}(\P^d_K, \cF)$ is even a $G$-module, this complex coincides 
with
$$\Big( {\rm Ind}^{\infty,G}_{P_{(d+1-j,1^j)}} K \leftarrow \!\!\!\!\!\!\!\!\!\! \bigoplus_{I\subset \Delta \atop { \#I=d-j+1 \atop \alpha_0,\ldots, \alpha_{d-j-1}\in I }}\!\!\!\!\!\!\!\!\!\! {\rm Ind}^{\infty, G}_{P_I} K
\leftarrow
\ldots  \leftarrow \bigoplus_{I\subset \Delta \atop {\#(\Delta\setminus I)=1 \atop \alpha_0,\ldots ,\alpha_{d-j-1} \in I}} {\rm Ind}^{\infty,G}_{P_I} K \Big)
\otimes H^{j}(\P^d_K, \cF)' .$$
Let
$$v^G_{P_{(d+1-j,1^j)}}(K):= {\rm Ind}^{\infty,G}_{P_{(d+1-j,1^j)}}K/\sum_{Q\supsetneq P_{(d+1-j,1^j)} } {\rm Ind}^{\infty,G}_{Q} K$$ be the smooth generalized Steinberg representation  with respect to the parabolic subgroup  $P_{(d+1-j,1^j)}.$
It is known that this is an irreducible smooth $G$-representation, cf.\, \cite{BW} ch. X.
Put
$$v^G_{P_{(d+1-j,1^j)}}( H^{j}(\P^d_K, \cF)'):= v^G_{P_{(d+1-j,1^j)}}(K) \otimes  H^{j}(\P^d_K, \cF)'$$
The only non-vanishing cohomology groups of the dual complex $(\varprojlim_n K_{j,n}^{\bullet,''})'$  are therefore given by
$$H^\ast((\varprojlim_n K_{j,n}^{\bullet,''})')= v^G_{P_{(d+1-j,1^j)}}(H^{j}(\P^d_K, \cF)')\oplus H^{j}(\P^d_K, \cF)'  \mbox{ for }\; j\geq 2$$
resp.
$$H^\ast((\varprojlim_n K_{j,n}^{\bullet,''})')={\rm Ind}^{\infty,G}_{P_{(d,1)}}  H^{1}(\P^d_K, \cF)'  \;\mbox { for } \; j=1.$$

 Now, we turn to the complexes $\varprojlim_n K_{j,n}^{\bullet,'}.$ Each entry in $\varprojlim_n K_{j,n}^{\bullet,'}$ is a compact projective
 limit of $K$-Fr\'echet spaces, hence nuclear, cf. \cite{S2} Proposition 19.9 (compare also the example at the end of chapter 16).
 By loc.cit. Corollary 19.3 these objects are reflexive. The duality functor is exact on the category of $K$-Fr\'echet spaces, cf. \cite{Ba} ch I, Cor. 1.4. So, it suffices to show
that the dual complexes (see \cite{S2} Prop. 16.10)
$$\varinjlim_{n} {\rm Ind}^{G_0}_{P_{(d+1-j,1^j)}^n} \tilde{H}^j_{\P^{d-j}_K(\epsilon_n)}(\P^d_K,\cF)' \leftarrow  \!\!\!\!\!\!\!\!\!\!\bigoplus_{I\subset \Delta \atop {\#I=d-j+1 \atop \alpha_0,\ldots, \alpha_{d-j-1}\in I, \alpha_{d-j} \notin I }} \!\!\!\!\!\!\!\!\!\! \varinjlim_{n} {\rm Ind}^{G_0}_{P_I^n} \tilde{H}^j_{\P^{d-j}_K(\epsilon_n)}(\P^d_K,\cF)' $$
$$\leftarrow  \!\!\!\!\!\!\!\!\!\!\bigoplus_{I\subset \Delta \atop { \#I=d-j+2 \atop \alpha_0\ldots \alpha_{j-1} \in I, \alpha_{d-j} \notin I}}  \!\!\!\!\!\!\!\!\!\! \varinjlim_{n} {\rm Ind}^{G_0}_{P_I^n} \tilde{H}^j_{\P^{d-j}_K(\epsilon_n)}(\P^d_K,\cF)' \leftarrow \ldots  \leftarrow \varinjlim_{n} {\rm Ind}^{G_0}_{P_{(d+1-j,j)}^n}  \tilde{H}^{j}_{\P^{d-j}_K(\epsilon_n)}(\P^d_K,\cF)'$$
consisting of compact inductive limits of locally convex $K$-vector spaces are exact apart from the very left and right position.
But this follows by the exactness of $\varinjlim$ from Lemma \ref{Steinberg}, as well.  Thus Proposition \ref{E_1 acyclic} is proved. \qed

\begin{Remark}
{\rm Alternatively, on could prove the acyclicity of $\varprojlim_n K_{j,n}^{\bullet,'}$ as follows. By the same reasoning as in the proof of Lemma \ref{Steinberg} one shows that even the complex $K_{j,n}^{\bullet,'}$ is exact apart from the very left and right position. Then we apply the topological Mittag-Leffler condition  to the projective limit to deduce the claim, cf. Lemma \ref{TopML}.} \qed
\end{Remark}

\medskip
By Proposition  \ref{E_1 acyclic} we deduce that 
the only non-vanishing entries in $E_2^{p,q}$ are given by the indices $(p,q)=(-j+1,j),\; j=1,\ldots,d,$ and $(p,q)=(0,j),\; j\geq 2.$
For the latter indices, we get $$E_2^{0,q}=H^q(\P^d,\cF).$$
For the other indices, we obtain
$$E_2^{-j+1,j}={\rm ker}(E_1^{-j+1,j} \rightarrow E_1^{-j+2,j})=$$
$$ {\rm ker}\,\Big(\varprojlim_{n} {\rm Ind}^{G_0}_{P_{(d+1-j,1^j)}^n} H_{\P_K^{d-j}(\epsilon_n)}^{j}(\P^d_K, \cF) \rightarrow
\varprojlim_{n} \!\!\!\bigoplus_{I\subset \Delta \atop { \# I=
d-j+1  \atop \alpha_0,\ldots, \alpha_{d-j-1} \in I}}
\!\!\!\!\!\!\!\!\!\!{\rm Ind}^{G_0}_{P_I^n}
M_{g,I}^j\Big).$$ 
Thus our spectral sequence has apart from stretching the	y-axis the same structure as in the case of constant coefficients, cf. p.70 \cite{SS}. Further, the composed maps $E_2^{0,s} \rightarrow H_{\cY}^s(\P^d_K,\cF)\rightarrow H^s(\P^d_k,\cF)$ where the first map
is the edge homomorphism are isomorphisms for $s>1$ and surjective for $s=1.$
By the same reasoning as in loc.cit. we conclude that our  spectral sequence degenerates at $E_2.$
By duality, i.e., by taking the strong dual of
these $K$-Fr\'echet spaces, we get locally analytic (cf. \cite{ST3} Cor. 3.3) $G_0$-representations
$$(E_2^{-j+1,j})' = {\rm coker}\,\Big(\varinjlim_{n} \!\!\!\!\! \bigoplus_{I\subset \Delta \atop { \# I= d-j+1  \atop \alpha_0,\ldots, \alpha_{d-j-1} \in I}}\!\!\!\!\! {\rm Ind}^{G_0}_{P_I^n} (M_{g,I}^j)' \rightarrow \varinjlim_{n} {\rm Ind}^{G_0}_{P_{(d+1-j,1^j)}^n} H_{\P_K^{d-j}(\epsilon_n)}^{j}(\P^d_K, \cF)'\Big) ,$$
\noindent which are by (\ref{kurzeexakteSequenz})  extensions of locally analytic $G_0$-representations
\bigskip
\begin{multline}\label{extensions}
\!\!\!\!\!0\rightarrow v^{G}_{{P_{(d+1-j,1^j)}}}(H^j(\P^d_K,\cF)')
\rightarrow (E_2^{-j+1,j})' \rightarrow \varinjlim_{n} {\rm
Ind}^{G_0}_{P_{(d+1-j,j)}^n}( \tilde{H}^j_{
\P^{d-j}_K(\epsilon_n)}(\P^d_K, \cF)' \otimes {\rm St}_j) \rightarrow 0.
\end{multline} 

\noindent Here, ${\rm St}_j=v^{GL_j}_{B \cap GL_j}(K)$ denotes the smooth Steinberg representation of the $p$-adic Lie group ${\rm GL}_j$ viewed as one of the factors of the Levi subgroup $L_{(d-j+1,j)}.$ Now, we plug in the result of Proposition \ref{Duality}.  We obtain isomorphisms of locally analytic $G_0$-representations:\bigskip
\begin{eqnarray*}
\varinjlim_{n} {\rm Ind}^{G_0}_{P_{(d+1-j,j)}^n}(\tilde{H}^j_{ \P^{d-j}_K(\epsilon_n)}(\P^d_K,\cF)' \otimes {\rm St}_j)\!\!\!  &\cong & \!\!\! \varinjlim_{n}  {\rm Ind}^{G_0}_{P^n_{(d-j+1,j)}}(\cO(U^{+,n}_{(d-j+1,j)},N_{d-j}')^{{\frd_{d-j}}}\otimes {\rm St}_{j}) \\
\!\!\! &  =  & \!\!\! \varinjlim_{n}  {\rm Ind}^{G_0}_{P^n_{(d-j+1,j)}}(\cO(U^{+,n}_{(d-j+1,j)},N_{d-j}') \otimes {\rm St}_{j})^{{\frd_{d-j}}}. 
\end{eqnarray*}
The latter equality holds since ${\rm St}_j$ is smooth. On the other hand,
\begin{eqnarray*}
\varinjlim_{n}  {\rm Ind}^{G_0}_{P^n_{(d-j+1,j)}}(\cO(U^{+,n}_{(d-j+1,j)},N_{d-j}') \otimes {\rm St}_{j})\!\!\! & \cong & \!\!\!C^{an}(G,P_{(d-j+1,j)};N_{d-j}' \otimes {\rm St}_{j}).
\end{eqnarray*}

\noindent Here $C^{an}(G,P_{(d-j+1,j)};N_{d-j}' \otimes {\rm St}_{j})$ denotes the locally analytic induced $G$-representation \cite{F} with values in $N_{d-j}' \otimes {\rm St}_j$:
\begin{eqnarray*}
C^{an}(G,P_{(d-j+1,j)};N_{d-j}' \otimes {\rm St}_{j})  = & \Big\{ & \mbox{locally analytic maps } f:G \rightarrow N_{d-j}' \otimes {\rm St}_{j} \mid  \\
 & &f(g\cdot p) =p^{-1}f(g) \; \forall g\in G, p\in P_{(d-j+1,j)} \Big\}.
\end{eqnarray*}
In the above formula, we have made use of the canonical identity
$$C^{an}(G_0,{\bf P_{(d-j+1,j)}}(O_K);N_{d-j}' \otimes {\rm St}_{j}) =  C^{an}(G,P_{(d-j+1,j)};N_{d-j}' \otimes {\rm St}_{j})$$
as locally analytic $G_0$-representations. This identity follows from the Iwasawa decomposition $G=G_0\cdot  P_{(d-j+1,j)}$ and by \cite{F} 4.1.4. Hence, the above space possesses even the structure of a locally analytic $G$-representation. We claim that  the above isomorphisms and the resulting extensions (\ref{extensions}) are even $G$-equivariant.
In fact, this follows essentially from the Iwasawa decomposition and Remark \ref{Bemerkung}.
Altogether, we get an isomorphism 
$$\varinjlim_{n} {\rm Ind}^{G_0}_{P_{(d+1-j,j)}^n}(\tilde{H}^j_{ \P^{d-j}_K(\epsilon_n)}(\P^d_K,\cF)' \otimes {\rm St}_j) \cong
C^{an}(G,P_{(d-j+1,j)};N_{d-j}' \otimes {\rm St}_{j})^{{\frd_{d-j}}}.$$

In the following, we also use for a parabolic
subgroup $P\subset G$ and a locally analytic $P$-representation $U$, the symbol ${\rm
Ind}^{an,G}_{P}(U)$ instead of $C^{an}(G,P; U).$

\begin{Lemma}\label{Filtration}
Let $V^\bullet= V^{-d+1} \supset V^{-d+2} \supset \cdots\supset  V^{-1} \supset V^0 \supset V^1= (0)$ be the canonical filtration on  $V^{-d+1}=H^1_{\cY}(\P^d_K,\cF)$  defined by our spectral sequence
(\ref{ss}). Then the subspaces $V^j$ are closed in the $K$-Fr\'echet space $H^1_{\cY}(\P^d_K,\cF)$. 
\end{Lemma}

\proof Recall the definition of $V^\bullet$. Let $K^\bullet$ be the simple complex attached to the double complex $E_0^{\bullet,\bullet}-$  the $0$-th term of our spectral sequence. Let $F^p(K^\bullet)=\sum_{j\geq p} E_0^{j,\bullet}$ be the  sum of the columns $E_0^{\bullet,j}$ with $j\geq p.$ Then $V^{-d+i}$ is by definition the kernel of the natural homomorphisms $H^1(K^\bullet) \rightarrow H^1(K^\bullet/F^{-d+i}(K^\bullet))$, 
or equivalently, it is the image of the natural homomorphism  $H^1(F^{-d+i}(K^\bullet)) \rightarrow H^1(K^\bullet)$, cf.  11.2.2 \cite{EGAIII}.
We consider the pull-back of $V^\bullet$ to $H^0(\cX,\cF)$ under the  boundary map $H^0(\cX,\cF)\rightarrow H^1_{\cY}(\P^d_K,\cF)$.
On $\cX$ the vector bundle $\cF$ splits, cf. (\ref{split}). Since the canonical  filtrations with respect to $\cO^{{\rm rk} \cF}$ and $\cF$
differ topologically at most by finite-dimensional Hausdorff $K$-vector spaces,
it is enough to treat the case of $\cF=\cO^{{\rm rk} \cF}$ and hence of $\cF=\cO.$  In this case one sees by the definition of
our spectral sequence that the algebraic part of $V^{-d+i}$ consists of functions having at most $d-i+1$ poles, cf. also Proposition \ref{cohStrukturgarbe}. Thus, our filtration coincides up to the numbering with Pohlkamps (\ref{Pohlkampfiltration}) which consists of closed subspaces.  \qed

Summarising the computation of this chapter we obtain the following theorem.

\begin{Theorem} 
Let $$V^{-d+1} \supset V^{-d+2} \supset \cdots\supset  V^{-1} \supset V^0 \supset V^1= (0)$$ be the canonical $G$-equivariant  filtration by closed $K$-Fr\'echet spaces  on $H^1_{\cY}(\P^d_K,\cF)$ defined by the spectral sequence (\ref{ss}).
For $j=1,\ldots,d-1$, there are extensions of locally analytic $G$-representations
$$0 \rightarrow v^{G}_{{P_{(d-j,1^{j+1)}}}}(H^{j+1}(\P^d_K,\cF)') \rightarrow (V^{-j}/V^{-j+1})' \rightarrow  {\rm Ind}^{an,G}_{P_{(d-j,j+1)}}(N_{d-j-1}' \otimes {\rm St}_{j+1})^{{\frd_{d-j-1}}} \rightarrow 0. $$
In the case $j=0$, there is an extension
$$0 \rightarrow {\rm Ind}^{\infty, G}_{{P_{(d,1)}}}(H^{1}(\P^d_K,\cF)') \rightarrow (V^0)' \rightarrow  {\rm Ind}^{an,G}_{P_{(d,1)}}(N_{d-1}' \otimes {\rm St}_{1})^{\frd_{d-1}} \rightarrow 0. $$
\end{Theorem}

\proof This follows from the above computation \qed
\vspace{0.5cm}

Consider the topological exact $G$-equivariant sequence of $K$-Fr\'echet spaces
$$0\rightarrow H^0(\P^d_K,\cF) \rightarrow H^0(\cX,\cF)\stackrel{p}{\rightarrow} H^1_{\cY}(\P^d_K,\cF) \rightarrow H^1(\P^d_K,\cF)\rightarrow 0.$$
For $i=0,\ldots,-d,$ we set
$$W^i:=p^{-1}(V^{i+1}).$$
Thus we get a $G$-equivariant filtration by closed $K$-Fr\'echet spaces 
$$W^{-d} \supset W^{-d+1} \supset \cdots \supset W^{-1} \supset W^0$$ on $W^{-d}=H^0(\cX,\cF).$

\begin{Corollary}
For $j=1,\ldots,d,$ there are extensions of locally analytic $G$-represen\-tations
$$0 \rightarrow v^{G}_{{P_{(d-j+1,1^{j})}}}(H^{j}(\P^d_K,\cF)') \rightarrow (W^{-j}/W^{-j+1})' \rightarrow  {\rm Ind}^{an,G}_{P_{(d+1-j,j)}}(N_{d-j}' \otimes {\rm St}_{j})^{{\frd_{d-j}}} \rightarrow 0. $$
In the case $j=0$, we get
$$W^0 = H^0(\P^d,\cF).$$ \qed
\end{Corollary}

For vector bundles, where the unipotent radical ${\bf U_{(1,d)}}$ acts trivially on the fibre, we can make  our result more precise using the main result Corollary \ref{irreduciblelagbraicrepr}.
\begin{Theorem}\label{Theorem2} 
Let $\cF=\cF_\lambda$ be a homogeneous vector bundle on $\P^d_K$  corresponding to the  $L_{(1,d)}$-dominant weight $\lambda \in \Z^{d+1}$.
Let $i_0\in \N$ be the unique integer with $w_{i}\ast \lambda \succeq  w_{i+1}\ast \lambda$
for all $i\geq i_0$, and
$w_{i}\ast \lambda \prec w_{i+1}\ast \lambda$
for all $i<i_0$. For $1\leq j\leq d,$ let
$$\mu_{j,\lambda}:=\left\{ \begin{array}{ccc} w_{j-1}\ast \lambda  &: j \leq i_0 \\  w_{j} \ast \lambda & : j > i_0 \end{array} \right. .$$
Write
$\mu_{j,\lambda}=(\mu',\mu'')$  with $\mu'\in \bbZ^{j}$ and $\mu''\in \bbZ^{d-j+1}.$
Further, set 
\begin{eqnarray*}
\Psi_{j,\lambda} =  \bigcup\limits_{k=0}^{|\mu''|} &\Big\{&  \!\!\!\! (\mu''+(c_1,\ldots, c_{d-j+1}),\mu'-(d_j,\ldots,d_1)) \mid \textstyle{\sum}_i c_i=\textstyle{\sum}_i d_{i}=k, c_1=0 \\ & & \!\!\!\mbox{ or }  d_1=0,\;  c_{i+1}\leq \mu_i'' -\mu_{i+1}'',  \;
i=1,\ldots,d-j,\; d_{i+1} \leq \mu'_{j-i} - \mu'_{j-i+1}, \,\\ & &  i=1,\ldots,j-1 \, \Big\}
\end{eqnarray*} and let 
$$N_{j,\lambda} = \bigoplus_{\mu \in \Psi_{j,\lambda}} V_{\mu} \;\;\;\; \subset \;\;\;\; K[X_{(m,n)}\mid m\leq d-j,\; n\geq d-j+1] \otimes  V_{(\mu'',\mu')} $$ 
be the sum of the irreducible
algebraic $P_{(d+1-j,j)}$-representations $V_{\mu}$ attached to $\mu.$
Let 
$$q_j: K[X_{(m,n)}\mid m\leq d-j,\; n\geq d-j+1] \otimes  V_{(\mu'',\mu')}  \rightarrow \tilde{H}^{j}_{\P^{d-j}_K}(\P^d_K,\cF_\lambda)$$ be the quotient map of  Corollary \ref{irreduciblelagbraicrepr}.
Then we can choose $N_{d-j}$ to be $q_j(N_{j,\lambda}).$
\end{Theorem}

\proof This follows from Corollary \ref{irreduciblelagbraicrepr}  .\qed

\begin{Remark} {\rm
By renumbering the filtration $W^\bullet$ on $H^0(\cX,\cF),$ i.e., if we set
$$\cF(\cX)^i:=W^{-d+i}$$ for $i=0,\ldots d,$ we get the filtration on $H^0(\cX,\cF)$ mentioned in the introduction.
\qed}
\end{Remark}

\noindent {\it Conclusion:} Although we have generalized the cases of $\cF=\Omega^d$ \cite{ST1} respectively $\cF=\cO$ \cite{P} to arbitrary homogeneous vector bundles on $\P^d_K$, our result has the lack that it does not yield explicit isomorphisms (partial boundary value maps)
as in loc. cit. We hope to determine these isomorphisms in a future paper.

\section{Examples}

\subsection{$\cF = \cO_{\P^d_K}$}

In this chapter we compare our result in the case $\cF =\cO= \cO_{\P^d_K}$ to that of \cite{P}.
For this purpose, we have to recall some more notation used there.

In loc.cit. $\!$there is defined a filtration by closed $K$-Fr\'echet spaces on $\cO(\cX)$ as follows.
For a character $\mu=\sum_{l=0}^d m_l\epsilon_l \in X_\ast(\overline{T}),$ we set
$$J_-(\mu)=\{l\,\mid\, 0 \leq l \leq d \mbox{ and } m_l <0\}.$$
For a subset $J\subset \{0,\ldots,d\},$ we put
$$\fra_J:= \sum_{\{\mu \,\mid\, J_-(\mu)\subset J \}} K\cdot \mu.$$
This is a $U(\frg)$-submodule of
 $$\cO_{\rm inf}(\cX):=\sum_{\mu \in X^\ast(\overline{T})} K\cdot \Xi_\mu \subset \cO(\cX)$$
and we have
$$J' \subset J \iff \fra_{J'} \subset \fra_J.$$
The extreme cases are
$$\fra_\emptyset=K\cdot \Xi_0 \,\mbox{ resp. }\, \fra_{\{0,\ldots,d\}}=\cO_{{\rm inf}}(\cX).$$
We set $\fra_J^{<}:=\sum_{J'\subsetneq J} \fra_{J'}$ for $J\neq \emptyset$
and $\fra_\emptyset^<=0.$  We obtain  $K$-vector space isomorphisms
$$\fra_J/\fra_J^>\; \stackrel{\sim}{\longrightarrow} \sum_{\{\mu\,\mid\, J=J_-(\mu)\}}K\cdot \Xi_\mu .$$
Put
\begin{eqnarray*}
A(J) & := & \{\mu \in X^\ast(\overline{T})\mid m_j=-1 \mbox{ for } j\in J_-(\mu) \mbox{ and } J_-(\mu)=J\} \\ \\
M_J & := & \sum _{\mu \in A(J)} K\cdot \Xi_\mu + \fra_J^< / \fra_J^< \\ \\
\frp_J & = & \{ (g_{i,j})\mid g_{i,j}=0 , \mbox{ if } i \not\in J
\mbox{ and } j \in J \}.
\end{eqnarray*}
Then $M_J \subset  \fra_J/\fra_J^<$ is a $\frp_J$-submodule of $\cO_{{\rm inf}}(\cX).$
For $0 \leq j \leq d,$ we set
$$\fra_j=\sum_{\{\mu\, \mid\, \# J_-(\mu) \,\leq \,j\}} K\cdot \Xi_\mu.$$
We get a filtration  by $\frg$-submodules
\begin{eqnarray}\label{LieFiltration}
K=\fra_{0} \subset \fra_{1} \subset \cdots \subset \fra_{d-1} \subset \fra_d=\cO_{{\rm inf}}(\cX)
\end{eqnarray}
 on $\cO_{{\rm inf}}(\cX).$
For $0\leq j \leq d,$ we put
$$\overline{j}:=\{d-j+1,\ldots,d\}.$$
Let $P_{\overline j}\subset G$ be the parabolic subgroup with
Lie algebra $\frp_{\overline j}.$
Let  $U_{\overline j}\subset P_{\overline j}$ be its unipotent radical and let $L_{\overline j}\subset P_{\overline j}$
be its standard Levi subgroup. The group $L_{\overline j}$ splits into a product
$$L_{\overline j}=L(\overline{j}) \times L'(\overline{j})$$ with
$L'(\overline{j})\cong {\rm GL}_{j}$ and $L(\overline{j}) \cong {\rm GL}_{d-j+1}.$
We get $$\frl_{\overline j}=\frl(\overline{j}) \times \frl'(\overline{j}).$$
In \cite{P}, chapter  1, it is shown that $\frl(\overline{j}) $ acts on $M_{\overline{j}}$ via the $j$-th symmetric power
on the $K$-vector space $K^{d+1-j}.$
Further  $\frl'(\overline{j}) $ acts on $M_{\overline{j}}$ by the trace character, i.e., by
$$\frt\cdot \Xi_\mu = (-\sum\nolimits_{i\in \overline{j}}\frt_i)\cdot \Xi_\mu ,\;\frt \in \frl'(\overline{j}).$$
Denote by $\frd_{\overline{j}}:={\rm ker}\; \varphi_{\overline{j}}\;$ the kernel of the epimorphism
$$\varphi_{\overline{j}}: U(\frg)\otimes_{U(\frp_{\overline{j}})} M_{\overline{j}} \rightarrow \fra_{\overline{j}}/\fra_{\overline{j}}^<.$$

\noindent Let $\cO_{{\rm alg}}(\cX)\subset \cO(\cX)$ be the subspace of algebraic functions.
More concretely,
\begin{eqnarray*}
\cO_{{\rm alg}}(\cX) & = & \Big\{ F=\frac{P}{Q}\mid Q= \prod_{j=1}^r (\sum_{i=0}^d c_{ij} X_i)^{l_j}\mid c_{ij} \in K, \;   \\
& & \mbox{ $P$ is a homogeneous polynomial of degree } \sum_{j=1}^r l_j  \Big\}.
\end{eqnarray*}
This is a dense subset of the $K$-Fr\'echet algebra $\cO(\cX),$ comp. \cite{ST1} 3.3.
We denote by
$$\imath: \cO_{{\rm alg}}(\cX) \rightarrow \bbN$$ the index function.
This is a $G$-invariant map taking values in the interval  $[0,d].$
For a function $F=\frac{P}{Q} \in \cO_{{\rm alg}}(\cX),$ with
 $(P,Q)=1 \mbox{ and  pairwise different } (c_{0,j},\ldots, c_{d,j})\in K^{d+1},\,  j=1,\ldots,r ,$
 it is defined by $\imath(\frac{P}{Q})=\imath_o(\frac{P}{Q})=r\;.$
  In general, it is given by
$$\imath(F)= \min \max \Big\{\imath_o(\frac{P_k}{Q_k}) \,\mid \, k\Big\},$$
where the minimum is taken over all representations $F=\sum_k \frac{P_k}{Q_k}$ with
$\frac{P_k}{Q_k} \in \cO_{{\rm alg}}(\cX).$
Put 
\begin{equation}\label{Filtration_alg}
\cO_{{\rm alg}}(\cX)_j:= \Big\{ F \in \cO_{{\rm alg}}(\cX)\mid \imath(F) \leq j \Big\}.
\end{equation}
We obtain a filtration
$$K=\cO_{{\rm alg}}(\cX)_0 \subset \cO_{{\rm alg}}(\cX)_1 \subset \cdots \subset \cO_{{\rm alg}}(\cX)_{d-1} \subset \cO_{{\rm alg}}(\cX)_d = \cO_{{\rm alg}}(\cX)$$
on $\cO_{{\rm alg}}(\cX).$
The relation between the filtration (\ref{LieFiltration}) on $\cO_{{\rm inf}}(\cX)$ and this one is
$$\cO_{{\rm alg}}(\cX)_j = \sum_{g\in G} g\cdot \fra_j, \; j=0,\ldots,d.$$
This follows from some of the  results in \cite{GV} and is explained in \cite{P} resp. \cite{ST1}.
Finally, for $j=1,\ldots,d$, let
$$\cO(\cX)_j:= \overline{\cO_{{\rm alg}}(\cX)_j} \subset \cO(\cX)$$
be the topological closure of $\cO_{{\rm alg}}(\cX)_j$ in $ \cO(\cX).$
We get a $G$-equivariant filtration by closed $K$-Fr\'echet spaces
\begin{equation}\label{Pohlkampfiltration}
K=\cO(\cX)_0 \subset \cO(\cX)_1 \subset \cdots \subset \cO(\cX)_d = \cO(\cX)
\end{equation}
on $\cO(\cX).$ Each subquotient is a reflexive $K$-Fr\'echet space with a continuous $G$-action (see \cite{ST1} Prop. 6)
and its dual is a locally analytic $G$-representation. Similarly to \cite{ST1}, Pohlkamp constructs for $1\leq j \leq d,$ isomorphisms
 $$(\cO(\cX)_j/\cO(\cX)_{j-1})' \stackrel{\thicksim}{\longrightarrow} C^{an}(G,P_{\overline{j}}\;; M_{\overline{j}}'\otimes {\rm St}_j)^{\frd_{\overline{j}}=0}$$
 of locally analytic $G$-representations. Here, the unipotent radical of $P_{\overline{j}}$ acts trivially on $M_{\overline{j}}'\otimes {\rm St}_j.$ The  group
 $L(\overline{j})$ acts in the obvious way. The action of  $L'(\overline{j})$ is given by the inverse of the determinant character and on ${\rm St}_j$ by the Steinberg representation.

We return to our computation. First of all, the cohomology of the structure sheaf $\cO$  is given by $$H^\ast(\P^d_K,\cO)=H^0(\P^d_K,\cO)=K.$$ So, in this
case all the contributions $v^{G}_{{P_{(d+j+1,1,\ldots,1)}}}(H^{-j}(\P^d_K,\cO)),\; j=-1,\ldots, -d,$ in Theorem \ref{Theorem2}  vanish.
Moreover, we have $\cO(\cX)_0=K=W^0.$ It remains to compute
the locally analytic part of our formula.
The structure sheaf corresponds to the weight $\lambda=(0,\ldots,0)\in  \Z^{d+1}.$ We get\footnote{Instead of a further comma, we use the symbol $\,\mid\,$ for a better distinction of the individual vectors.}
$$w_j\ast \lambda=(\underbrace{-1,\ldots,-1}_j\mid j,0,\ldots,0)$$ and
$$ \mu_{j,\lambda}= w_j\ast \lambda \; \mbox{ for all } j=1,\ldots,d.$$ Further we compute
$$\Phi_{j,\lambda}=\bigcup_{k=0}^j \{(\underbrace{-1,\ldots,-1}_{j-1},-1-k\mid j,k,0\ldots,0) \}$$
resp.
$$\Psi_{j,\lambda}=\bigcup_{k=0}^j \{(j,k,0\ldots,0\mid \underbrace{-1,\ldots,-1}_{j-1},-1-k) \}.$$
By Corollary \ref{irreduciblelagbraicrepr}  we deduce that the  $U(\frg)$-module $\tilde{H}^{d-j}_{ \P^j_K}(\P^d_K, \cO)=H^{d-j}_{ \P^j_K}(\P^d_K, \cO)$
is generated by a quotient of the ${\bf P_{(j+1,d-j)}}$-representation $N_{d-j,\lambda}=\bigoplus_{\mu \in \Psi_{d-j,\lambda}} V_\mu$, where $V_\mu$ is the irreducible algebraic ${\bf L_{(j+1,d-j)}}$-representations with highest weight $\mu.$
Actually, by the following proposition and the $U(\frg)$-structure with respect to $\cO$, cf. (\ref{Strukturgarbe}), it turns out that the representation  to the one with highest weight
$$z^{-1}_{d-j}\cdot \mu_{d-j,\lambda}= (d-j,0,\ldots,0\mid \underbrace{-1,\ldots,-1)}_{d-j}$$
generates $H^{d-j}_{\P_K^j}(\P_K^d,\cO)$ as $U(\frg)$-module.
\begin{Proposition}\label{cohStrukturgarbe} For $0\leq j \leq d-1,$ we have
$$H^{d-j}_{\P_K^j}(\P_K^d,\cO)=\bigoplus_{k_0,\ldots, k_j \geq 0  \atop {k_{j+1},\ldots, k_d < 0  \atop  k_0+ \ldots + k_d=0}} K\cdot X_0^{k_0}X_1^{k_1}\cdots X_d^{k_d}$$
\end{Proposition}

\proof The proof follows easily from the formula (\ref{IndLim}). \qed

Hence we can choose $N_j$ to be the ${\bf P_{(j+1,d-j)}}$-module isomorphic to the outer tensor product
$${\rm Sym}^{d-j}(K^{j+1})\boxtimes \det{}^{-1}.$$
By the proposition above, we may identify   $N_{j}$ with the $K$-vector space generated by the elements
$$\frac{P}{X_{j+1}\cdots X_{d}} \in \bigoplus_{k_0,\ldots, k_j \geq 0  \atop {k_{j+1},\ldots, k_d < 0  \atop  k_0+ \ldots + k_d=0}} K\cdot X_0^{k_0}X_1^{k_1}\cdots X_d^{k_d},$$
where $P$ is homogeneous polynomial in $X_{0},\ldots,X_{j}$ of degree $d-j.$
These are precisely the elements
$$\Xi_\mu \in \cO(\cX) \; \mbox{ where } \mu \in  A(\overline{d-j}).$$
In particular, we have
$$N_j=M_{\overline{d-j}}$$ as subspaces of $\cO_{{\rm inf}}(\cX).$
Further, the set of differential equations $\frd_j = \frd_{\overline{d-j}}$ coincide.
By  Theorem \ref{Theorem2}  we have a filtration $W^{-d} \supset W^{-d+1} \supset \cdots \supset W^{-1}\supset W^0$ on $H^0(\cX,\cO)$
with
$$(W^j/W^{j+1})' \cong  C^{an}(G,P_{(d+j+1,j)};N_{d+j}' \otimes {\rm St}_{-j})^{{\frd_{d+j}}} ,\; j=-1,\ldots,-d.$$

Hence we have shown that the graded pieces of our filtration  coincides with that in \cite{P}.
Moreover, by definition of the canonical filtration on $H^1_{\cY}(\P^d_K,\cO)$ induced by the spectral sequence (\ref{ss}) (cf. Lemma \ref{Filtration}), we see that both filtrations are the same. In fact, this follows from Proposition  \ref{cohStrukturgarbe}   and (\ref{LieFiltration}).

\subsection{$\cF = \Omega_{\P^d_K}^d$}

In this chapter we compare our result in the case $\cF = \Omega^d= \Omega_{\P^d_K}^d$ to that of Schneider and Teitelbaum. In \cite{ST1} the authors define a $G$-equivariant decreasing filtration by closed $K$-Fr\'echet spaces
 $$\Omega^d(\cX)^0 \supset \Omega^d(\cX)^1 \supset \cdots \supset \Omega^d(\cX)^{d-1} \supset \Omega^d(\cX)^d \supset \Omega^d(\cX)^{d+1}=0$$ on $\Omega^d(\cX)^0=H^0(\cX,\Omega^d).$
As in \cite{P} this construction involves an index function $\imath $  on
the algebraic differential forms $\Omega^d_{{\rm alg}}(\cX)$, which counts the negative prime divisors without multiplicities of a given differential form.
For any integer $j \in \N,$ the filtration step $\Omega^d(\cX)^{j}$ is defined by the topological closure in $H^0(\cX,\Omega^d)$ of its algebraic differential forms
$$\Omega_{{\rm alg}}^d(\cX)^{j}:=\Big\{ \eta \in \Omega_{{\rm alg}}^d(\cX)\mid \imath(\eta) \leq d+1-j \Big\}.$$
Furthermore, they construct explicit isomorphisms ({\it boundary value maps})
$$ (\Omega^d(\cX)^j/\Omega^d(\cX)^{j+1})' \stackrel{\thicksim}{\longrightarrow} C^{an}(G,P_{\underline{j}}; M_{\underline{j}}' \otimes {\rm St}_{d+1-j})^{\deju =0}$$
of locally analytic $G$-representations.
Here, $P_{\underline{j}}$=$P_{(j,d+1-j)}\subset G$ is the (lower) standard-parabolic  subgroup to the decomposition $(j,d+1-j)$
of $d+1$ and 
$$\underline{j}:=\{0,\ldots,j-1\}.$$
The $K$-vector space $M_{\underline{j}}\subset U(\frg)$ is given by the sum
$$M_{\underline{j}}= \sum_{\mu \in B(\underline{j})} K\cdot L_\mu,$$
where
$$B(\underline{j})= \Big\{ \mu= \sum\nolimits_k m_k \epsilon_k \in X^\ast(\overline{T})\,\mid\, m_k=1 \mbox{ for } k \in \underline{j},\, m_k \leq 0 \mbox{ for } k \notin \underline{j} \Big\}.$$
Further, the symbol $L_\mu$ denotes a (sorted) element of $U(\frg)$, cf. \cite{ST1} p. 31. It is of weight $\mu$ and satisfies $L_\mu\cdot \xi= - \Xi_\mu\cdot \xi,$
where
$$\xi =  \frac{X_0^d}{X_1\cdots X_d}\cdot d(\frac{X_1}{X_0})\wedge d(\frac{X_2}{X_0})\wedge \cdots \wedge  d(\frac{X_d}{X_0}).$$
The unipotent radical of $P_{\underline{j}}$ acts trivially on the tensor product $M_{\underline{j}}' \otimes {\rm St}_{d+1-j}$.
The second factor of the  Levi subgroup $L_{(j,d+1-j)}= {\rm GL}_j \times {\rm GL}_{d+1-j}$ acts on $M_{\underline{j}}'$ via the symmetric power ${\rm Sym}^j(K^{d+1-j}).$ On ${\rm St}_{d+1-j}$ it acts via the Steinberg representation. The action of ${\rm GL}_j$ is given by the inverse of the determinant character.   In particular the case $j=0$ yields the Steinberg representation ${\rm St}_{d+1}$.

We turn to our computation.
We have
$$H^\ast(\P^d_K,\Omega^d)=H^d(\P^d_K,\Omega^d)=K\cdot \xi.$$
Consequently, all the contributions $v^{G}_{{P_{(d+j+1,1,\ldots,1)}}}(H^{-j}(\P^d_K,\Omega^d)),\; j=-1,\ldots, -d,$ of Theorem \ref{Theorem2}   vanish except for
$j=-d.$ In the latter case we obtain the Steinberg representation $v^{G}_{{P_{(1,1,\ldots,1)}}}(K)={\rm St}_{d+1}.$
The canonical bundle corresponds to the homogeneous vector bundle $\cF_\lambda$ of weight $\lambda=(-d,1,\ldots,1)\in  \Z^{d+1},$ cf. \ref{Exampleirred}.
We have
$$w_j\ast \lambda=(0,\ldots,0,-d+j\mid\underbrace{1,\ldots,1}_{d-j})\;\; \mbox{ for } j=1\ldots,d.$$
It follows that
$$\mu_{j,\lambda}= w_{j-1}\ast \lambda\;\; \mbox{ for } j=1,\ldots,d$$
and
$$\Phi_{j,\lambda}=\{\mu_{j,\lambda} \}.$$
We deduce that  $\tilde{H}^{d-j}_{\P^{j}_K}(\P^d_K, \Omega^d)$
is generated as  $U(\frg)$-module by the ${\bf P_{(j+1,d-j)}}$-represen\-tation $N_{j}=N_{d-j,\lambda}$ corresponding to the irreducible ${\bf L_{(j+1,d-j)}}$-representation with highest weight
$$z^{-1}_{d-j}\cdot \mu_{d-j,\lambda}= (\underbrace{1,\ldots,1}_{j+1}\mid 0,\ldots,0,-j-1).$$
It follows that the ${\bf P_{(j+1,d-j)}}$-module $N_j$ is isomorphic to
$$\det \;\;\boxtimes \;\; {\rm Sym}^{j+1}(K^{d-j})'.$$
Again, we want to realize the corresponding representation concretely.

\begin{Proposition} For $j \geq 0,$ we have
$$\tilde{H}^{d-j}_{\P^j_K}(\P_K^d,\cO(-d-1))=\bigoplus_{k_0,\ldots, k_j \geq 0  \atop {k_{j+1},\ldots, k_d < 0  \atop  k_0+ \ldots + k_d=-d-1}} K\cdot X_0^{k_0}X_1^{k_1}\cdots X_d^{k_d}.$$
\end{Proposition}

\proof The proof is similar to that of Proposition \ref{cohStrukturgarbe}  \qed

By identifying $\Omega^d$ with the twisted sheaf $\cO(-d-1),$ the element $\xi$ corresponds to the fraction $\frac{1}{X_0\cdots X_d}.$
It follows that $N_{j}$ is the $K$-vector space generated by the elements
$$\frac{X_0\cdots X_j}{X_{j+1}^{m_{j+1}}\cdots X_d^{m_d}} \cdot \xi,$$
with $m_{j+1} + \cdots + m_d = j+1.$
The fractions $\frac{X_0\cdots X_j}{X_{j+1}^{m_{j+1}}\cdots X_d^{m_d}}$ are exactly the elements $\Xi_\mu$, with
$\mu \in B(\underline{j+1}),$ cf. \cite{ST1}, p. 65. We get for $j=0,\ldots,d-1,$ isomorphisms
$$M_{\underline{j+1}} \stackrel{\sim}{\longrightarrow} N_{j}$$
$$L_\mu \mapsto L_\mu \cdot \xi = - \Xi_\mu \cdot  \xi .$$

Moreover, the set of differential equations $\frd_{j}$ and $\frd_{\underline{j+1}}$ are the same.
By Theorem \ref{Theorem2}  we have a filtration $W^{-d} \supset W^{-d+1} \supset \cdots \supset W^{-1} \supset W^0$ on $H^0(\cX,\Omega^d)$
where $(W^{-d}/W^{-d+1})'$ is an extension
\begin{eqnarray}\label{ExakteSequenz}
0\rightarrow  v^{G}_{{P_{(1,1,\ldots,1)}}}(K) \rightarrow (W^{-d}/W^{-d+1})' \rightarrow C^{an}(G,P_{(1,d)};N_0'\otimes {\rm St}_d)^{{\frd_0}}\rightarrow 0
\end{eqnarray}
and
$$(W^{j}/W^{j+1})' \cong C^{an}(G,P_{(d+1+j,-j)};N_{d+j}' \otimes {\rm St}_{-j})^{{\frd_{d+j}}} .$$
Thus we see that the graded pieces of our filtration coincide with that of \cite{ST1}.
Moreover, by looking at the pole order of sections in $\Omega^d(\cX)$ as in the case of the structure sheaf we see that the filtrations are the same apart from the first filtration step. The difference is just given by the extension above.
In other words, we have an extension
$$0\rightarrow  \Omega^{d}(\cX)^1/\Omega^{d}(\cX)^2 \rightarrow W^{-d}/W^{-d+1} \rightarrow \Omega^{d}(\cX)^0/\Omega^{d}(\cX)^1 \rightarrow 0,$$
such that its dual coincides with (\ref{ExakteSequenz}).

\subsection{$\cF = \Omega_{\P^d_K}^1$}

This chapter provides another example for our computation. It treats the cotangent bundle $\cF = \Omega^1=\Omega_{\P^d_K}^1$ on $\P^d_K.$

We have  $$H^\ast(\P^d_K,\Omega^1)=H^1(\P^d_K,\Omega^1)=K.$$ Therefore, all the contributions $v^{G}_{{P_{(d+j+1,1^{-j})}}}(H^{-j}(\P^d_K,\Omega^d))$ in Theorem \ref{Theorem2}  vanish except for
$j=-1.$ In the latter case we obtain the generalized Steinberg representation $v^{G}_{P_{(d,1)}}(K).$
The cotangent bundle corresponds to the homogeneous vector bundle $\cF_\lambda$ given by the weight $\lambda=(-1,1,0,\ldots,0)\in  \Z^{d+1}.$
By Theorem \ref{Theorem2}  we have a filtration $W^{-d} \supset W^{-d+1} \supset \cdots \supset W^{-1} \supset W^0$ on $H^0(\cX,\Omega^1)$
where $(W^{-1}/W^{0})'$ is an extension
$$0\rightarrow  v^{G}_{{P_{(d,1)}}}(K) \rightarrow (W^{-1}/W^{0})' \rightarrow C^{an}(G,P_{(d,1)};N_{d-1}'\otimes {\rm St}_1)^{{\frd_{d-1}}}$$
and
$$(W^{j}/W^{j+1})' \cong C^{an}(G,P_{(d+1+j,-j)};N_{d+j}' \otimes {\rm St}_{-j})^{{\frd_{d+j}}} $$
for $j\neq -1.$
A computation shows that
$$\mu_{j,\lambda}:=\left\{ \begin{array}{cc} \lambda  &: j = 1 \\ \\ (0,\underbrace{-1,\ldots,-1}_{j-1}\mid j-1,0\ldots,0) & : j > 1. \end{array} \right.$$
Further, for $j>1$
$$\Phi_{j,\lambda}=\{\mu_{j,\lambda}\} \cup \bigcup_{k=1}^{j-1}\{(l,-1,\ldots,-1,-1-l-k\mid j-1,k,0\ldots,0)\mid  l=0,-1 \}$$
resp.
$$\Psi_{j,\lambda}=\{z_i^{-1}\cdot \mu_{j,\lambda}\} \cup \bigcup_{k=1}^{j-1}\{(j-1,k,0\ldots,0\mid l,-1,\ldots,-1,-1-l-k )\mid  l=0,-1 \}.$$
In the case $j=1$, we compute
$$\Phi_{1,\lambda}=\{\mu_{1,\lambda}, (-2\mid 1,1,0\ldots,0) \}$$
resp.
$$\Psi_{1,\lambda}=\{z_i^{-1}\cdot \mu_{1,\lambda}, (1,1,0\ldots,0\mid-2) \}$$
We will see that for $n>1$  the weights with $k\leq 1$ yield a generating system of $\tilde{H}^{d-j}_{\P^{j}_K}(\P^d_K, \Omega^1)$.

We deduce that  $\tilde{H}^{d-j}_{\P^{j}_K}(\P^d_K, \Omega^1)$
contains the irreducible algebraic ${\bf L_{(j+1,d-j)}}$-representation $V_\mu$ with highest weight
$$\mu=z^{-1}_{d-j}\cdot \mu_{d-j,\lambda}:=\left\{ \begin{array}{cc} (1,0,\ldots,0\mid -1)  &: j=d-1  \\ \\ (d-j-1,0\ldots,0\mid 0,\underbrace{-1,\ldots,-1}_{d-j-1}) & : j < d-1. \end{array} \right.$$
It follows that $V_{z^{-1}_{d-j}\cdot \mu_{d-j,\lambda}}$ is isomorphic to
$${\rm Sym}^{d-j-1}(K^{j+1})  \boxtimes (K^{d-j} \otimes \det{}^{-1}) \mbox{ for } j< d-1.$$
For $j=d-1,$ we get $V_{z^{-1}_1  \cdot \mu_{1,\lambda} }\cong  K^d \;\boxtimes\; \det{}^{-1}.$
We shall give an explicit  realization of $V_{z^{-1}_{d-j}\cdot \mu_{d-j,\lambda}}$.

Let $V_{j}$ be the finite-dimensional $K$-vector space generated by the elements
 $$\frac{P\cdot X_k^2}{X_{j+1}\cdots X_d}\cdot d(\frac{X_l}{X_k}),\; k\in\{0,\ldots,j\}, l \in
\{j+1,\ldots,d\},$$
where $P$ is a homogeneous polynomial  of degree $d-j-2$ in the indeterminates $X_0,\ldots, X_{j}.$ In the case $j=d-1$, let $V_{d-1}$ be generated by the elements $d(\frac{X_k}{X_d}),$ $k\in \{0,\ldots,d-1 \}.$
Consider the $K$-linear map
\begin{eqnarray*}
{\rm Sym}^{d-j-2}(K^{j+1}) \otimes K^{j+1} \otimes K^{d-j} & \longrightarrow & V_{j} \\
 P\otimes e_k\otimes e_l & \mapsto & \frac{P\cdot X_k^2}{X_{j+1}\cdots X_d}\cdot d(\frac{X_l}{X_k}).
\end{eqnarray*}
This map  is clearly a surjective linear map of $K$-vector spaces.
Consider the following action of ${\bf L_{(j+1,d-j)}}$ on the LHS.
On $K^{j+1}$ the action of ${\bf L(j+1)}$ is the standard representation.
On ${\rm Sym}^{d-j+2}(K^{j+1})$ it acts via the $(d-j+2)$-th symmetric power.
On $K^{d-j}$ the operation of ${\bf L(j+1)}$ is the trivial one.
The action of ${\bf L(d-j)}$ on $K^{d-j}$ is the standard representation.
On the other factor it operates via the inverse of the determinant character.
From the identity
$$\frac{P\cdot X_k^2}{X_{j+1}\cdots X_d}\cdot d(\frac{X_l}{X_k}) = - \frac{P\cdot X_l^2}{X_{j+1}\cdots X_d}\cdot d(\frac{X_k}{X_l}),\; k\in\{0,\ldots,j\}, l \in \{j+1,\ldots,d\},$$
it follows that $V_j$ is a finite-dimensional  ${\bf L_{(j+1,d-j)}}$ -module and that the map above is a surjection of ${\bf L_{(j+1,d-j)}}$-representations.
Inside the representation ${\rm Sym}^{d-j-2}(K^{j+1}) \otimes K^{j+1}$ we have the irreducible ${\bf L(j+1)}$-subrepresen\-tation ${\rm Sym}^{d-j-1}(K^{j+1}),$
which corresponds to the highest weight $(d-j-1,0,\ldots,0)$ of ${\bf L(j+1)}.$ A computation shows that the above map restricts
to an isomorphism
$$ {\rm Sym}^{d-j-1}(K^{j+1}) \boxtimes (K^{d-j} \otimes \det{}^{-1}) \stackrel{\sim}{\longrightarrow} V_{j}.$$

For $j<d-2$, the representation $V_j=V_{z^{-1}_{d-j}\cdot \mu_{d-j,\lambda}}$ is a quotient of some representation containing  the representation ${\rm Sym}^{d-j-2}(K^{j+1}).$ We have realized the latter one as a subrepresentation of $K[X_0,\ldots, X_j]$.
Thus $K[X_0,\ldots, X_j] \cdot ({\rm Sym}^{d-j-2}(K^{j+1}) \otimes K^{j+1} \otimes K^{d-j}) \subset K[X_0,\ldots, X_j] \cdot (K^{j+1} \otimes K^{d-j}).$
By the discussion in 1.4 it suffices to consider in Lemma \ref{Lemma_boxtimes} the representations $K^{j+1}\otimes K^{d-j}$ instead of ${\rm Sym}^{d-j-2}(K^{j+1}) \otimes K^{j+1} \otimes K^{d-j}$. Since the highest weight of $K^{j+1} \otimes K^{d-j}$ is 
$(1,0\ldots,0\mid 0,-1,\ldots,-1)$ we deduce that the weights in $z_{d-j}^{-1} \Phi_{d-j,\lambda}$ with $k\leq 1$ yield a generating system. 

As for  the other irreducible subrepresentations, we  note that in the case $j=d-1$ the irreducible ${\bf P}_{(d,1)}$-representation corresponding to the weight
$(1,1,0,\ldots,0\mid -2)$ is generated by the expressions 
$$\frac{X_i}{X_d}d(\frac{X_j}{X_d}) - \frac{X_j}{X_d}d(\frac{X_i}{X_d}).$$
For $j<d-1,$ we write down a highest weight vector  for the remaining  weights.
In the case of 
$(d-j-1,1,0\ldots,0 \mid -1,-1,\ldots,-1,-1)$ it is given by  
$$\frac{X_0^{d-j-1}}{X_{j+2}\cdots X_d} d(\frac{X_1}{X_{j+1}}) - \frac{X_0^{d-j-2}X_1}{X_{j+2}\cdots X_d} d(\frac{X_0}{X_{j+1}}).$$
In the case of
$(d-j-1,1,0\ldots,0 \mid 0,-1,\ldots,-1,-2)$ it is given by $$\frac{X_0^{d-j-2}X_{j+1}^2}{X_{j+1}\cdots X_d}\big(\frac{X_0}{X_d} d(\frac{X_1}{X_{j+1}}) - \frac{X_1}{X_d}d(\frac{X_0}{X_{j+1}})\big).$$

\newpage
\section{Appendix}
In this final part of our paper we present another approach for the main computation.
We replace the acyclic complex of Theorem \ref{Theorem1} by a similar complex avoiding
adic spaces. It is based purely on rigid analytic varieties. The construction
is similar to that of \cite{SS} in the case of constant coefficients. In that case, the main difference is essentially that in loc.cit. the authors regard intersections of hyperplanes,  whereas we deal
with arbitrary subspaces directly. Our approach follows the construction of \cite{O3}.
For this purpose, we have further to investigate the neighborhoods of the closed subvarieties $Y_U=\P(U) \subset \P^d_K,$ where $U\subset K^{d+1}$ is a linear subspace.

\setcounter{subsection}{+1}

Let
$$\Lambda=\bigoplus_{i=0}^d O_K\cdot  e_i$$
be the $O_K$-lattice of $V$ generated by our fixed basis $\{e_0,\ldots,e_d\}$.
Let $O$ be the ring of integers in $\Cp.$ For $n\in \N$, we put 
$$O_K^{(n)}:=O_K/\pi^nO_K,\; O^{(n)}:=O/\pi^nO,\;  \Lambda^{(n)}:=\Lpn. $$
Consider for $n\in \N,$ the mod-$n$ reduction map
$${\rm red}_n:\P_K^d(\Cp) \longrightarrow \P(\Lambda)(O^{(n)}).$$
If $L_x \subset \Cp^{d+1}$ is a line representing a point $x \in \P^d_K(\Cp),$
then
$${\rm red}_n(x)= (L_x \cap (\Lambda \otimes_{O_K} O)) \otimes O^{(n)}.$$
Let $I_U \subset K[T_0,\ldots,T_d]$  be  the vanishing ideal of  the linear subvariety $Y_U.$
The schematic closure of $Y_U$ in $\mathbb
P^d_{O_K}$ is defined by the the ideal
$$\tilde{I}_U = I_U \cap O_K[T_0,\ldots,T_d].$$
Consider the $\epsilon_n$-neighborhood $Y_U(\epsilon_n)$
 for a given non-trivial $K$-subspace $U\subsetneq V.$
Let $f \in \tilde{I}_U$  be a polynomial, such that at least one coefficient is a unit. Then we have for $x\in \P_K^d(\Cp)$ (take any unimodular representative of $x$),
\begin{eqnarray*}
x\in Y_U(\epsilon_n)(\Cp) & \aequi & |f(x)|\leq \epsilon_n \\
                        & \aequi & \pi^n  \mid f(x) \\
                        & \aequi & f(x) = 0 \mbox{ (mod } \pi^n).
                        \end{eqnarray*}
So $Y_U(\epsilon_n)$ is nothing else but
$$ Y_U(\epsilon_n)=\Big\{ x
\in (\P_K^d)^{rig}\mid {\rm red}_n(x) \in \tilde{Y}_U(O^{(n)}) \Big\},$$
where $\tilde{Y}_U$ is the closed subscheme of $\P^d_{O_K}$ defined
by $\tilde{I}_U.$ Moreover, we see that
\begin{eqnarray}\label{tubes=}
Y_U(\epsilon_n)=Y_{U'}(\epsilon_n)
\end{eqnarray}  if
$U \cap \Lambda \equiv U'\cap \Lambda$ mod $\pi^n,$ compare also Lemma 2, ch. 1 in \cite{SS}.

Let $$U_\bullet=\big((0) \subsetneq U_1 \subsetneq U_2  \subsetneq \cdots \subsetneq U_s \subsetneq \Lambda^{(n)}\big)$$ be a filtration of free  $O_K^{(n)}$-modules.  Choose a lift of $U_\bullet$ to a $K$-filtration
$$\tilde{U}_\bullet=\big((0) \subsetneq \tilde{U}_1 \subsetneq \tilde{U}_2 \subsetneq \dots \subsetneq \tilde{U}_s \subsetneq V\big)$$ of $K$-subspaces of $V.$
We set
$$Y_{U_\bullet}:= Y_{\tilde{U}_1}(\epsilon_n).$$
This definition depends by (\ref{tubes=}) only on the $O_K^{(n)}$-module $U_1.$ Hence, the set $Y_{U_\bullet}$
is a well-defined rigid analytic open subvariety of $(\P^d_K)^{rig}$. Thus we have associated to
each filtration $U_\bullet$ of  $\Lambda^{(n)}$ by free $O_K^{(n)}$-submodules  a quasi-compact rigid analytic subset $Y_{U_\bullet}\subset (\P_K^d)^{rig}.$
For $0\leq i \leq d,$ let
$$\Ln_i=\bigoplus_{j=0}^i \,O_K^{(n)} \cdot e_j \subset \Lambda^{(n)}$$
be the free  $O_K^{(n)}$-submodule generated by our first $i+1$  basis vectors.
For any subset  $I= \{\alpha_{i_1},\ldots,\alpha_{i_r}\} \subset \Delta,$ we put
$$\Ln_I= \big((0) \subset \Ln_{i_1} \subset \Ln_{i_2} \subset \dots \subset \Ln_{i_r}\subset \Lambda^{(n)}\big).$$
Then  we get $Y_{\Ln_I} = Y_I(\epsilon_n).$
Thus, analogously to (\ref{complex}),  we can
construct for every \'etale (resp. Zariski) sheaf $\cG$ on ${\cY}_{n}$ the following
complex of \'etale (resp. Zariski) sheaves on ${\cY}_n$:
\begin{multline}\label{complex2}
0 \rightarrow \cG \rightarrow\!\!\! \bigoplus_{I \subset \Delta \atop
|\Delta\setminus I|=1} \bigoplus_{g \in G_0/P^n_I}
(\phi^n_{g,I})_*(\phi^n_{g,I})^* \cG \rightarrow \!\!\!\bigoplus_{I \subset
\Delta \atop |\Delta \setminus I|=2} \bigoplus_{g \in
G_0/P^n_I}(\phi^n_{g,I})_*(\phi^n_{g,I})^* \cG \rightarrow \\
\dots \rightarrow \!\!\!\bigoplus_{I \subset \Delta \atop |\Delta \setminus I|=d-1}
\bigoplus_{g \in G_0/P^n_I} (\phi^n_{g,I})_*(\phi^n_{g,I})^* \cG \rightarrow
\bigoplus_{g \in G_0/P^n_\emptyset}
(\phi^n_{g,\emptyset})_*(\phi^n_{g,\emptyset})^* \cG \rightarrow 0,
\end{multline}
where $\phi^n_{g,I}$ denotes the open embedding $Y_{g\Ln_I} \hookrightarrow
{\cY}_n$  of rigid analytic varieties.
In contrast to $(\ref{complex})$, this complex is not acyclic as the following example shows. This was pointed out to me by P. Schneider
some years ago.

\medskip
\begin{Example}
{\rm Let $d=2.$ Then the above complex is nothing else but
$$  0 \rightarrow \cG \rightarrow\!\!\! \bigoplus_{g \in G_0/P^n_{\{\alpha_1\}}} (\phi^n_{g,\{\alpha_1\}})_*(\phi^n_{g,\{\alpha_1\}})^* \cG \;\; \oplus  \bigoplus_{g \in
G_0/P^n_{\{\alpha_2\}}}(\phi^n_{g,\{\alpha_2\}})_*(\phi^n_{g,\{\alpha_2\}})^* \cG  $$
$$\rightarrow \bigoplus_{g \in G_0/P_\emptyset^n}
(\phi^n_{g,\emptyset})_*(\phi^n_{g,\emptyset})^* \cG \rightarrow 0.$$
We have
\begin{eqnarray*}
Y_{\{\alpha_1\}}(\epsilon_n)(\Cp) & = & \big\{ x\in \mathbb P^d_K(\Cp)\mid {\rm red}_n(x)=\Ln_1\otimes_{O_K^{(n)}} O^{(n)} \big\} \\ \\
Y_{\{\alpha_2\}}(\epsilon_n)(\Cp) & = & \big\{ x\in \mathbb P^d_K(\Cp)\mid {\rm red}_n(x)\subset \Ln_2\otimes_{O_K^{(n)}} O^{(n)}\big\} \\ \\
Y_\emptyset(\epsilon_n)(\Cp)  & = & Y_{\{\alpha_1\}}(\epsilon_n)(\Cp).
\end{eqnarray*}
Let $x\in \P^2_K(\Cp)$  be a point such that the corresponding line $L_x \subset \Cp^{3}$ has the shape
$L_x=\Cp\cdot(\pi^{n-1}e_0+ae_2), \; a \in O^\times,$ with $[a]\in (O^{(n)})^\times \setminus (O_K^{(n)})^\times.$
 Consider the planes
$$E=O_K^{(n)}\cdot (e_0 + \pi e_1) \oplus O_K^{(n)}\cdot e_2 \mbox{ and } E'= O_K^{(n)}\cdot e_0 \oplus O_K^{(n)}\cdot e_2$$
in $\Lambda^{(n)}.$ Then
$${\rm red}_n(x)\in E\otimes_{O_K^{(n)}} O^{(n)}\; \cap\; E'\otimes_{O_K^{(n)}} O^{(n)},$$ but ${\rm red}_n(x)\neq L\otimes_{O_K^{(n)}} O^{(n)}$ for all $L \in \mathbb P(\Lambda)(O_K^{(n)}).$ So, localizing the above complex in $x$ yields a sequence
$$0\rightarrow \cG_x\rightarrow \cG_x^r  \rightarrow 0$$
with $r=\#\{E\subset \Lambda^{(n)}\mid\,E \mbox{ is free with }
{\rm rk}(E)=2, \, {\rm red}_n(x)\subset E \otimes_{O_K^{(n)}} O^{(n)}\} \geq 2.$
Hence, the complex (\ref{complex2}) cannot be acyclic in general.}
\qed
\end{Example}

The above example suggests to fill the gaps in (\ref{complex2}).
Let $U$ be a $O_K^{(n)}$-submodule of $\Lambda^{(n)},$ not necessarily free. We define
$$Y_U:=\Big\{x \in (\P^d_K)^{rig}\mid {\rm red}_n(x) \in U \otimes_{O_K^{(n)}} O^{(n)} \Big\}.$$
If $U$ is a free $O_K^{(n)}$-module then this definition coincides with the previous one.
Put
$${\rm rk}(U)=\dim_{O_K/\pi O_K} (U+\pi\Ln/\pi\Ln).$$
This is just the rank of the torsion free submodule of $U.$
We also have the ordinary rank
$${\rm rk}'(U):=\min \{n \in \N\mid \mbox{ there are } m_1,\ldots,m_n \in U \mbox{ which generate } U\}. $$
We have ${\rm rk}(U) \leq {\rm rk}'(U)$ where the equality holds  if and only if $U$ is free ($\Leftrightarrow$ torsion free).

\begin{Proposition}\label{Tuben}
The set $Y_U$ is a quasi-compact rigid analytic open subset of $(\P_K^d)^{rig}.$
\end{Proposition}

\proof Compare also Lemma 5 in \cite{SS}. Without loss of generality, we may suppose that $U$ is generated by
$$e_0,e_1,\ldots, e_i,\pi^{n_{i+1}}e_{i+1},\ldots,\pi^{n_{j}}e_j $$
for certain integers $0 < n_{k} < n, \, k=i+1,\ldots,j.$ Then $Y_U$ is just the quasi-compact subset
$$\Big\{x=[x_0:\cdots:x_d] \in (\P_K^d)^{rig}\mid   |x_k| \leq |\pi^n|, k=1,\ldots i,\;   |x_k| \leq |\pi^{n-n_k}| , \;k=i+1,\ldots,j \Big\}$$
of $(\P^d_K)^{rig}.$
\qed

\medskip
Let $U_\bullet=\big((0)\subsetneq U_1 \subsetneq U_2 \subsetneq \dots \subsetneq U_s \subsetneq \Lambda^{(n)}\big)$ be a filtration of $\Lambda^{(n)}$ by $O_K^{(n)}$-submodules. As in the case of filtrations consisting of free $O_K^{(n)}$-submodules, we put
$$Y_{U_\bullet}= Y_{U_1}.$$
Note that $Y_{U_\bullet} = \emptyset$ is the empty set unless ${\rm rk}(U) \geq 1.$
With these newly-created rigid analytic subsets of $(\P_K^d)^{rig}$ we modify the complex (\ref{complex2}) as follows:\medskip
\begin{eqnarray}\label{complex3}
\nonumber 0 & \rightarrow & \cG \rightarrow\!\!\! \bigoplus_{{U_\bullet= (0)\subsetneq U \subsetneq \Lambda^{(n)} \atop {\rm rk}'(U) \leq d}\atop {\rm rk}(U)\geq 1 }
(\phi^n_{U_\bullet})_*(\phi^n_{U_\bullet})^* \cG \rightarrow \!\!\!\bigoplus_{{U_\bullet = (0)\subsetneq  U_1 \subsetneq U_2 \subsetneq \Lambda^{(n)} \atop {\rm rk}'(U_2) \leq d }\atop {\rm rk}(U_1)\geq 1}(\phi^n_{U_\bullet})_*(\phi^n_{U_\bullet})^* \cG \rightarrow \\ \\ \nonumber & \dots &
\rightarrow \!\!\!\!\!\!\!\!\!\!\!\!\bigoplus_{{U_\bullet= (0) \subsetneq U_1 \subsetneq \ldots \subsetneq U_p \subsetneq \Lambda^{(n)} \atop {\rm rk}'(U_p) \leq d
}\atop {\rm rk}(U_1)\geq 1 }\!\!\!\!\!\!\!\!\!\!\!\!(\phi^n_{U_\bullet})_*(\phi^n_{U_\bullet})^* \cG \rightarrow \ldots  \rightarrow\!\!\!\!\!\!\!\!\!\!\!\!\!\!\!\bigoplus_{{U_\bullet= (0) \subsetneq U_1 \subsetneq \ldots \subsetneq U_{n\cdot (d-1)} \subsetneq \Lambda^{(n)} \atop {\rm rk}'(U_{n\cdot (d-1)})\leq d} \atop {\rm rk}(U_1)\geq 1}\!\!\!\!\!\!\!\!\!(\phi^n_{U_\bullet})_*(\phi^n_{U_\bullet})^* \cG\rightarrow 0,
\end{eqnarray} where
$\phi^n_{U_\bullet}$ denotes the open embedding $Y_{U_\bullet} \hookrightarrow {\cY}_n.$
Why we impose  on $U$  the condition ${\rm rk}'(U) \leq d$ will become clear later on, cf. Prop. \ref{Homotopieaequiv}.
It simply  means that $U$ is contained in a proper free submodule of $\Lambda^{(n)}.$

As for the following theorem, we refer to \cite{JP} for the notion of an overconvergent sheaf on a rigid analytic variety. The crucial property of such a sheaf is that it vanishes if and only if all its stalks vanish.

\begin{Theorem}\label{Theorem4} Let $\cG$ be an overconvergent \'etale sheaf on
${\cY}_{n}.$ Then the complex (\ref{complex3}) is acyclic.
\end{Theorem}

The proof of this theorem is similar to Satz 5.3 in \cite{O1}. The main difference is that we are now dealing with modules instead of vector spaces. For proving Theorem \ref{Theorem4}, we have to introduce some more notation.

Let $X=(X,\prec)$ be a partially ordered set (poset). We associate to $X$ a
simplicial  complex
$$ X^\bullet = \bigcup_{n \in \mathbb N} X^n,$$
where a $n$-simplex $\tau \in X^n$ is given by a
$n+1$-tuple
$$\tau=(x_0 \prec x_1 \prec \cdots \prec x_n)$$
with elements $x_i \in X, \; i=0,\ldots, n.$ In particular, for the 0-simplices $X^0$ of $X^\bullet$ we have $X^0=X.$

A morphism $f: X \longrightarrow Y$ of posets is a map which preserves the order. Such a morphism
induces  a simplicial map
$$f^\bullet :X^\bullet \longrightarrow Y^\bullet$$
of simplicial complexes. Thus we obtain a functor from the category of posets to the category of simplicial complexes.

\begin{Example} a) \rm  Let $R=O_K^{(n)},O_K,K$ resp. $M=\Ln,\Lambda, V.$
Put
$$T_R(M):=\Big\{ R\mbox{-submodules of $M$}\mid \; {\rm rk}(U)\geq 1 \mbox{ and } {\rm rk}'(U) \leq d\Big\}.$$
We supply this set with the structure of a poset by considering the canonical order given by inclusion.
Thus a
$n$-simplex $\tau \in T_R(M)^n$ is  a flag
$$\tau =\big(\,(0) \subsetneq U_0 \subsetneq U_1 \subsetneq \cdots \subsetneq U_n \subsetneq M \,\big)$$
of submodules with ${\rm rk}(U_0) \geq 1$    and ${\rm rk}'(U_n) \leq d. $\\ \\
b) In the situation above, let $T^f_R(M)$ be the subposet consisting of all non-trivial  $R$-modules such that $M/U$ is free. We get a morphism of posets $T_R^f(M) \hookrightarrow T_R(M).$ \\ \\
c) Let $R=K.$ Then $T_K(V)$ is nothing else but the Tits complex of  ${\rm GL}(V),$ compare also \cite{Q}. \qed
\end{Example}

\begin{Proposition}\label{Homotopieaequiv} 
The morphism of posets
\Abb{\psi:}{T_K(V)}{T_{O_K}(\Lambda)}{W}{W \cap \Lambda}
induces a homotopy equivalence
$$T_{O_K}(\Lambda)^\bullet \simeq T_K(V)^\bullet$$
\end{Proposition}

\proof Using Proposition 1.6 of \cite{Q} it is enough to show that
for every proper submodule  $U\in T_{O_K}(\Lambda)$  the subposet
$\{{W \in T_K(V)\mid U \subset W\cap \Lambda}\}$ is contractible.
Note that these subsets are non-empty since ${\rm rk}'(U) \leq d.$ The contractibility
follows easily from the next proposition. \qed

\begin{Proposition}\label{Quillen} (\cite{Q}, 1.5) Let $X$ be a  poset and let $x_0 
\in X$ be a fixed element. Further, let  $f: X \longrightarrow X$  be an endomorphism of posets  with $$x\succeq f(x) \preceq x_0 \;\; \mbox{, for all } x \in X.$$ Then $X^\bullet$ is contractible.
\end{Proposition}

\begin{Remark} {\rm
It is easily seen that we may identify  $T_K(V)$ with $T^f_{O_K}(\Lambda).$ Under this identification the morphism $\psi$ corresponds to the inclusion $T^f_{O_K}(V) \hookrightarrow  T_{O_K}(\Lambda).$ \qed}
\end{Remark}

We continue with the proof of Theorem \ref{Theorem4}.

{\bf \noindent  Proof of Theorem \ref{Theorem4}:}  Since $\cG$ is overconvergent and all the
morphisms $\phi^n_{U_\bullet}$ are quasi-compact, we conclude (cf. \cite{JP} 3.5)
that all the appearing \'etale sheaves in the complex are
overconvergent, as well. Thus it is enough to show the acyclicity of the
localized complex with respect to any \'etale point  of ${\cY}_{n},$ cf.
loc.cit. 3.4. So let $e$ be an \'etale point of ${\cY}_{n}.$ By
definition this is just a separable closure $H_e$ of some valued field
$F_a$ depending on an analytic point $a$ lying below $e.$ Let $F_e$ be the
completion of $H_e.$ The \'etale point $e$ corresponds to a general
morphism $${\rm Spm }(F_e) \longrightarrow {\cY}_{n}$$ of rigid analytic varieties,
that is, to a morphism
 $${\rm Spm}(F_e) \longrightarrow {\cY}_{n}\hat{\otimes}_K F_e$$
of rigid analytic varieties,  hence to a line $L_e \in \P_K^d(F_e).$
Localizing the complex (\ref{complex3}) in $e$ yields a chain
complex with values in $F_e.$ The pull back of the complex (\ref{complex3}) to
$\P^d_K \times_{K}{\rm Spm}(F_e)$ is just the complex, where all the
appearing objects of (\ref{complex3}) are defined  with respect to the base field
$F_e.$ Localizing of this complex in $e$  would give the same chain
complex. Hence we can assume without loss of generality that $F_e =\Cp.$
The chain complex is  induced by  a subcomplex $C^\bullet$ of $T^\bullet_{O_K^{(n)}}(\Lambda^{(n)})$, which is generated by its $0$-dimensional simplices $C^0$, a subposet of $T_{O_K^{(n)}}(\Lambda^{(n)})$. Its simplices are given by
$$C^\bullet=\Big\{ U_\bullet \in T^\bullet_{O_K^{(n)}}(\Lambda^{(n)})\mid {\rm red}_n(L_e) \in Y_{U_\bullet}\Big\}.$$ By the next lemma, we will see that $C^\bullet$ is contractible. Theorem \ref{Theorem4}  follows, since the  \'etale point was arbitrary.
\qed

\begin{Lemma} The simplicial complex $C^\bullet$ is contractible
\end{Lemma}
\proof  Put $\overline{T_{O_K^{(n)}}(\Lambda^{(n)})}=T_{O_K^{(n)}}(\Lambda^{(n)})\cup \{0,\Lambda^{(n)}\}$ and supply this set with the canonical order. Thus we have realized $T_{O_K^{(n)}}(\Lambda^{(n)})$ as a subposet of $\overline{T_{O_K^{(n)}}(\Lambda^{(n)})}.$ Let $U_0\in T_{O_K^{(n)}}(\Lambda^{(n)})$ be a fixed element.
Consider the map
\Abb{f:}{C^0}{\overline{T_{O_K^{(n)}}(\Lambda^{(n)})}\;\;.}{U}{U_0 \cap U}
This map is a morphism of posets with
$$ U \succeq f(U) \preceq U_0 \;\; \mbox{ for all } \,U \in C^0     .$$
\noindent By  Proposition \ref{Quillen}, it suffices to show that the image of $f$ is contained in $C^0.$
But this follows from the inclusion
$$Y_U \cap Y_{U_0} \subset Y_{U\cap U_0},$$
which itsself follows from the flatness of $O_K^{(n)} \hookrightarrow O^{(n)} .$ \qed

\bigskip
Let $\cF$ be our fixed homogeneous vector bundle on $\P^d_K$. 
Consider the projective limit of posets $\varprojlim\nolimits_n T_{O_K^{(n)}}(\Lambda^{(n)})$ which identifies with $T_{O_K}(\Lambda).$
Let $$(U^n_\bullet)_{n\in \N} \in\varprojlim\nolimits_n T_{O_K^{(n)}}(\Lambda^{(n)})$$
be any element which corresponds to $U^\bullet \in T_{O_K}(\Lambda).$
By (the proof in) Proposition \ref{Tuben}, we have inclusions $Y_{U^{n+1}_\bullet} \subset Y_{U^n_\bullet}$ and therefore  homomorphisms
$$H^\ast_{Y_{U^{n+1}_\bullet}}(\P^d_K, \cF) \rightarrow H^\ast_{Y_{U^n_\bullet}}(\P^d_K, \cF)$$ for all $ n\in \N.$
In particular, we get a projective system
$(H^\ast_{Y_{U^n_\bullet}}(\P^d_K, \cF))_{n\in \N}$ of $K$-vector spaces. As in Lemma \ref{Frechetspaces} one verifies that this system
consists of $K$-Fr\'echet spaces and the transition maps are dense. Let $U^f_\bullet$ the largest subobject in $U^\bullet$
such that $U^f_\bullet \in T^f_{O_K}(\Lambda)= T_{K}(V).$

\begin{Proposition}\label{letzteProp}
There is a topological isomorphism of $K$-Fr\'echet spaces
\begin{equation*}\varprojlim_{n \in \N} H^\ast_{Y_{U^n_\bullet}}(\P^d_K, \cF) \cong  H^\ast_{Y_{U^f_\bullet \otimes K}}(\P^d_K, \cF).
\end{equation*}
\end{Proposition}

\proof This follows from the description  of the analytic varieties $Y_{U_\bullet} $ given in the proof of Proposition \ref{Tuben} together with  Proposition \ref{projlimI}   \qed

We consider the  acyclic complex of Theorem \ref{Theorem4} in the case $\cG=\Z.$ Then the resulting complex induces for each $n \in \mathbb N,$ a
spectral sequence
$$E_1^{-p,q,n}=  \bigoplus_{U_\bullet \in T(\Lambda^{(n)})^p} H^q_{Y_{U_\bullet}}(\P^d_K,\cF) \Longrightarrow H^{-p+q}_{{\cY}_n}(\P^d_K,\cF).$$
This spectral sequence has a similar shape as the one in section \ref{Berechnung}. Its rows are given by
$$E_1^{\bullet,j,n}: \bigoplus_{U_\bullet \in T(\Lambda^{(n)})^{n\cdot j} \atop {\rm rk}(U_1)=d-j} H_{Y_{U_\bullet}}^{j}(\P^d_K, \cF)\rightarrow \!\!\!\!\!\!\!\  \bigoplus_{U_\bullet \in T(\Lambda^{(n)})^{n\cdot j-1} \atop {\rm rk}(U_1)=d-j} H_{Y_{U_\bullet}}^{j}(\P^d_K, \cF)$$
$$\rightarrow \!\!\!\!\!\!\! \bigoplus_{U_\bullet \in T(\Lambda^{(n)})^{n\cdot j-2} \atop {\rm rk}(U_1)=d-j} H_{Y_{U_\bullet}}^{j}(\P^d_K, \cF) \rightarrow \ldots  \rightarrow \bigoplus_{U_\bullet \in T(\Lambda^{(n)})^1 \atop {\rm rk}(U_1)=d-j} H_{Y_{U_\bullet}}^{j}(\P^d_K, \cF)$$

\noindent Passing to the limit and applying Proposition \ref{letzteProp} resp. Proposition \ref{projlimI} to it yields a spectral sequence
$$E_1^{-p,q} \Longrightarrow H_{\cY}^{-p+q}(\P^d_K, \cF).$$ 
Its rows are given by
$$E_1^{\bullet,j}: \varprojlim_{n\in \N}\bigoplus_{U_\bullet \in T(\Lambda)^{n\cdot j} \atop {\rm rk}(U_1)=d-j} H_{Y_{U^f_\bullet\otimes K}}^{j}(\P^d_K, \cF)\rightarrow  \varprojlim_{n\in \N} \bigoplus_{U_\bullet \in T(\Lambda)^{n\cdot j-1} \atop {\rm rk}(U_1)=d-j} H_{Y_{U^f_\bullet \otimes K}}^{j}(\P^d_K, \cF)$$
$$\rightarrow   \varprojlim_{n\in \N} \bigoplus_{U_\bullet \in T(\Lambda)^{n\cdot j-2} \atop {\rm rk}(U_1)=d-j} H_{Y_{U^f_\bullet\otimes K}}^{j}(\P^d_K, \cF) \rightarrow \ldots  \rightarrow \varprojlim_{n\in \N} \bigoplus_{U_\bullet \in T(\Lambda)^1 \atop {\rm rk}(U_1)=d-j} H_{Y_{U^f_\bullet\otimes K}}^{j}(\P^d_K, \cF).$$

\noindent Now we apply Proposition \ref{Homotopieaequiv}. We see that the spectral sequence (\ref{ss}) of section \ref{Berechnung} is homotopy equivalent
to $E_1^{\bullet,\bullet}.$ Thus  from now on, we may carry on with the computation there. \qed

\end{document}